\documentclass[english,11pt]{article}
\usepackage{float}
\usepackage{algorithm,algorithmic}
\usepackage{amsmath}
\usepackage{amssymb}
\usepackage{graphicx}
\usepackage{bm}
\usepackage{amsthm}
\usepackage{subcaption}
\usepackage{enumitem}
\usepackage[title]{appendix}
\usepackage{babel}
\usepackage{stfloats}
\usepackage{wrapfig}

\usepackage[utf8]{inputenc} 
\usepackage[T1]{fontenc}    
\usepackage{url}            
\usepackage{booktabs}       
\usepackage{amsfonts}       
\usepackage{nicefrac}       
\usepackage{microtype}  
\usepackage{mathrsfs}

\usepackage{graphicx}
\usepackage{geometry}
\newgeometry{
	top= 1in,
	bottom =  1.4 in,
	headheight= 5pt,
	headsep=25pt,
	footskip=30pt,
	lmargin=1  in,  rmargin= 1.5in
}

\numberwithin{equation}{section}
\usepackage{hyperref}
\hypersetup{
	colorlinks=true,
	linkcolor=red,
	filecolor=magenta,      
	citecolor=blue,
}
\newtheorem{lem}{Lemma}
\newtheorem{prop}{Proposition}

\newtheorem{defi}{Definition}

\numberwithin{theo}{section}
\numberwithin{lem}{section}
\numberwithin{prop}{section}
\numberwithin{alg}{section}
\numberwithin{assum}{section}
\numberwithin{defi}{section}
\numberwithin{coro}{section}
\numberwithin{claim}{section}
\numberwithin{example}{section}

\theoremstyle{remark}
\newtheorem{remark}{Remarks}
\numberwithin{remark}{section}

\usepackage{mathtools}
\def\Tiny{\fontsize{4pt}{4pt}\selectfont}
\newcommand*{\eqdef}{\ensuremath{\overset{\mathclap{\text{\Tiny def}}}{=}}}

\title{General Optimal Step-size  for ADMM-type Algorithms: Domain Parametrization   and  Optimal   Rates} 
\author{Yifan Ran}
\date{}
\begin{document}
	\maketitle
	\makeatletter{\renewcommand*{\@makefnmark}{}
		\footnotetext{
			The	author is with the department of Electrical and Electronics Engineering, Imperial College London, UK (e-mail: y.ran18@imperial.ac.uk).  \makeatother}}

\maketitle

\begin{abstract}
	In this work,  we  solve a  49-year open  problem, the general optimal step-size for  ADMM-type algorithms.
	For a  convex program: $\text{min.} \,\, f(\bm{x}) + g(\bm{z})$,  $\text{s.t.}\, \mathcal{A}\bm{x} - \mathcal{B}\bm{z} = \bm{c}  $, given an arbitrary fixed-point  initialization 
	$  \bm{\zeta}^0 $, an optimal step-size choice is given by a root of the following polynomial:
	\begin{equation*} 
	\rho^4\Vert \mathcal{A}\bm{x}^\star\Vert^2  - \rho^3\langle \mathcal{A}\bm{x}^\star, \bm{\zeta}^0\rangle + \rho\langle \bm{\lambda}^\star,\bm{\zeta}^0\rangle -  \Vert\bm{\lambda}^\star\Vert^2 = 0,
	\end{equation*}
	with $ \rho \neq 0 $ a domain  step-size,
	which  relates to  the classical positive one  via $ \gamma = \rho^2$.   	
	We denote  by   $ \cdot^\star $  the   optimal solution, by  $ \bm{\lambda} $ the Lagrange  multiplier associated with the  equality  constraint (dual variable).
	
	The above  polynomial always admits a closed-form solution.
	The optimality  is in the sense  that a worst-case fixed-point convergence rate  is minimized,  which  is a balance  of   the normalized  primal  and  dual   iterates convergence  speed  (reciprocally related). 
	In cases where  either the  primal or dual solution  is trivial (a zero vector),  improvement  can  be made by  accelerating    the non-trivial sequence only.
	
	For practical use,   adaptively replace the above optimal solutions with the current iterates, which are known at every iteration. Numerically, it  exhibits  almost identical performance as the  theoretical one  (after a few iterations),  similar to the underlying best fixed step-size (found by exhaustive grid  search).

\vspace{11pt}
\noindent 
$ \textbf{Keywords:} $\,  
Proximal operator, Duality,	Fixed-point theory, Bounded linear operators,	Range and domain parametrizations,  Primal-dual  solution angle, Initialization angle, Firm non-expansiveness,	Alternating  direction method of multipliers (ADMM).

\vspace{11pt}

\end{abstract}

\section{Introduction}	
An ongoing trend  referred to as `Big Data' is drawing an increasing amount of attention  from both industry and academia. How to efficiently handle the explosion of data is  of central importance. In the review paper by Stephen Boyd et al. \cite{boyd12}, the  
Alternating  Direction Method of Multipliers (ADMM) is recommended: `\textit{It is  well suited to distributed convex optimization, and in particular to large-scale problems arising in statistics, machine learning, and related areas.'} 
The ADMM algorithm was  introduced by Glowinski and Marrocco \cite{glowinski1975approximation} in 1975, and Gabay and Mercier \cite{gabay1976dual} in 1976, with roots in the 1950s.
The algorithm was studied throughout the 1980s, and by the mid-1990s, most  of the theoretical results were established \cite{boyd12}. 

One  theoretical  consequence is that ADMM will converge  with  arbitrary positive step-size choices.
Interestingly, all these  choices, along with their corresponding iteration number complexity, seem   to  roughly form a  convex graph, implying  the feasibility of determining the  best choice by convex  optimization techniques \cite{rockafellar1997convex,boyd2004convex}.   
However,  how  to formulate  the general  convex program and then solve for  its  optimal solution has  been a  mystery.
Challengingly, numerical experiments show that the underlying best fixed  step-size  can take any positive values, from extremely small to arbitrarily large, related to problem structures and  data sizes. 
Particularly, the  best  choice does not generally translate to alternative  formulations or data sets. 
It is   hard  to  see   any pattern  and  hence difficult to guess a  good choice. 
Meanwhile, a  good  step-size can largely improve  the algorithm's   efficiency.
It is therefore of  both theoretical  and  practical  interest to determine  the  optimal step-size choice.

Unfortunately,  the  ADMM optimal step-size  is one of the very few issues that have not been addressed in the last 49 years. As pointed out in a recent paper by Bartolomeo Stellato and Stephen Boyd  et al. \cite[sec. 1.2, First-order methods]{stellato2020osqp} `\textit{Despite some recent theoretical
	results \cite{banjac2018tight,giselsson2016linear}, it remains unclear how to select those parameters to optimize the algorithm convergence rate}'.
The same open  problem is also mentioned  by  Ernest K. Ryu and Wotao Yin in \cite[Sec. 8, Parameter selection]{ryu2022large}.

In this work,  we provide a general, worst-case   optimal choice, without  the  need for any  tailored structure information.  It is simple, universally applicable, highly efficient,  and practically useful. 
The key is to introduce the \textit{domain parametrization}, i.e.,   associating the step-size parameter  to  the domain of a  function.
This is  in  contrast to the  current literature  where one   assigns the parameter to  the range, i.e.,   the left-hand side.
We argue that the  classical  range association is not  mathematically natural, in the sense that  an extra, implicit, step-size-related  scaling  arises, which  shifts  the algorithm sequence  to a  parallel variant.
Such  parallelism does not affect  any practical uses, since all convergence  properties (such as rates, step-size choices) will be  the same. However,
it would cause  a  strong   barrier to finding  the  general theoretical optimal step-size.
While this issue is not  obvious, there does exist a hint
from the well-known Moreau identity --- its original symmetry  is lost after the classical parametrization.
Such  extra scaling/shifting will be naturally avoided by the \textit{domain parametrization}.
In fact, there are  more benefits to our new parametrization. We will present  some new, simple, highly  unified analyses for the ADMM-type algorithms,  and draw several  useful conclusions.  

After finding the  optimal step-size, we can invoke it to obtain some  optimal convergence rates, which reveal  that (i) under zero initialization, the ADMM efficiency is intrinsically characterized by the angle between the primal and  dual solutions; 
(ii)  if  either  the primal or dual solution is trivial  (a zero vector),  then under zero initialization, the non-trivial iterate convergence  speed is independent of the step-size  choice;
(iii) meanwhile, in the same partly-trivial  setting,  if instead employ  an arbitrary non-zero  initialization. 
Then, invoke a corresponding optimal step-size, we show  that the non-trivial iterate convergence rate will be improved by a factor of $ \sin^2 \omega \in [0,1]$, where $ \omega $ is the angle between  the initialization and the  non-trivial solution.  
Particularly, when $ \omega = 0 $ or $ \omega = \pi $, we obtain  $ \sin^2 \omega  =  0$, and the  algorithm   converges instantly, implying  the  initialization  is  optimal.
This last case also implies that only the initial angle matters, due to the distance will be automatically adjusted by a closed-form optimal step-size.

Our optimal step-size  changes with different  initializations, and a  good initialization choice can  be shown  to improve the algorithm efficiency.
In the literature, how to estimate a good initialization  has   been intensively studied, referred to as the `warm start'.  
Typical  success is regarding the type of problems that sequentially apply ADMM. 
The reason for the repeat application of ADMM can be due to decisions need updating as the scenario changes dynamically.  
In this case, one can directly adopt the solution from the previous ADMM as an initialization for the next one, see e.g. \cite{chen2020online}. 
Another class of problems is clustering and dictionary learning. Typically, there is an outer and an inner loop. One can directly use the output of each inner loop as a warm start for the next one, see e.g. \cite{ye2014scaling}.  
At last, one  can  apply machine learning to find a near-optimal initialization in a general sense, typical work can be found by Bartolomeo Stellato et. al. \cite{sambharya2023end,sambharya2023learning}, see also  \cite{mak2023learning} on a power flow problem.

At last,  our results are transferable owing  to the  equivalence  of ADMM  and many other algorithms. 
The most well-known equivalence is perhaps with the Douglas-Rachford splitting (DRS) \cite{douglas1956numerical} from numerical analysis. It is widely recognized that ADMM iterates are equivalent to applying DRS to its dual problem, see e.g. \cite{eckstein1992douglas,fortin2000augmented}. Recently,  Daniel O'Connor and Lieven Vandenberghe \cite{o2020equivalence} show that the   Primal-Dual Hybrid Gradient (PDHG) method \cite{chambolle2011first,esser2010general,pock2009algorithm} is  also equivalent to ADMM. Apart from being  equivalent, there are many closely related algorithms such as Spingarn’s method of partial inverses, Dykstra’s alternating projections method, and Bregman iterative algorithms for $ l_1 $ problems in signal	processing,  see more details from the review paper \cite{boyd12}.

\subsection{Notations}	
We use blackboard bold letters  such  as  	$ \mathbb H $  to denote a real Hilbert space  with  inner product $ \langle\cdot, \cdot\rangle $ and norm  $\Vert\cdot\Vert  $ equipped.
The uppercase calligraphic letters  are used to denote  bounded linear operators.  Particularly, $ \mathcal{I} $ always denotes the  identity operator. 
An operator $ \mathcal{A} $ on a Hilbert space is a point-to-set mapping $ \mathcal{A}:\mathbb H \rightarrow 2^\mathbb H $.
We denote by $  \Gamma_0 (\mathbb H) $  the  space  of  convex, closed  and proper  (CCP) functions from $ \mathbb H $  to the extended real line  $ (-\infty, +\infty]$,
by $ \mathfrak{B}(\mathbb{H}, \mathbb{K}) $ the space of bounded linear operators from domain $ \mathbb{H} $ to range $ \mathbb{K} $,
by $\circ$ the operator or function composition, by $ \text{dom}(\cdot) $ the  domain, by  $ \text{ran}(\cdot) $ the  range, by  $ \emptyset $  the  empty set.
The  adjoint of an  operator $ \mathcal{A}  \in  \mathfrak{B}(\mathbb{H}, \mathbb{K}) $  is  the  unique operator  $ \mathcal{A}^*  \in  \mathfrak{B}( \mathbb{K},  \mathbb{H}) $
satisfying $ \langle \mathcal{A} \bm{x}, \bm{y} \rangle = \langle  \bm{x}, \mathcal{A}^*\bm{y} \rangle , \,  \forall \bm{x} \in \mathbb{H}, \forall \bm{y} \in \mathbb{K}$. 
The Fenchel  conjugate  function is defined as $ f^*(\cdot) = \sup\, \langle \bm{z}, \cdot \rangle - f(\bm{z})$.
The uppercase bold, lowercase bold, and not bold letters are used for matrices, vectors, and scalars, respectively.

\subsection{ADMM algorithm}
ADMM solves the following general convex problem:
\begin{align}\label{pro00}
\underset{\bm{x},\bm{z}}{\text{minimize}}& \quad   f(\bm{x}) + g(\bm{z}), \nonumber\\ 
\text{subject\,to}& \quad \mathcal{A}\bm{x} - \mathcal{B}\bm{z} = \bm{c} ,
\end{align}
with variables $ \bm{x} \in  \mathbb H $, $ \bm{z} \in  \mathbb P $, $ \bm{c} \in  \mathbb K$ and functions  $ f\in \Gamma_0 (\mathbb{H}) $, $ g\in \Gamma_0 (\mathbb{P}) $, where $ \mathcal{A} \in  \mathfrak B (\mathbb H,  \mathbb K)$ and  $ \mathcal{B} \in  \mathfrak B (\mathbb P,  \mathbb K)$ are injective.
A  solution is assumed exists.

The augmented Lagrangian is
\begin{equation}\label{lag_class}
\mathcal{L}_\gamma(\bm{x},\bm{z},\bm{\lambda}) 
\,=\,  f(\bm{x}) + g (\bm{z})  + \frac{\gamma}{2}\Vert\mathcal{A}\bm{x} - \mathcal{B}\bm{z} - \bm{c} + \bm{\lambda}/\gamma\Vert^2, 
\end{equation}	
with $\gamma > 0$ being a positive step-size.	 The  classical ADMM iterates are 
\begin{align}\label{admm_iterates}
\bm{x}^{k+1} =\,\,& \underset{\bm{x}}{\text{argmin}} \,\, 	\mathcal{L}_\gamma(\bm{x},\bm{z}^k,\bm{\lambda}^k) , \tag{primal}\\				
\bm{z}^{k+1} =\,\, & \underset{\bm{z}}{\text{argmin}} \,\, 	\mathcal{L}_\gamma(\bm{x}^{k+1},\bm{z},\bm{\lambda}^k) ,\tag{auxiliary}\\
\bm{\lambda}^{k+1} =\,\,  &\bm{\lambda}^{k} + \gamma( \mathcal{A}\bm{x}^{k+1} - \mathcal{B}\bm{z}^{k+1} - \bm{c}), 	\tag{dual}
\end{align}	
where $ k = 0,1,\dots $  denotes the iteration number counter.

Additionally, the  above 3-steps algorithm can be rewritten into  a 1-step  fixed-point  characterization:
\begin{equation}
\bm{\zeta}^{k+1}_1 =  \mathcal{F}_\gamma  \bm{\zeta}^k_1 ,
\end{equation}
where $ \mathcal{F}_\gamma  $  is a certain $\gamma$-related operator detailed later, and where
$ \bm{\zeta}_1 $  denotes a fixed-point variable,  with a classical expression associated  with \eqref{pro00} being
\begin{equation}
\hspace{3cm}\bm{\zeta}^{k+1}_1 = \mathcal{A}\bm{x}^{k+1} + \bm{\lambda}^{k}/\gamma. \tag{classical fixed-point}
\end{equation}
Moreover, the above $  \mathcal{F}_\gamma  $ is well-known to be firmly  non-expansive (defined below), see e.g.  \cite{ryu2022large,poon2019trajectory}.
\begin{defi}\label{def01}
	Let $ \mathcal{F}: \mathbb H \rightarrow \mathbb H$. Then,  $ \mathcal{F} $ is firmly non-expansive if 
	\begin{equation}
	\Vert \mathcal{F}\bm{x} - \mathcal{F}\bm{y}  \Vert^2 \leq  \langle \mathcal{F}\bm{x} - \mathcal{F}\bm{y}, \bm{x} - \bm{y} \rangle, \quad\forall \bm{x},\bm{y} \in \mathbb H.
	\end{equation} 
\end{defi}

\subsection{Key results}

\subsubsection{open problem (theoretical)}
The general optimal step-size  remains open in the last 49 years for ADMM  (and 68  years for DRS), see \cite[Sec. 1.2, First-order methods]{stellato2020osqp}, \cite[Sec. 8, Parameter selection]{ryu2022large}.
In this work, we provide a  universally applicable optimal choice as a closed-form root of the following polynomial:
\begin{equation} 
\rho^4\Vert \mathcal{A}\bm{x}^\star\Vert^2  - \rho^3\langle \mathcal{A}\bm{x}^\star, \bm{\zeta}^0\rangle + \rho\langle \bm{\lambda}^\star,\bm{\zeta}^0\rangle -  \Vert\bm{\lambda}^\star\Vert^2 = 0,
\end{equation}
with $\rho\neq 0$  denotes a domain step-size,  which  relates  to the  classical one as  in  \eqref{lag_class}  via  $\gamma  = \rho^2$, where $  \bm{\zeta}^0  = \rho_0\mathcal{A}\bm{x}^0 + \bm{\lambda}^0/\rho_0$ is  an arbitrarily  fixed-point initialization, and  where $ \bm{x}^0 $, $ \bm{\lambda}^0 $,  and  $\rho_0$  denote  the primal  iterate,  dual iterate, and  the step-size  initial  choice, respectively.

The  above step-size selection does not incorporate any tailored structure information, and  is therefore  only optimal in a worst-case  sense.  Nevertheless, numerically,   we observe that it  is  already   close to the  underlying best fixed choice  in all our experiments (found by exhaustive grid  search), seemly sufficient for basic practical usages.

\subsubsection{practical  use}
For practical use, solve the following polynomial at every $ (k+1) $-th iteration: 
\begin{equation}
\rho^4\Vert \mathcal{A}\bm{x}^{k+1}\Vert^2  - \rho^3\langle \mathcal{A}\bm{x}^{k+1}, \bm{\zeta}^0  \rangle + \rho\langle \bm{\lambda}^{k+1}, \bm{\zeta}^0  \rangle -  \Vert\bm{\lambda}^{k+1}\Vert^2 = 0,  
\end{equation}
see more  details  from Sec. \ref{sec_prac}. Numerically,  it  exhibits highly similar performance  as the above theoretical one  (almost  identical after a few iterations, typically $ k = 10 $).


For ease  of  use, below  we present the  simplest  version, expressed in terms of  the  classical step-size $\gamma  >  0$.
\begin{algorithm}[H]
	\caption{ ADMM  with  adaptive step-sizes (Simplified version)}
	\begin{algorithmic}[1]
		\REQUIRE  Set $\,\bm{z}^0= \bm{0}, \,\bm{\lambda}^0 = \bm{0},\,   \gamma_{0}  = 1 $.
		\WHILE{iteration $ k  =  0, 1, 2, \dots   $} \STATE
		\begin{align}
		\bm{x}^{k+1} =\,\,&  \, \underset{\bm{x}}{\text{argmin}} \,\, f(\bm{x})+ \frac{ \gamma_k}{2}\Vert  \mathcal{A}\bm{x} -   \mathcal{B}\bm{z}^k  - \bm{c} + \bm{\lambda}^k/ \gamma_k  \Vert^2, \\				
		\bm{z}^{k+1} =\,\, &\, \underset{\bm{z}}{\text{argmin}} \,\, g(\bm{z})+ \frac{ \gamma_k}{2}\Vert  \mathcal{A}\bm{x}^{k+1} -   \mathcal{B}\bm{z}  - \bm{c}  + \bm{\lambda}^k/ \gamma_k  \Vert^2, \\
		\bm{\lambda}^{k+1} =\,\,  &\,\bm{\lambda}^{k} +  \gamma_k \big( \mathcal{A}\bm{x}^{k+1} - \mathcal{B}\bm{z}^{k+1} - \bm{c} \big), \\
		\gamma_{k+1} 	=\,\,  &\,      \Vert \bm{\lambda}^{k+1} \Vert / \Vert \mathcal{A}\bm{x}^{k+1} \Vert.
		\end{align}		
		\ENDWHILE
		\ENSURE primal  and dual solutions $ \bm x^\star $  and  $ \bm \lambda^\star $, respectively.
	\end{algorithmic}
\end{algorithm}

\subsubsection{unscaled fixed-point}
In  the  literate,  there exist  two fixed-point expressions  for ADMM-type algorithms:
\begin{equation}\label{class_fix}
\underbracket{\bm{\zeta}_1^\star = \mathcal{A}\bm{x}^\star +\bm{\lambda}^\star/\gamma}_{\text{for primal  problem}}, 
\qquad \,\,
\underbracket{\bm{\zeta}_2^\star  = \gamma \mathcal{A}\bm{x}^\star + \bm{\lambda}^\star}_{\text{for dual  problem}}.
\end{equation}
see $ \bm{\zeta}_1^\star $  from  e.g. \cite[sec. 3.1]{ryu2022large}, corresponding to the primal problem, and $ \bm{\zeta}_2^\star $  from e.g. \cite[sec. 3.2]{ryu2022large}, \cite{poon2019trajectory}, corresponding to the Fenchel dual problem.

In  this work,  we   show  that the above expressions  lead to self-contradicted step-size choices,  see Sec. \ref{sec_fail}.   
This is due to  that they  are both implicitly scaled/shifted, and  the unscaled version  is
\begin{equation}
\bm{\zeta}^\star = \rho\mathcal{A}\bm{x}^\star +\bm{\lambda}^\star/\rho,  \tag{unscaled}
\end{equation}
where  $\rho  = \pm \sqrt{\gamma}$ is our domain  step-size, see derivation from  Proposition \ref{fpi}, and proof of being unscaled in Sec. \ref{self_d} via self-duality.

\subsubsection{Some  optimal rates}
Below, we   present some   convergence rates,  induced by certain optimal step-sizes. 

$\bullet$ First, a general case  with  zero initialization, see  theoretical details  from   Sec. \ref{sec_itr_conv1} and  numerical verifications  through Lasso in  Sec. \ref{num_lasso}:

(i) (Balanced  rates.) Suppose   $ \bm{x}^\star \neq  \bm{0} $,  $ \bm{\lambda}^\star \neq  \bm{0} $, i.e.,
$ \Vert \mathcal{A}\bm{x}^\star\Vert \neq  0, \, \, \Vert\bm{\lambda}^\star\Vert  \neq 0$ ($ \mathcal{A} $ assumed  injective).
Set  initialization $ \bm{\zeta}^0  =  \bm{0} $.  Let  
\begin{equation}
\hspace{3cm}    \theta  =  \arccos  \,  \frac{ \langle \mathcal{A}\bm{x}^\star, \bm{\lambda}^\star\rangle }{\Vert\mathcal{A}\bm{x}^\star \Vert \Vert \bm{\lambda}^\star \Vert}.  \tag{intrinsic   angle}
\end{equation}
Choose domain step-size  to be
$   \rho^\star = \pm  \sqrt{\Vert\bm{\lambda}^\star\Vert/\Vert\mathcal{A}\bm{x}^\star\Vert}$.
Then, 
\begin{align}
\frac{\Vert \mathcal{A}\bm{x}^{k+2} - \mathcal{A}\bm{x}^{k+1} \Vert^2   }{\Vert\mathcal{A}\bm{x}^\star \Vert^2 }
\,\,\leq\,\,  &  \,\,\frac{2}{k+1}   \,   (1  + \cos \theta ),   \tag{primal  }\\
\frac{\Vert \bm{\lambda}^{k+1} - \bm{\lambda}^{k} \Vert^2   }{\Vert\bm{\lambda}^\star \Vert^2 }
\quad\,\, \leq\,\, &  \,\,\frac{2}{k+1}   \,  (1  + \cos \theta ),   \tag{dual  }\\
\frac{\Vert \bm{\zeta}^{k+1} - \bm{\zeta}^{k} \Vert^2   }{\Vert\mathcal{A}\bm{x}^\star \Vert\Vert\bm{\lambda}^\star \Vert }
\quad\,\, \leq\,\, &  \,\,\frac{2}{k+1}   \,  (1  + \cos \theta ).  \tag{fixed-point} 
\end{align}

$\bullet$  Next, a special case  where the dual solution is  a  zero vector,  see   theoretical details  from   Sec. \ref{sec_special} and   numerical verifications  through a feasibility problem in  Sec. \ref{feasi}: 

(ii) (Fixed non-trivial rate.) 
Suppose  $ \Vert \mathcal{A}\bm{x}^\star\Vert \neq  0,  \, \Vert\bm{\lambda}^\star\Vert  = 0 $. 
Set initialization $ \bm{\zeta}^0  =  \bm{0} $.  Then,
\begin{align}
\qquad	{\Vert \mathcal{A}\bm{x}^{k+2} - \mathcal{A}\bm{x}^{k+1} \Vert^2  }
\,\,\leq\,\,\,\,    & \frac{1}{k+1}\cdot{ \Vert  \mathcal{A}\bm{x}^\star \Vert^2 } , \tag{non-trivial primal}  \\
{\Vert \bm{\lambda}^{k+1} - \bm{\lambda}^{k} \Vert^2   }
\,\, \leq  \,\,\,\,  &  \frac{1}{k+1}  \cdot \Vert\mathcal{A}\bm{x}^\star\Vert^2 \cdot  \rho^4    ,  \tag{trivial dual}  \\
\Vert \bm{\zeta}^{k+1} - \bm{\zeta}^k  \Vert^2 
\,\, \leq  \,\,\,\,   &  \frac{1}{k+1} \cdot \Vert\mathcal{A}\bm{x}^\star\Vert^2  \cdot  \rho^2,  \tag{fixed-point} 
\end{align}
where the  optimal step-size is  $ \rho^\star  \rightarrow  0 $, which  only affect the dual and  fixed-point  sequences.

(iii) (Guaranteed non-trivial improvement.)
Suppose  $ \Vert \mathcal{A}\bm{x}^\star\Vert \neq  0$,  $\Vert\bm{\lambda}^\star\Vert  = 0$.  Set  any  initialization $ \bm{\zeta}^0  \neq \bm{0} $.  
Let 
\begin{equation}
\hspace{2.5cm}   \omega_1  =  \arccos  \, \frac{\langle \mathcal{A}\bm{x}^\star, \bm{\zeta}^0  \rangle }{\Vert\mathcal{A}\bm{x}^\star\Vert\Vert \bm{\zeta}^0\Vert}. \tag{initial  angle}
\end{equation}
Choose domain step-size  to be $ \rho^\star_\text{pri} =  {\Vert\bm{\zeta}^0 \Vert}/  (\Vert\mathcal{A}\bm{x}^\star\Vert\cos \omega_1 )$.
Then,
\begin{align}
{\Vert \mathcal{A}\bm{x}^{k+2} - \mathcal{A}\bm{x}^{k+1} \Vert^2  }
\,\, \leq  \,\,\,\,  &     \frac{1}{k+1}  \cdot \,\, \Vert  \mathcal{A}\bm{x}^\star \Vert^2\cdot \sin^2 \omega_1 ,  \tag{non-trivial primal}   \\
{\Vert \bm{\lambda}^{k+1} - \bm{\lambda}^{k} \Vert^2   }
\,\, \leq  \,\,\,\,  &  \frac{1}{k+1} \cdot  \frac{ \Vert  \bm{\zeta}^0\Vert^4}{  \Vert \mathcal{A}\bm{x}^\star \Vert^2} \cdot \frac{ \tan^2  \omega_1 }{\cos^2 \omega_1},  \tag{trivial dual} \\
\Vert \bm{\zeta}^{k+1} - \bm{\zeta}^k  \Vert^2 
\,\, \leq  \,\,\,\,  & \frac{1}{k+1} \cdot  \,\,\Vert  \bm{\zeta}^0\Vert^2  \cdot \tan^2  \omega_1 .  \tag{fixed-point} 
\end{align}
Similar results hold for a symmetric situation  where $ \Vert \mathcal{A}\bm{x}^\star\Vert =  0,  \, \Vert\bm{\lambda}^\star\Vert  \neq 0$.

\subsection{Organization}
In  Sec. \ref{sec_2}, we introduce a new domain-type  parametrization and several useful consequences on the  proximal operator and  the ADMM  algorithm. 
Sec.  \ref{sec_3} establishes some convergence rates.  Sec. \ref{sec_4} presents the optimal step-sizes by minimizing  worst-case convergence rates.
Sec. \ref{sec_5} conducts two numerical experiments to verify our results.

\section{Two  parametrizations: range and domain types}\label{sec_2} 
In this  section, we introduce  a new domain-type  parametrization, parallel to  the  classical  range-type parametrization.  
We show that the  classical way  introduces  an  implicit,  extra,  and  step-size-related scaling, causing  the 49-year optimal step-size  open problem.
This issue will be naturally avoided by  our new  parametrization. 

\subsection{Fundamental  ground}
Throughout the paper, we will repeatedly use   a  standard  tool,  the proximal operator, which  is  widely  used  across  fields, see   a survey paper  \cite{proxi_algs}.
It is developed as  a  generalization of the projector (onto convex sets) by Moreau  \cite{moreau1965proximite}, defined as (without parametrization):
\begin{defi}
	Given $ f \in \Gamma_0 (\mathbb H) $,   $\text{dom}(f) \neq \emptyset$, the  proximal operator is defined as
	\begin{equation}
	\textbf{Prox}_{ f}(\bm{v}) =  \underset{\bm{z}}{\text{argmin}} \,\,   f(\bm{z}) + \frac{1}{2}\Vert \bm{z} - \bm{v}\Vert^2,
	\end{equation}
	where $ \Vert \cdot \Vert $ denotes the Euclidean norm, and  where $ \text{dom}(\cdot) $  denotes the domain.
\end{defi}	
An important  relation is the following Moreau identity:
\begin{equation}\label{moreu_bas}
\bm{v} =\,\textbf{Prox}_{f} (\bm{v}) +  \textbf{Prox}_{f^*} (\bm{v}),
\end{equation}
where   $ f^*(\cdot) \eqdef \sup\, \langle \bm{z}, \cdot \rangle - f(\bm{z})$  is  the Fenchel  conjugate  function. 
The above is an important link for primal and  dual problems.

\subsection{Classical range parametrization}
It is standard to consider a parametrized version of the   proximal operator.
\begin{defi}\label{prox_lam}
	Given $ f \in \Gamma_0 (\mathbb H) $,   $\text{dom}(f) \neq \emptyset$. Let    $\gamma > 0 $  denote a  positive step-size.  Then, the $ \gamma $-parametrized proximal operator is 
	\begin{equation}\label{class_prox}
	\textbf{Prox}_{\gamma f}(\bm{v}) =  \underset{\bm{z}}{\text{argmin}} \,\,  \gamma f(\bm{z}) + \frac{1}{2}\Vert \bm{z} - \bm{v}\Vert^2
	= \underset{\bm{z}}{\text{argmin}} \,\,   f(\bm{z}) + \frac{1}{2\gamma}\Vert \bm{z} - \bm{v}\Vert^2.
	\end{equation}
\end{defi}
The Moreau identity becomes
\begin{align}\label{moreau}
\bm{v} 
\,=\,  \textbf{Prox}_{\gamma f}(\bm{v}) + \bigg(\gamma\mathcal{I} \bigg)\circ \textbf{Prox}_{\frac{1}{\gamma}f^*} \circ \bigg(\frac{1}{\gamma} \mathcal{I} \bigg)  (  \bm{v} )  ,
\end{align}
where   $ \mathcal{I}: \mathbb H \rightarrow  \mathbb H $  denotes the identity  operator.
\begin{remark}[hint:  extra  scaling]
	Compare \eqref{moreau} and  \eqref{moreu_bas}, we  see that --- after parametrization, the Moreau identity components are  no longer symmetric, with an  extra $ \gamma  \, (\cdot)\,    \frac{1}{\gamma} $  type structure.
	This seemly inconsequential  phenomenon contains an important  hint on   the step-size open  problem --- an extra scaling.
\end{remark}

\subsection{Domain  parametrization}\label{sec_domain}
We note that  conventionally one always  assigns the step-size  to the range of a function. 
Here, we  introduce a symmetric  domain-type parametrization.
\begin{prop}\label{new_def1}
	Given $ f \in \Gamma_0 (\mathbb H) $,   $\text{dom}(f) \neq \emptyset$ and  scalar  $\rho \neq 0 $, the domain-parametrized proximal operator is 
	\begin{align}\label{def_right}
	\textbf{Prox}_{f \rho} (\bm{v}) 
	=& \,\,\,\,\underset{\bm{z}_1\in  \mathbb{H}}{\text{argmin}} \,\, f(\rho\bm{z}_1 )+ \frac{1}{2}\Vert \bm{z}_1- \bm{v}\Vert^2, \nonumber\\
	=&\, \frac{1}{\rho}\,\underset{\bm{z}_2\in  \mathbb{H}}{\text{argmin}} \,\, f(\bm{z}_2)+ \frac{1}{2}\Vert \frac{\bm{z}_2}{\rho}- \bm{v}\Vert^2.
	\end{align}
\end{prop}
\begin{proof}
	The second equality is  obtained via a variable substitution
	$ \bm{z}_2 = \rho\bm{z}_1 $.
\end{proof}

\begin{lem}[translation rule]\label{lem_rel}
	The classical  proximal operator  as  in \eqref{class_prox}  and our  domain-parametrized version \eqref{def_right} are related via
	\begin{align}
	&\textbf{Prox}_{f \rho}(\bm{v})  = \quad\bigg(\frac{1}{\rho}\mathcal{I}\bigg)   \circ	\textbf{Prox}_{\rho^2 f}  \circ \bigg(\rho\,\mathcal{I}\bigg) \,\,(\bm{v}), \,\,\, \rho \neq  0  \label{re1}\\
	&\textbf{Prox}_{\gamma f}(\bm{v})  = \,\bigg(\pm\sqrt{\gamma}\,\mathcal{I}\bigg) \circ	
	\textbf{Prox}_{f  (\pm\sqrt{\gamma})}\circ\bigg(\pm\frac{1}{\sqrt{\gamma}}\mathcal{I}\bigg)\,(\bm{v}),  \,\,\,  \gamma  >  0.  \label{re2}
	\end{align}
\end{lem}
\begin{proof}
	For the first relation \eqref{re1}. From the right-hand  side,  we obtain
	\begin{align}
	\bigg(\frac{1}{\rho}\mathcal{I}\bigg)\circ	\textbf{Prox}_{\rho^2 f}  \circ  	\bigg(\rho\,\mathcal{I}\bigg)  \,\bm{v}
	=&   \frac{1}{\rho}  \,  \underset{\bm{z}}{\text{argmin}} \,\,  f(\bm{z}) + \frac{1}{2\rho^2}\Vert \bm{z} - \rho\bm{v}\Vert^2,\nonumber\\
	=&   \frac{1}{\rho}  \,  \underset{\bm{z}}{\text{argmin}} \,\,  f(\bm{z}) + \frac{1}{2}\Vert \frac{\bm{z}}{\rho} - \bm{v}\Vert^2, 
	\end{align}
	which coincides with definition \eqref{def_right}.

	For the second relation \eqref{re2}, consider the positive sign  case. From the right-hand  side, we obtain
	\begin{align}
	\bigg(\sqrt{\gamma}\,\mathcal{I}\bigg) \circ	\textbf{Prox}_{f  \sqrt{\gamma}}\circ\bigg(\frac{1}{\sqrt{\gamma}}\mathcal{I}\bigg)\,\bm{v}
	=&   \frac{\sqrt{\gamma}}{\sqrt{\gamma}}\,\underset{\bm{z}_2}{\text{argmin}} \,\, f(\bm{z}_2)+ \frac{1}{2}\Vert \frac{\bm{z}_2}{\sqrt{\gamma}}- \frac{\bm{v}}{\sqrt{\gamma}}\Vert^2  ,\nonumber\\
	=&   \quad\,\,\,\,   \underset{\bm{z}_2}{\text{argmin}} \,\, f(\bm{z}_2)+ \frac{1}{2\gamma}\Vert\bm{z}_2- \bm{v}\Vert^2 ,  
	\end{align}
	which coincides with definition \eqref{class_prox}.
	For the negative sign, similarly,
	\begin{align}
	&		\bigg(-\sqrt{\gamma}\,\mathcal{I}\bigg) \circ	\textbf{Prox}_{f  (-\sqrt{\gamma})}\circ\bigg(\frac{1}{-\sqrt{\gamma}}\mathcal{I}\bigg)\,\bm{v}   \nonumber\\
	&=  \frac{-\sqrt{\gamma}}{-\sqrt{\gamma}}\,\underset{\bm{z}_2}{\text{argmin}} \,\, f(\bm{z}_2)+ \frac{1}{2}\Vert \frac{\bm{z}_2}{-\sqrt{\gamma}}- \frac{\bm{v}}{-\sqrt{\gamma}}\Vert^2    
	=  \underset{\bm{z}_2}{\text{argmin}} \,\, f(\bm{z}_2)+ \frac{1}{2\gamma}\Vert\bm{z}_2- \bm{v}\Vert^2,
	\end{align}
	which is the same as definition \eqref{class_prox}.
	The proof is now concluded.
\end{proof}

\begin{remark}[step-size relation]
	In  view  of the above lemma, clearly our domain step-size relates to the classical range step-size via
	\begin{equation}
	\gamma  \,\,=\,\, \rho^2.
	\end{equation}
\end{remark}

\subsection{Extra scaling}
By now, we  repeatedly observe $ \rho\mathcal{I}\circ (\cdot) \circ \frac{1}{\rho}\mathcal{I} $  type of operations,
recall Lemma  \ref{lem_rel} and Moreau identity \eqref{moreau}. 
Here, we show that this operation implies  shifting. How to remove it is key to our technique.

\subsubsection{parallelism}
To start, we specify what we mean by `extra scaling'.

Let process $ \bm{\zeta}^{k+1} = \mathcal{F} \bm{\zeta}^k$ generate a sequence $ \{\bm{\zeta}^{k}\} $,  with    $ k = 0,1,\dots $ denotes  the iteration number  counter.  
Then, given any  $ \rho \neq 0  $, we have
\begin{align}
\,\,   \bm{\zeta}^{k+1}  =  \mathcal{F}  \, \bm{\zeta}^k,
\iff 
&\,\,\quad\, \bm{\zeta}^{k+1} \,=\, \frac{1}{\rho}\mathcal{I}\circ \bigg(\rho\mathcal{I}\circ \mathcal{F} \circ \frac{1}{\rho}\mathcal{I}\bigg)  \circ  \rho\mathcal{I}  \,\, (\bm{\zeta}^{k}), \nonumber\\
\iff 
&\,\,  (\rho \bm{\zeta}^{k+1}) \,=\,  \bigg(\rho\mathcal{I}\circ \mathcal{F} \circ \frac{1}{\rho}\mathcal{I}\bigg) \  ( \rho \bm{\zeta}^{k}).\label{rel00}
\end{align}

\begin{lem}\label{lem_01} 
	Let $ \mathcal{F}: \mathbb H \rightarrow \mathbb H$  be firmly non-expansive. Then,  the composited mapping $ \rho\mathcal{I}\circ \mathcal{F} \circ \frac{1}{\rho}\mathcal{I} $  with any $\rho\neq  0$ is  also  firmly  nonexpansive.
\end{lem}
\begin{proof}
	Recall the definition of  firm non-expansiveness, restated here  as
	\begin{defi}
		Let $ \mathcal{F}: \mathbb H \rightarrow \mathbb H$. Then,  $ \mathcal{F} $ is firmly non-expansive if 
		\begin{equation}
		\Vert \mathcal{F}\bm{x} - \mathcal{F}\bm{y}  \Vert^2 \leq  \langle \mathcal{F}\bm{x} - \mathcal{F}\bm{y}, \bm{x} - \bm{y} \rangle, \quad\forall \bm{x},\bm{y} \in \mathbb H.
		\end{equation} 
	\end{defi}	
	It follows that, 
	$  \forall  \bm{\zeta}, \bm y \in \mathbb{H} $, the following holds:
	\begin{align} 
	&\,\,\quad \Vert \mathcal{F} \bm{\zeta} -  \mathcal{F} \bm y\Vert^2 
	\,\leq \quad \langle  \mathcal{F}\bm{\zeta} - \mathcal{F} \bm y,  \bm{\zeta} -  \bm{y} \rangle, \nonumber\\
	\iff
	&\,\,\rho^2 \Vert \mathcal{F}  \bm{\zeta} -  \mathcal{F} \bm{y} \Vert^2 
	\,\leq \rho^2\langle \mathcal{F}\bm{\zeta} - \mathcal{F} \bm y, \bm{\zeta} - \bm{y} \rangle , \,\, \quad \rho \neq 0, \nonumber\\
	\iff
	&\,\,\quad \Vert \widetilde{\mathcal{F}}  \widetilde{\bm{\zeta}} -  \widetilde{\mathcal{F}}\widetilde{\bm{y}}\Vert^2 
	\,\leq \quad\langle \widetilde{\mathcal{F}}  \widetilde{\bm{\zeta}} - \widetilde{\mathcal{F}} \widetilde{\bm{y}},  \widetilde{\bm{\zeta}} - \widetilde{\bm{y}} \rangle,
	\end{align}
	with substitutions:
	\begin{equation}
	\widetilde{\bm{\zeta}} = \rho \bm{\zeta}, \qquad\widetilde{\bm{y}} = \rho \bm{y}, \qquad \widetilde{\mathcal{F}} = \rho\mathcal{I}\circ  \mathcal{F}\circ \frac{1}{\rho}\mathcal{I}.
	\end{equation}
	The proof is now concluded.
\end{proof}
\begin{lem}[$\rho$-distanced  parallel sequences]
	Let process $ \bm{\zeta}^{k+1} = \mathcal{F} \bm{\zeta}^k$ generate a sequence $ \{\bm{\zeta}^{k}\} $,  with  iteration  number counter  $ k= 0,1,2\dots $.

	Then,  there  exist the following $\rho$-distanced parallel variants:
	\begin{equation}
	\big\{\bm{y}^{k} \, |\, \bm{y}^{k}  = \rho\bm{\zeta}^{k}, \,\,\,  \forall \rho \in [-\infty, +\infty] / \{ 0 \} \big\} ,
	\end{equation}
	corresponding  to process  $ \bm{y}^{k+1} =  \widetilde{\mathcal{F}} \bm{y}^k$, 
	with  $ \widetilde{\mathcal{F}} = \rho\mathcal{I}\circ \mathcal{F} \circ \frac{1}{\rho}\mathcal{I}$.
\end{lem}
\begin{proof}
	The proof follows straightforwardly from \eqref{rel00} and Lemma \ref{lem_01}.
\end{proof}


\subsubsection{hidden step-size challenge}
Given the above $\rho$-distanced parallel  sequences,  clearly, all of  them share the  same  convergence properties, including  rates and the best step-size choice.
Then, we would determine  such a best choice by optimizing the convergence rate.

A standard way to find a worst-case  optimal step-size is to minimize an upper bound.  
However,  with the above parallelism,  there  exists  a hidden obstacle. 
To see this, consider the following convergence rate bounds (hold whenever  $\mathcal{F}$ is firmly non-expansive, see later Proposition  \ref{prop_convergence}):
\begin{align}
\Vert \bm{\zeta}^{k+1} - \bm{\zeta}^k  \Vert^2 		\,\,&\leq\,\, \frac{1}{k+1} \Vert \bm{\zeta}^\star - \bm{\zeta}^0  \Vert^2, \label{form01} \\  
\Vert \bm{y}^{k+1} - \bm{y}^k  \Vert^2 				\,\,&\leq\,\,  \frac{1}{k+1} \Vert \bm{y}^\star - \bm{y}^0  \Vert^2.
\end{align}
The  obstacle is that the above two upper bounds are different in general.  
If   the $\rho$-distanced parallelism is not removed, then directly  minimizing  the upper bounds with  respect to any  $\rho$-related  variable will lead to contradicted step-size choices in general. 

As an example, recall from \eqref{class_fix} the  two fixed-point expressions in the current literature, restated here as
\begin{equation}
\bm{\zeta}_1^\star = \mathcal{A}\bm{x}^\star +\bm{\lambda}^\star/\gamma, 
\qquad \,\,
\bm{\zeta}_2^\star  = \gamma \mathcal{A}\bm{x}^\star + \bm{\lambda}^\star,
\end{equation}
which are related via $ \bm{\zeta}_2^\star  = \gamma \bm{\zeta}_1^\star$, leading to  the following two optimization problems:
\begin{align}
\underset{\gamma >0}{\text{minimize}}\,\, & 	\Vert \mathcal{A}\bm{x}^\star +\bm{\lambda}^\star/\gamma  - \bm{\zeta}_1^0\Vert^2 , \\
\underset{\gamma >0}{\text{minimize}}\,\,  &	\Vert \gamma\mathcal{A}\bm{x}^\star +\bm{\lambda}^\star - \bm{\zeta}_2^0\Vert^2 ,
\end{align}
which will yield different optimal solutions $ \gamma^\star $ that are problematic,  see more details in  Sec.  \ref{sec_fail}.

\subsection{Unscaling guarantee}
The key to our success is to guarantee no step-size related extra scaling   exists when  we optimize a convergence rate upper bound.

\subsubsection{a conventional success}
Before providing our solution, we  briefly review a typical success in  the literature, which implicitly  avoided the above  scaling issue.

Indeed, a natural way to address 
such an extra scaling  is by division. 
Consider optimizing the following ratio:
\begin{equation}\label{fac}
\delta
= \underset{\{k | \bm{\zeta}^{k+1}\neq \bm{\zeta}^{k}\}}{\sup} \, \frac{\Vert \bm{\zeta}^{k+1} - \bm{\zeta}^\star \Vert}{\Vert \bm{\zeta}^{k} - \bm{\zeta}^\star \Vert}
= \underset{\{k | \bm{y}^{k+1}\neq \bm{y}^{k}\}}{\sup} \, \frac{\Vert \bm{y}^{k+1} - \bm{y}^\star \Vert}{\Vert \bm{y}^{k} - \bm{y}^\star \Vert},
\end{equation}
which holds for any non-zero $\rho$-distanced parallel sequences $\bm{y}^{k} = \rho\bm{\zeta}^{k} $.
Such a ratio-factor optimizing approach is proposed in \cite{ghadimi2014optimal} for quadratic programming, and  \cite{shi2014linear} for  decentralized consensus problems. 
However, there are  two notable drawbacks:  

\vspace{5pt}

$ \bullet $ (i) First,  the factor $ \delta $ needs to be the same for every iteration $ k $, otherwise  optimizing it only improves a  few middle  steps, and  the whole convergence  rate does not  necessarily  improve.  
This requirement can be satisfied when  the convergence rate is linear. 
However, a  linear rate  does not hold in general  (requires a strong assumption, see later Proposition \ref{prop_convergence}). That said,  the  step-size choice obtained from here is not universally applicable.

\vspace{5pt}
$ \bullet $ (ii) More importantly, their successes require analytical tractability, i.e., the fixed-point operator $ \mathcal{F} $ (equivalently, all algorithm iterates/sub-problems) admits an explicit analytical form.  
This requirement is not feasible for a large number of convex programs, which again limits its scope.

\vspace{5pt}

In view of these  drawbacks, a general  step-size selection scheme seems not attainable following this  ratio-factor optimizing path.

\subsubsection{idea sketch}\label{sec_2_5_2}
A natural question  is --- apart from the above division approach, is there another way to guarantee no extra scaling?

The answer is affirmative, by appealing to  symmetry.
Given a certain  fixed-point  operator (an observation), suppose an extra scaling exists. 
Then, it can be characterized into 
\begin{equation}\label{sym_set}
\mathcal{F}_{\text{obs}} = \rho\mathcal{I}\circ  \mathcal{F}_0 \circ \frac{1}{\rho}\mathcal{I},\quad \rho \neq 0,
\end{equation}
where  $ \mathcal{F}_0 $  denotes the underlying unscaled operator.

The insight is that the  extra scaling structure $ \rho\mathcal{I}\circ  (\cdot)\circ \frac{1}{\rho}\mathcal{I} $ is  not symmetric. 
Then,  consider  its  adjoint operator 
$  \mathcal{F}^*_{\text{obs}}:   \mathbb H \rightarrow \mathbb H $, which can be abstractly written as
\begin{equation}\label{adj}
\mathcal{F}^*_{\text{obs}} = \bigg(\rho\mathcal{I}\circ  \mathcal{F}_0 \circ \frac{1}{\rho}\mathcal{I}\bigg)^* =  \frac{1}{\rho} \mathcal{I}\circ   \mathcal{F}^*_0  \circ \rho\mathcal{I}.
\end{equation}

Compare \eqref{sym_set}  and \eqref{adj}, we see that the extra  scaling  structure  changes. 
Now, suppose the following always holds for all choices of step-size parameter $ \rho \neq 0 $:
\begin{equation}\label{claim}
\mathcal{F}_{\text{obs}}  =  \mathcal{F}^*_{\text{obs}}, \tag{to prove}
\end{equation}
which  clearly includes  the  case $ \mathcal{F}_0   =  \mathcal{F}^*_0 $, since we assume it holds $ \forall \rho \neq 0 $.

Then, under this assumption,  clearly we can conclude that  $  \mathcal{F}_{\text{obs}} $ does not contain the  extra scaling structure.
Relation \eqref{claim}  will be concretely proved  later  in Sec. \ref{self_d}, and  the proof  relies on some  properties  of our domain  parametrization, which we first  establish them below.

\subsection{Domain parametrization: useful properties}
\subsubsection{non-zero scalar parameter}
Here, we  show that the  domain-parametrized proximal  operator is  well-defined. 
\begin{prop}\label{prop_firm2}
	The domain-parametrized  proximal operator $ \textbf{Prox}_{f \rho}  $	defined in \eqref{def_right} is firmly nonexpansive and single-valued.	
\end{prop}
\begin{proof}
	To start, we need some lemmas.
	\begin{lem}\cite[Proposition 6.19]{bauschke2017convex}\label{lem1}
		Let $ C$ be a convex subset of $ \mathbb H $, let $ \mathbb K $ be a real Hilbert space, let  $  \mathcal{L} \in \mathfrak{B}(\mathbb{H}, \mathbb{K}) $
		and let  $D$ be a convex subset of $ \mathbb K $. Let    $\text{int}(\cdot) $ and $ \text{sri}(\cdot)$ denote  the interior and the strong relative interior, respectively.
		Suppose $ D \bigcap \text{int}\,\, \mathcal{L}(C) \neq \emptyset$ or $ \mathcal{L}(C) \bigcap \text{int}\,\, D \neq \emptyset$. Then, $ 0\in \text{sri}\, (D - \mathcal{L}(C) ) $.
	\end{lem}
	
	\vspace{1pt}
	
	\begin{lem}\cite[Corollary 16.53]{bauschke2017convex}\label{dom_ran}
		Let $ f $ be a CCP function and $  \mathcal{L} \in \mathfrak{B}(\mathbb{H}, \mathbb{K}) $. Suppose $ 0\in \text{sri}\, (\text{dom}(f) -\text{ran}(\mathcal{L}) ) $. Then, 
		\begin{equation*}
		\partial (f\circ\mathcal{L}) = \mathcal{L}^*\circ\partial f\circ\mathcal{L}.
		\end{equation*}
	\end{lem}

	\begin{lem}\cite[Proposition 23.25]{bauschke2017convex}\label{coro_maxm}
		Let $ \mathbb K $ be a real Hilbert space,  suppose  $ \mathcal{L} \in \mathfrak{B}(\mathbb{H}, \mathbb{K}) $ is such that $ \mathcal{L}\mathcal{L}^* $ is invertible, let $\mathcal{A}: \mathbb{K} \rightarrow2^\mathbb{K}$ be maximal monotone, and set $ \mathcal{B} = \mathcal{L}^*\mathcal{A}\mathcal{L} $. Then, 	
		$\mathcal{B}: \mathbb{H}\rightarrow 2^\mathbb{H}$ is maximal monotone.	
	\end{lem}

	For the above 3 lemmas,
	substitute  the mapping  $ \mathcal{L}  $ there with our scalar one $ \rho \mathcal{I} \,\,(\rho \neq 0) $.  
	First, a non-zero scaling is clearly bijective, the two spaces in the above lemmas hence coincide, i.e., $ \mathbb{H}  = \mathbb{K} $.
	Moreover, the range of  a  non-zero scalar map is clearly  the entire Hilbert space, i.e., $ \text{ran}(\rho \mathcal{I})  =  \mathbb{K} =   \mathbb{H}$.

	In view of Lemma \ref{lem1}, substitute  the set  $ \mathcal{D}  $ there with $ \text{dom}(f)  $ (subset of $ \mathbb{H} $).
	Then,  
	$ \text{ran}( \rho \mathcal{I} ) \bigcap \text{int}\,\, \text{dom}(f)  = \mathbb{H} \, \bigcap \text{int}\,\, \text{dom}(f) \neq \emptyset$, 	 which  holds unless $ \text{dom}(f)  = \emptyset$  (omitted by assumption).
	Therefore, we arrive at  $ 0\in \text{sri}\, (\text{dom}(f) -\text{ran}( \rho \mathcal{I} ) ) $. 
	
	Then, invoke  Lemma \ref{dom_ran}, we have
	$ \partial (f\rho) =   \rho \mathcal{I} \circ (\partial f) \circ  \rho \mathcal{I}  $. 
	By Lemma \ref{coro_maxm}, the operator $ \partial (f\rho) $ is  maximal monotone.	
	
	Now, invoke  definition  \eqref{def_right},
	\begin{align}
	\textbf{Prox}_{f \rho} (\bm{v}) 
	= \underset{\bm{z}}{\text{argmin}} \,\, f(\rho\bm{z} )+ \frac{1}{2}\Vert \bm{z}- \bm{v}\Vert^2
	\Longleftrightarrow& 	\,\bm{v} - \bm{z}  \in  \partial (f\rho) \,(\bm{z}),   \nonumber\\
	\Longleftrightarrow&    \,\bm{z}   = \big(\mathcal{I} + \partial (f\rho)  \big)^{-1}\,(\bm{v}).
	\end{align}
	Owing to $ \partial (f\rho) $  being  maximal monotone, the proximal operator is single-valued and  firmly nonexpansive, see e.g. \cite{bauschke2017convex}. The proof is now concluded.
\end{proof}

\subsubsection{injective operator parameter}\label{sec_inj}
Here, we show that an  injective parameter is feasible in two special cases. This aspect enables some  simple, unified analysis frameworks for ADMM and primal-dual problems in the next sections.
\begin{prop}\label{pro_inj}[single-valued]
	Given $ f \in \Gamma_0 (\mathbb H) $,   $\text{dom}(f) \neq \emptyset$ and $ \mathcal D  \in \mathfrak B (\mathbb H,  \mathbb K) $ being injective.  
	Let      $ \mathcal D^*  \in \mathfrak B (  \mathbb K, \mathbb H) $  and  $ \mathcal D^{-1}  \in \mathfrak B (  \mathbb K, \mathbb H) $   denote    the adjoint and inverse mappings,  respectively. Then,
	the following  proximal operators are firmly non-expansive and single-valued:
	\begin{align}\label{eq_inv}
	\textbf{Prox}_{f\mathcal{D}^{-1}} (\bm{v})
	=&  \,\,\underset{\bm{z}_1 \in \mathbb K }{\text{argmin}} \,\, f\left(\mathcal{D}^{-1}\bm{z}_1\right)+ \frac{1}{2}\Vert \bm{z}_1 - \bm{v}\Vert^2  ,\\
	\textbf{Prox}_{f\mathcal{D}^*} (\bm{v})
	=&  \,\,\underset{\bm{z}_1 \in \mathbb K }{\text{argmin}} \,\, f\left(\mathcal{D}^*\bm{z}_1\right)+ \frac{1}{2}\Vert \bm{z}_1 - \bm{v}\Vert^2  .
	\end{align}
\end{prop} 
\begin{proof}
	The proof follows  similar procedures as  Proposition \ref{prop_firm2}, with different substitutions.
	We will  invoke the  same lemmas  and  substitute  the mapping  $ \mathcal{L}  $ there with   $ \mathcal D^{-1}   $ and  $ \mathcal D^*   $. To avoid repeating, below we no  longer state such lemmas.
	
	For both cases, $ \mathcal L   =  \mathcal D^{-1}  $ and $ \mathcal L   =  \mathcal D^*  $,
	we always have $ \text{ran}(\mathcal L) =   \text{dom}\, \mathcal  D =   \mathbb H $.
	It follows that   $ \text{ran}(\mathcal L) \bigcap \text{int}\,\, \text{dom}(f)  =   \mathbb H   \bigcap \text{int}\,\, \text{dom}(f) \neq \emptyset$,
	which  holds unless $ \text{dom}(f)  = \emptyset$  (omitted by assumption).

	Now, invoke  Lemma \ref{dom_ran}, we have
	$ \partial (f\circ \mathcal L) =  \mathcal L^* \, (\partial f) \, \mathcal L $. 	
	All what  left is regarding  Lemma \ref{coro_maxm}, where we need   $ \mathcal{L}\mathcal{L}^* $ there being invertible.
	Let $ \mathcal{L} =  \mathcal{D}^{-1}$, then
	\begin{equation}
	\mathcal{L}\mathcal{L}^* = \mathcal{D}^{-1}  ({\mathcal{D}}^*)^{-1} = (\mathcal{D}^*\mathcal{D})^{-1}, 
	\end{equation}
	Let $ \mathcal{L} =  \mathcal{D}^*$, then
	\begin{equation}
	\mathcal{L}\mathcal{L}^* = \mathcal{D}^*  (\mathcal{D}^*)^* = \mathcal{D}^*\mathcal{D}.
	\end{equation}
	Since $ \mathcal{D} $ is assumed  injective, the above two cases  are  clearly  both  invertible. 
	The proof is therefore concluded.
\end{proof}

Next, we establish some alternative/equivalent  characterizations  for the above $\textbf{Prox}_{f\mathcal{D}^{-1}} $ and  $ 	\textbf{Prox}_{f\mathcal{D}^*}  $, which will be used later.
To start, we need some lemmas.
\begin{defi}\label{defi0}
	The  orthogonal  complement  of a subset $ C $ of $ \mathbb{H} $ is: $  C^\perp = \{ \bm{x} \in  \mathbb{H} \,|\,  \langle \bm{x}, \bm{y}  \rangle = 0, \, \forall \bm{y} \in  C \}  $.
\end{defi}
\begin{lem}\cite[Fact 2.25]{bauschke2017convex}\label{lem_ker}
	Let $ \mathcal{T} \in   \mathfrak B (\mathbb H,  \mathbb K) $, and let $\text{ker}\, \mathcal{T} = \{ \bm{x} \in \mathbb{H} \, |\,\mathcal T \,\bm{x} = \bm{0}\} $ be the kernel  of $ \mathcal{T} $. 
	Then, $ \text{ker}\, \mathcal{T}^* = (\text{ran}\, \mathcal T)^\perp $.
\end{lem}
\begin{lem}\label{lem_inj}
	Let $ \mathcal D  \in \mathfrak B (\mathbb H,  \mathbb K) $ be injective. Then, its adjoint inverse mapping $ 	(\mathcal{D}^*)^{-1}: \mathbb H \rightarrow \mathbb K $  is also injective.
\end{lem}
\begin{proof}
	First, let us  note that
	\begin{equation}\label{ran_dom}
	\text{ran}\, \mathcal \mathcal D^{-1}  = \text{dom}\, \mathcal  D  = \mathbb{H}.
	\end{equation}
	Then, by Definition \ref{defi0}, we obtain
	\begin{align}
	(\text{ran}\, \mathcal \mathcal D^{-1})^\perp  
	= (\text{dom}\, \mathcal \mathcal D)^\perp 
	= \{ \bm{x}\in  \mathbb{H}  \,|\,  \langle \bm{x}, \bm{y}  \rangle = 0, \,\forall \bm{y} \in  \text{dom}\, \mathcal{D}   = \mathbb{H}  \}.
	\end{align}
	Since the above defined $ \bm{y} $ takes value in the entire space $ \mathbb{H} $, the above set clearly only contains the zero vector, i.e., 
	\begin{equation}
	(\text{ran}\, \mathcal \mathcal D^{-1})^\perp = \{\bm{0}\}.
	\end{equation}
	By Lemma \ref{lem_ker}, the above implies 
	\begin{equation}
	\text{ker}\, 	(\mathcal{D}^*)^{-1}
	=  (\text{ran}\, \mathcal \mathcal D^{-1})^\perp = \{\bm{0}\}.
	\end{equation}
	That implies $ (\mathcal{D}^*)^{-1}  $ being injective, which   concludes the proof.
\end{proof}

\begin{lem}[transformation]\label{prop_trans}
	Given $ f \in \Gamma_0 (\mathbb H) $,   $\text{dom}(f) \neq \emptyset$ and $ \mathcal D  \in \mathfrak B (\mathbb H,  \mathbb K) $ being injective.  Then,
	\begin{align}
	\,\mathcal{D} \, &\underset{\bm{z}_2 \in \mathbb H}{\text{argmin}} \,\, f(\bm{z}_2)+ \frac{1}{2}\Vert \mathcal{D}\bm{z}_2 - \bm{v}\Vert^2
	= \underbracket{\underset{\bm{z}_1 \in \mathbb K }{\text{argmin}} \,\, f\left(\mathcal{D}^{-1}\bm{z}_1\right)+ \frac{1}{2}\Vert \bm{z}_1 - \bm{v}\Vert^2}_{\textbf{Prox}_{f\mathcal{D}^{-1}} (\bm{v})} ,\nonumber\\
	\,(\mathcal{D}^*)^{-1} \, &\underset{\bm{z}_2 \in \mathbb H}{\text{argmin}} \,\, f(\bm{z}_2)+ \frac{1}{2}\Vert (\mathcal{D}^*)^{-1} \bm{z}_2 - \bm{v}\Vert^2
	= \underbracket{\underset{\bm{z}_1 \in \mathbb K }{\text{argmin}} \,\, f\left(\mathcal{D}^*\bm{z}_1\right)+ \frac{1}{2}\Vert \bm{z}_1 - \bm{v}\Vert^2}_{	\textbf{Prox}_{f\mathcal{D}^*} (\bm{v})}  .\nonumber
	\end{align}
\end{lem}
\begin{proof}
	The proof is  straightforward  by a  variable substitution. 
	
	For the  first  relation, we prove from  its left-hand-side.  
	Let $ \bm{z}_1  = \mathcal{D} \bm{z}_2$, which is a  one-to-one correspondence, owing to  $ \mathcal{D} $  being  injective.  
	Hence, we obtain
	\begin{equation}
	\underset{\bm{z}_2 \in \mathbb H}{\text{argmin}} \,\, f(\bm{z}_2)+ \frac{1}{2}\Vert \mathcal{D}\bm{z}_2 - \bm{v}\Vert^2
	= \underset{(\mathcal{D}^{-1}\bm{z}_1) \in \mathbb H}{\text{argmin}} \,\, f(\mathcal{D}^{-1}\bm{z}_1)+ \frac{1}{2}\Vert \bm{z}_1 - \bm{v}\Vert^2,
	\end{equation}
	which yields  a  minimizer relation $ \bm{z}_2^\star  =   \mathcal{D}^{-1}\bm{z}_1^\star $, which can   be  rewritten into $\mathcal{D} \bm{z}_2^\star  =   \bm{z}_1^\star $. This  gives  the first relation above.
	
	Same arguments  hold for the  second case, where mapping $ (\mathcal{D}^*)^{-1} $ is injective by   Lemma \ref{lem_inj}.
	We omit the same arguments to avoid repeating.
	The  proof  is  now concluded.
\end{proof}

\begin{lem}[single-valued partial evaluation]\label{lam_sin_p}
	Given $ f \in \Gamma_0 (\mathbb H) $ and $\text{dom}(f) \neq \emptyset$. Let $ \mathcal D  \in \mathfrak B (\mathbb H,  \mathbb K) $ be injective.  Then, the  following  two evaluations are single-valued:
	\begin{align}
	&\underset{\bm{z}_2 \in \mathbb H}{\text{argmin}} \,\, f(\bm{z}_2)+ \frac{1}{2}\Vert \mathcal{D}\bm{z}_2 - \bm{v}\Vert^2 , \label{re001}\\
	&\underset{\bm{z}_2 \in \mathbb H}{\text{argmin}} \,\, f(\bm{z}_2)+ \frac{1}{2}\Vert (\mathcal{D}^*)^{-1} \bm{z}_2 - \bm{v}\Vert^2,  \label{re002}
	\end{align}
	or equivalently, the following two mappings are single-valued:
	\begin{align}\label{resol}
	(\mathcal{D}^*\mathcal{D} + \partial f )^{-1} \mathcal{D}^*, \quad
	((\mathcal{D}^*\mathcal{D})^{-1} + \partial f )^{-1} \mathcal{D}^{-1}.  
	\end{align}
\end{lem}
\begin{proof}
	Following from  Lemma \ref{prop_trans},
	\begin{equation}
	\mathcal{D} \,\, \underbracket{\underset{\bm{z}_2 \in \mathbb H}{\text{argmin}} \,\, f(\bm{z}_2)+ \frac{1}{2}\Vert \mathcal{D}\bm{z}_2 - \bm{v}\Vert^2}_{\eqref{re001}}
	\,\,=\,\,  	\textbf{Prox}_{f\mathcal{D}^{-1}} (\bm{v}),
	\end{equation}
	where the  right-hand side    is single-valued by Proposition \ref{pro_inj}. 
	Since $ \mathcal D  $ is assumed  injective, i.e., a one-to-one correspondence, then its  domain element, corresponding to
	\eqref{re001}, has to be single-valued.
	
	Same arguments hold for the symmetric case \eqref{re002} and is omitted to avoid repeating. 
	The expressions in \eqref{resol} follows instantly by the first-order optimality condition.
	The proof is now concluded.
\end{proof}

\subsubsection{Domain-parametrized Moreau identity}\label{sec_mor}
\begin{prop}\label{pmoreau}
	Given $ f \in \Gamma_0 (\mathbb H) $  and  $\text{dom}(f) \neq \emptyset$. Let $ \mathcal D  \in \mathfrak B (\mathbb H,  \mathbb K) $ be injective.  Then, the following relation holds:
	\begin{equation}\label{Moreau identity}
	\bm{v}  =\,\textbf{Prox}_{f\mathcal{D}^{-1}} (\bm{v}) +  \textbf{Prox}_{f^*{\mathcal{D}^*}} (\bm{v}), \quad  \forall \bm{v} \in  \mathbb K.
	\end{equation}
\end{prop}
\begin{proof}
	First, by Proposition \ref{pro_inj}, the above  $\textbf{Prox}_{f\mathcal{D}^{-1}} $ and  $ 	\textbf{Prox}_{f\mathcal{D}^*}  $ are well-defined. 
	Now,  let	
	\begin{align}
	\bm{x}^\star & = \underset{\bm{x}}{\text{argmin}} \,\, \,\, f(\bm{x}) \,+ \frac{1}{2}\Vert \mathcal{D}\bm{x} - \bm{v}\Vert^2 , \label{x_part}
	\\
	\bm{y}^\star & = \underset{\bm{y}}{\text{argmin}} \,\, f^*(\bm{y})+ \frac{1}{2}\Vert (\mathcal{D}^*)^{-1}\bm{y} - \bm{v}\Vert^2, \label{y_part}
	\end{align}
	which are both single-valued by  Lemma \ref{lam_sin_p}.
	By Lemma  \ref{prop_trans}, our goal can be converted    to proving the following relation:
	\begin{equation}\label{goal}
	\bm{v} = \mathcal{D}\bm{x}^\star +(\mathcal{D}^*)^{-1} \bm{y}^\star
	\iff 	\bm{y}^\star = \mathcal{D}^*(\bm{v} - \mathcal{D}\bm{x}^\star). \tag{goal}
	\end{equation}
	
	To this end, consider $ \bm{x}^\star $  first.  Invoke the first-order optimality condition for \eqref{x_part}, we obtain
	\begin{align}
	\mathcal{D}^*\bm{v} -\mathcal{D}^*\mathcal{D}\,\bm{x}^\star &\in \partial {f}(\bm{x}^\star) ,\\
	\iff \hspace{4.5cm}	 \bm{x}^\star &\in \partial f^*(  \mathcal{D}^*\bm{v} - \mathcal{D}^*\mathcal{D}\bm{x}^\star ), \\
	\iff\hspace{2cm}	 \mathcal{D}^{-1}(\bm{v} -  \bm{v} + \mathcal{D}\bm{x}^\star)  &\in \partial f^*(  \mathcal{D}^*\bm{v} - \mathcal{D}^*\mathcal{D}\bm{x}^\star ), \\
	\iff\,\,\,	 \mathcal{D}^{-1}(\bm{v} - (\mathcal{D}^*)^{-1}\underbracket{\mathcal{D}^*(  \bm{v} - \mathcal{D}\bm{x}^\star)}_{\eqdef \, \bm{t}}  ) &\in \partial f^*(  \underbracket{\mathcal{D}^*\bm{v} - \mathcal{D}^*\mathcal{D}\bm{x}^\star}_{\eqdef\, \bm{t}} ) ,\\
	\iff	 \hspace{2.1cm}	\mathcal{D}^{-1}(\bm{v} - 	(\mathcal{D}^*)^{-1}\bm{t}) &\in \partial f^*(\bm{t}), \\
	\iff	 \hspace{4.75cm}	\bm{t} &= ((\mathcal{D}^*\mathcal{D})^{-1} + \partial f^* )^{-1} \mathcal{D}^{-1} (\bm{v}),
	\label{v-sx}
	\end{align}
	with the last line equality owes to the single-valued  property, recall Lemma \ref{lam_sin_p}, where $ \bm{t} \eqdef \mathcal{D}^*(  \bm{v} - \mathcal{D}^*\bm{x})  $ is  a variable substitution.
	
	Now, consider $ \bm{y}^\star $.  Invoke the first-order optimality condition for \eqref{y_part}, we obtain
	\begin{align}
	\mathcal{D}^{-1}(\bm{v} - 	(\mathcal{D}^*)^{-1}\bm{y}^\star) &\in \partial f^*(\bm{y}^\star)  ,\\
	\iff \hspace{2.8cm}	  		\bm{y}^\star &= \big((\mathcal{D}^*\mathcal{D})^{-1} + \partial f^* \big)^{-1} \mathcal{D}^{-1}  (\bm{v}), \label{f^*_Y}
	\end{align}
	with the last line equality owing to the single-valued  property, recall Lemma \ref{lam_sin_p}.
	
	Compare \eqref{v-sx} and \eqref{f^*_Y}, the expressions are exactly the same.  That is,
	\begin{equation}
	\bm{y}^\star = \bm{t} = \mathcal{D}^*(\bm{v} - \mathcal{D}\bm{x}^\star),
	\end{equation}
	which proves \eqref{goal}. The proof is now concluded.
\end{proof}

\subsection{Domain-parametrized ADMM}\label{sec_admm}
Here, we express  the  ADMM algorithm in terms of our  domain-parametrized proximal operator.  We  see some useful new results, particularly,  a great  analysis convenience.

Consider a general convex problem:
\begin{align}\label{ab_prob}
\underset{\bm{x},\bm{z}}{\text{minimize}}& \quad   f(\bm{x}) + g(\bm{z}), \nonumber\\ 
\text{subject\,to}& \quad \mathcal{A}\bm{x} - \mathcal{B}\bm{z} = \bm{c} ,  \tag{problem}
\end{align}
with variables $ \bm{x} \in  \mathbb H $, $ \bm{z} \in  \mathbb P $, $ \bm{c} \in  \mathbb K$ and functions  $ f\in \Gamma_0 (\mathbb{H}) $, $ g\in \Gamma_0 (\mathbb{P}) $, where $ \mathcal{A} \in  \mathfrak B (\mathbb H,  \mathbb K)$ and  $ \mathcal{B} \in  \mathfrak B (\mathbb P,  \mathbb K)$ are injective.
We assume a solution always exists.

The augmented Lagrangian can be written as
\begin{equation}
\mathcal{L}_\rho(\bm{x},\bm{z},\bm{\lambda}) 
\,=\,  f(\bm{x}) + g (\bm{z})  + \frac{1}{2}\Vert \rho (\mathcal{A}\bm{x} - \mathcal{B}\bm{z} - \bm{c}) + \bm{\lambda}/\rho\Vert^2,  
\end{equation}	
where $  \rho \neq 0$ is a non-zero step-size parameter. 
The  ADMM iterates are
\begin{align}
\bm{x}^{k+1} =\,\,&  \, \underset{\bm{x} \in \mathbb H}{\text{argmin}} \,\, f(\bm{x})+ \frac{1}{2}\Vert  \rho(\mathcal{A}\bm{x} -   \mathcal{B}\bm{z}^k  - \bm{c})  + \bm{\lambda}^k/\rho  \Vert^2, \nonumber\\				
\bm{z}^{k+1} =\,\, &\, \underset{\bm{z} \in \mathbb P}{\text{argmin}} \,\, g(\bm{z})+ \frac{1}{2}\Vert  \rho(\mathcal{A}\bm{x}^{k+1} -   \mathcal{B}\bm{z}  - \bm{c})  + \bm{\lambda}^k/\rho  \Vert^2, \nonumber\\
\bm{\lambda}^{k+1} =\,\,  &\,\bm{\lambda}^{k} + \rho^2( \mathcal{A}\bm{x}^{k+1} - \mathcal{B}\bm{z}^{k+1} - \bm{c}). 	
\end{align}

\subsubsection{New proximal steps}
Given $ \mathcal{A} \in  \mathfrak B (\mathbb H,  \mathbb K)$ and  $ \mathcal{B} \in  \mathfrak B (\mathbb P,  \mathbb K)$ being injective.  
Clearly,  $ \rho\mathcal{A}  $ and $ \rho\mathcal{B}  $ are also injective   $\forall \rho  \neq 0$,  and  therefore  both  the inverse  and   the adjoint mappings are
well-defined parameters  by  Proposition \ref{pro_inj}.

Then,  by Lemma \ref{prop_trans},  we  can rewrite  the  ADMM iterates  into:
\begin{align}\label{R-ADMM}
\rho\mathcal{A}	\bm{x}^{k+1} \,\,=\,\,&  \,\textbf{Prox}_{f(\rho\mathcal{A})^{-1}} \left( \,\,\,\,  \rho   \bm{c}  + \rho\mathcal{B}\bm{z}^k  - \bm{\lambda}^k/\rho\right) , \nonumber\\				
\rho\mathcal{B}\bm{z}^{k+1} \,\,=\,\, &   \,\textbf{Prox}_{g(\rho\mathcal{B})^{-1}} \left(-\rho\bm{c} + \rho\mathcal{A}\bm{x}^{k+1}  + \bm{\lambda}^k/\rho\right), \nonumber\\
\bm{\lambda}^{k+1} \,\,=\,\,  &\,  \bm{\lambda}^{k} + \rho^2( \mathcal{A}\bm{x}^{k+1} - \mathcal{B}\bm{z}^{k+1} - \bm{c}). 	\tag{prox. ADMM}
\end{align}

We may further  simplify the  above via the following variable substitutions:
\begin{equation}
\widetilde{\bm{x}}^{k} \eqdef \rho\mathcal{A}\bm{x}^{k}, \qquad \widetilde{\bm{z}}^{k} \eqdef \rho\mathcal{B}\bm{z}^{k}, \qquad \widetilde{\bm{\lambda}}^{k} \eqdef \frac{1}{\rho} \bm{\lambda}^{k},\qquad
\widetilde{\bm{c}} \eqdef \rho\bm{c}.
\end{equation}
This leads to the following    \textit{scaled form}:
\begin{align}\label{scaled_from}
\widetilde{\bm{x}}^{k+1} =\,\,& \textbf{Prox}_{f(\rho\mathcal{A})^{-1}} \left( \,\,\, \widetilde{\bm{c}} \, + \,\widetilde{\bm{z}}^k  \,\,\, - \widetilde{\bm{\lambda}}^{k}\right),  \nonumber\\				
\widetilde{\bm{z}}^{k+1} =\,\, &  \textbf{Prox}_{g(\rho\mathcal{B})^{-1}} \left(-\widetilde{\bm{c}} + \widetilde{\bm{x}}^{k+1}  + \widetilde{\bm{\lambda}}^{k}\right), \nonumber\\
\widetilde{\bm{\lambda}}^{k+1} =\,\,  &\widetilde{\bm{\lambda}}^{k} + \widetilde{\bm{x}}^{k+1} - \widetilde{\bm{z}}^{k+1} - \widetilde{\bm{c}}. 	\tag{scaled form}
\end{align}

\begin{remark}[enabled proximal  steps]
	To our  knowledge, the above  proximal expression is  not available  in  the  current literature (regarding  the general constrained  form \eqref{ab_prob}), see comments from Radu Ioan Bot and  Erno Robert Csetnek \cite{boct2019admm}: `\textit{Generally, the minimization
		with respect to the variable x does not lead to a proximal step}', with  workarounds on  this  issue in e.g. \cite{shefi2014rate,banert2021fixing}.
\end{remark}

\subsubsection{Convergence: a direct establishment}\label{sec_conver}
In the literature,  a typical way to establish the ADMM convergence for the general  form \eqref{ab_prob}  is by first converting to an unconstrained form, either via an \textit{infimal post-composition} \cite[sec. 3.1]{ryu2022large}, or    dualization \cite[sec. 3.2]{ryu2022large}, \cite{poon2019trajectory}. 
Then, convergence is shown via  the   Douglas-Rachford splitting (DRS) iterates.

Here,  owing  to the direct proximal steps \eqref{scaled_from}, we are able to establish a  direct characterization (no conversion to  DRS).
To start, we need two lemmas. 
\begin{lem}\label{post_comp}
	Let $ h(\bm{z}) = f(\bm{z}) - \langle \bm{c}, \bm{z} \rangle$.  Then,
	\begin{equation}
	\textbf{Prox}_{h\rho^{-1}} (\bm{v}) = \textbf{Prox}_{{f}\rho^{-1}}  \big(\,  \bm{c} + \bm{v} \,\big).
	\end{equation}
	\begin{proof}
		The proof follows straightforwardly from definition:
		\begin{align}
		\textbf{Prox}_{h\rho^{-1}} (\bm{v}) 
		=&  \underset{\bm{z}_1}{\text{argmin}} \,\, h(\rho^{-1}\bm{z}_1)+ \frac{1}{2}\Vert\bm{z}_1 - \bm{v}\Vert^2 \nonumber\\
		=&  \underset{\bm{z}_1}{\text{argmin}} \,\,  f(\rho^{-1}\bm{z}_1) - \langle \bm{c}, \bm{z}_1 \rangle+ \frac{1}{2}\Vert\bm{z}_1 - \bm{v}\Vert^2 \nonumber\\
		=&  \underset{\bm{z}_1}{\text{argmin}} \,\,  f (\rho^{-1}\bm{z}_1) + \frac{1}{2}\Vert\bm{z}_1 -\bm{c} - \bm{v}\Vert^2 \nonumber\\
		=&  \textbf{Prox}_{{f}\rho^{-1}} (  \bm{c} + \bm{v}  ) ,
		\end{align}
		which concludes the proof.
	\end{proof}
\end{lem}

\begin{lem}[1/2-averaging, see e.g.  \cite{ryu2022large,bauschke2017convex}]\label{1/2-ave}
	Let $ \mathcal{T}: \mathbb H \rightarrow \mathbb H$ be non-expansive. Then, the 1/2-averaging $ \frac{1}{2}\mathcal{I} + \frac{1}{2}\mathcal{T} $ is firmly non-expansive. 
	Reversely,  let $ \mathcal{F} : \mathbb H \rightarrow \mathbb H$ be firmly non-expansive, then $ 2\mathcal{F} - \mathcal{I}  $ is non-expansive.
\end{lem}

Now,  we are in position to establish  the convergence.
\begin{prop}[convergence]\label{fpi}
	\eqref{R-ADMM} and  \eqref{scaled_from} admit the following fixed-point characterization:
	\begin{equation}\label{fpi_admm}
	\bm{\zeta}^{k+1} = \mathcal{F}_\text{ADMM} \,\, \bm{\zeta}^{k},
	\end{equation}
	with  $ 	\bm{\zeta}^{k+1} \,\,=\,\,   \rho\mathcal{A}\bm{x}^{k+1} +  \frac{1}{\rho} \bm{\lambda}^{k} $,  where
	\begin{align}\label{fp}
	\mathcal{F}_{ADMM} \,\,&= \,\, \frac{1}{2}(2\textbf{Prox}_{\bar{f}} - \mathcal{I})\circ(2\textbf{Prox}_{ \widehat{g}} - \mathcal{I}) + \frac{1}{2}\mathcal I ,\nonumber
	\end{align}
	and where
	\begin{align}
	\bar{f}(\bm{y}) &\eqdef  f \circ (\rho\mathcal{A})^{-1}(\bm{y})  -  \langle \rho\bm{c}, \bm{y} \rangle,  \\  
	\widehat{g}(\bm{y}) &\eqdef g \circ   (\rho\mathcal{B})^{-1}\,(\bm{y}) +   \langle \rho\bm{c}, \bm{y} \rangle,  \quad  \forall\bm{y}  \in  \mathbb{K}.
	\end{align}
	The sequence $ \{\bm{\zeta}^{k}\} $ converges to a fixed-point $ \bm{\zeta}^\star \in \text{Fix}\,\, \mathcal{F}_\text{ADMM}  $ (if it exists).
\end{prop}

\begin{proof}
	For light  of  notations, we prove via the \eqref{scaled_from}, with substitutions  restated here  as:
	\begin{equation}\label{scaled_variable}
	\widetilde{\bm{x}}^{k} \eqdef \rho\mathcal{A}\bm{x}^{k}, \qquad \widetilde{\bm{z}}^{k} \eqdef \rho\mathcal{B}\bm{z}^{k}, \qquad \widetilde{\bm{\lambda}}^{k} \eqdef \frac{1}{\rho}  \bm{\lambda}^{k}
	\qquad \widetilde{\bm{c}} \eqdef\rho \bm{c}.
	\end{equation}
	In view of the  $ \bm{\lambda} $-update in  \eqref{scaled_from}, we have
	\begin{align}
	\widetilde{\bm{\lambda}}^{k}
	\,\,	=\,\,  \widetilde{\bm{\lambda}}^{k-1} + \widetilde{\bm{x}}^{k} - \widetilde{\bm{z}}^{k} - \widetilde{\bm{c}} 
	\,\,	=\,\,   \bm{\zeta}^{k}  - \widetilde{\bm{z}}^{k} - \widetilde{\bm{c}}. 
	\end{align}
	Then,   we arrive at
	\begin{align}\label{pri_fpi}
	\bm{\zeta}^{k+1} =\,\,\,\,  &  \widetilde{\bm{x}}^{k+1} + \widetilde{\bm{\lambda}}^{k}  , \\
	=\,\,\,\,  & \textbf{Prox}_{f  (\rho\mathcal{A})^{-1}} (2\widetilde{\bm{c}} +  2 \widetilde{\bm{z}}^{k}  - \bm{\zeta}^{k}) + \bm{\zeta}^{k}  - \widetilde{\bm{z}}^{k} - \widetilde{\bm{c}} , \nonumber\\
	=\,\,\,\,  & \textbf{Prox}_{ f (\rho\mathcal{A})^{-1}} \bigg(\widetilde{\bm{c}} + 2\textbf{Prox}_{ g(\rho\mathcal{B})^{-1}}(\bm{\zeta}^{k} - \widetilde{\bm{c}} )- (\bm{\zeta}^{k} - \widetilde{\bm{c}} ) \bigg)  +  \nonumber\\
	& (\bm{\zeta}^{k} - \widetilde{\bm{c}}) - \textbf{Prox}_{ g(\rho\mathcal{B})^{-1}}(\bm{\zeta}^{k} - \widetilde{\bm{c}} ), \nonumber\\
	=\,\,\,\,  & \textbf{Prox}_{  {f}(\rho\mathcal{A})^{-1} }\circ \bigg(\,\widetilde{\bm{c}} + \big(2\textbf{Prox}_{ {g}(\rho\mathcal{B})^{-1}} - \mathcal{I}\big)   (\bm{\zeta}^{k} - \widetilde{\bm{c}} )\bigg) -\nonumber\\
	& \frac{1}{2} \bigg( \,\widetilde{\bm{c}} +
	\big( 2\textbf{Prox}_{ {g}(\rho\mathcal{B})^{-1}} - \mathcal{I}\big)  (\bm{\zeta}^{k} - \widetilde{\bm{c}} ) \bigg)  + \frac{1}{2}\bm{\zeta}^{k}, \nonumber\\
	=\,\,\,\,  & \underbracket{\left(\frac{1}{2}(2\textbf{Prox}_{\bar{f}} - \mathcal{I})\circ(2\textbf{Prox}_{ \widehat{g}} - \mathcal{I}) + \frac{1}{2}\mathcal I    \right)}_{	\mathcal{F}_{ADMM} } \,   \bm{\zeta}^{k},\nonumber
	\end{align}
	where the last line  uses Lemma \ref{post_comp} and substitutions:
	\begin{align}
	\bar{f}(\bm{y}) &\eqdef  f \circ (\rho\mathcal{A})^{-1}(\bm{y})  -  \langle \rho\bm{c}, \bm{y} \rangle,  \\  
	\widehat{g}(\bm{y}) &\eqdef g \circ   (\rho\mathcal{B})^{-1}\,(\bm{y}) +   \langle \rho\bm{c}, \bm{y} \rangle,  \quad  \forall\bm{y}  \in  \mathbb{K}.
	\end{align}
	By Lemma \ref{1/2-ave}, the  above  $ \mathcal{F}_{ADMM} $ is firmly non-expansive.	
	It follows that sequence $ \{ \bm{\zeta}^k\} $ converges to a fixed-point $ \bm{\zeta}^\star \in Fix\,\, \mathcal{F}_{ADMM} $ (if it exists),  see e.g. \cite{bauschke2017convex,poon2019trajectory}.

	At last, invoke definition \eqref{scaled_variable}, we obtain the expression:
	\begin{equation}\label{ADMM_fix}
	\bm{\zeta}^{k+1} 
	=\,\, \widetilde{\bm{x}}^{k+1} + \widetilde{\bm{\lambda}}^{k}
	=\,\,  \rho\mathcal{A}\bm{x}^{k+1} +  \frac{1}{\rho} \bm{\lambda}^{k}.
	\end{equation} The proof is now concluded.
\end{proof}

\subsection{Duality}\label{admm_self}
Here, we consider the dual characterization, which is key  to establishing  the  above fixed-point expression \eqref{ADMM_fix} being unscaled,  corresponding  to Sec. \ref{sec_2_5_2}.

\subsubsection{Domain-parametrized primal-dual characterization}
Here,  we dualize \eqref{ab_prob},  restated here as (primal problem):
\begin{align}\label{conv_1}
\underset{\bm{x},\bm{z}}{\text{minimize}}& \quad   f(\bm{x}) + g(\bm{z}), \nonumber\\ 
\text{subject\,to}& \quad \mathcal{A}\bm{x} - \mathcal{B}\bm{z} = \bm{c} ,
\end{align}
with variables $ \bm{x} \in  \mathbb H $, $ \bm{z} \in  \mathbb P $, $ \bm{c} \in  \mathbb K$ and functions  $ f\in \Gamma_0 (\mathbb{H}) $, $ g\in \Gamma_0 (\mathbb{P}) $, where $ \mathcal{A} \in  \mathfrak B (\mathbb H,  \mathbb K)$ and  $ \mathcal{B} \in  \mathfrak B (\mathbb P,  \mathbb K)$ are injective. We assume a solution always exists.

To start, we need one lemma.
\begin{lem}\label{lem01}
	Given  function $ g  \in \Gamma_0(\mathbb  H) $, let 	$ h (\cdot ) = g (\cdot) +  \langle \bm{c}, \cdot \rangle$. Then,
	\begin{equation}
	h^* (\cdot) = g^*(\cdot - \bm{c}).
	\end{equation}
\end{lem}
\begin{proof}
	By the Fenchel conjugate definition, 
	\begin{equation}
	h^* (\bm v) = \underset{\bm{z}}{\sup}\, \langle \bm{z}, \bm v \rangle - g (\bm{z}) - \langle \bm{c}, \bm{z} \rangle  
	= \underset{\bm{z}}{\sup}\, \langle \bm{z}, \bm v - \bm{c}\rangle - g (\bm{z}) =  g^* (\bm{v} - \bm{ c}),
	\end{equation}	
	which concludes the proof.
\end{proof}

\begin{prop}\label{prop_uni_pd}
	The  convex  program \eqref{conv_1} admits the following  primal-dual  characterization:
	\begin{align}
	&\underset{\widetilde{\bm x}}{\text{minimize}}\qquad\quad\,\,\,\bar{f} \quad \big(\,  \widetilde{\bm x}  \,\big) \,\,\, + \,\, \widehat{g}\,\,\, \big(\, \widetilde{\bm x} - \widetilde{\bm c} \,\big) ,   \tag{primal}\\
	&\underset{\widetilde{\bm\lambda}}{\text{minimize}}\quad    \bar{f}^*\circ(-\mathcal{I})\big(\,\widetilde{\bm\lambda}\,\big)\,\,	+  \,\,  \widehat{g}^*\, \big(\,\widetilde{\bm\lambda} - \widetilde{\bm c}\,\big) ,  \tag{dual}
	\end{align}
	with variable defined as: 
	\begin{equation}
	\widetilde{\bm x}  \eqdef \rho\mathcal{A}\bm{x}  \,\,\in \mathbb{K},
	\qquad  \widetilde{\bm\lambda} \eqdef\bm\lambda/\rho\,\,\in  \mathbb{K},  
	\qquad \widetilde{\bm c} \eqdef \rho\bm c\,\,\in  \mathbb{K},
	\end{equation}
	and function defined as: 
	\begin{align}
	\bar{f}(\bm{y}) &\eqdef  f \circ (\rho\mathcal{A})^{-1}(\bm{y})  -  \langle \rho\bm{c}, \bm{y} \rangle,  \\  
	\widehat{g}(\bm{y}) &\eqdef g \circ   (\rho\mathcal{B})^{-1}\,(\bm{y}) +   \langle \rho\bm{c}, \bm{y} \rangle,  \quad  \forall\bm{y}  \in  \mathbb{K}.
	\end{align}
\end{prop}
\begin{proof}
	Problem \eqref{ab_prob} can be rewritten into
	\begin{align}
	&\underset{\bm x,\bm z}{\text{minimize}}\quad   f(\bm x) + g(\bm z) \qquad\qquad\qquad\qquad\quad\,\, \text{subject\,to}\,\,  \mathcal{A}\bm x - \mathcal{B}\bm z =  \bm c ,\nonumber \\
	\iff
	&\underset{\bm x,\bm z}{\text{minimize}}\quad    f(\rho\mathcal{A})^{-1}(\rho\mathcal{A}\bm x) + g(\rho\mathcal{B})^{-1}(\rho\mathcal{B}\bm z) \,\,\, \text{subject\,to}\,\,  \rho\mathcal{A}\bm x - \rho\mathcal{B}\bm z = \rho\bm c,\nonumber\\
	\iff
	&\underset{\widetilde{\bm x},\widetilde{\bm z}}{\text{minimize}}\quad      f(\rho\mathcal{A})^{-1}(\widetilde{\bm x}) + g(\rho\mathcal{B})^{-1}(\widetilde{\bm z}) - \Vert \widetilde{\bm c} \Vert^2
	\,\text{subject\,to}\,\,   \widetilde{\bm x} - \widetilde{\bm z} = \widetilde{\bm c} ,\nonumber \\
	\iff &\underset{\widetilde{\bm x},\widetilde{\bm z}}{\text{minimize}}\quad     \bar{f}(\widetilde{\bm x}) + \widehat{g}(\widetilde{\bm z}) \qquad\qquad\qquad\qquad\quad\,\,\,
	\text{subject\,to}\,\,   \widetilde{\bm x} - \widetilde{\bm z} = \widetilde{\bm c} . \nonumber
	\end{align} 
	with variable defined as: $ \widetilde{\bm x}=  \rho\mathcal{A}\bm{x},\,  \widetilde{\bm\lambda} = \bm\lambda/\rho,  \,\,  \widetilde{\bm c} = \rho\bm c$, 
	and function defined as: $ 	\bar{f}(\bm{x}) =   f \circ (\rho\mathcal{A})^{-1}(\bm{x})  -  \langle \rho\bm{c}, \bm{x} \rangle, \,\,  \widehat{g}(\bm{x}) =   g \circ (\rho\mathcal{B})^{-1}(\bm{x}) +   \langle \rho\bm{c}, \bm{x} \rangle$.
	
	Next, we derive the  dual problem.  Consider  the following Lagrangian:
	\begin{align}\label{lag01}
	\mathcal{L}(\bm{x},\bm{z},\bm{\lambda}) 
	\,=\,\,	&  f(\bm{x}) + g (\bm{z})  + \langle   \bm{\lambda}, \, \mathcal{A}\bm{x} - \mathcal{B}\bm{z} - \bm{c}  \rangle   ,\nonumber\\
	\,=\,\,	&  f(\bm{x}) + g (\bm{z})  + \big\langle   (\rho)^{-1}  \bm{\lambda}, \, \rho \big( \mathcal{A}\bm{x} - \mathcal{B}\bm{z} - \bm{c}  \big)\big\rangle  ,\nonumber\\
	\,=\,\,	& \underbracket{\,\widetilde{f}(\widetilde{\bm{x}}) + \widetilde g (\widetilde{\bm{z}})  + \langle \widetilde  {\bm\lambda}, \,  \widetilde{\bm{x} }- \widetilde{\bm{z}} - \widetilde{\bm c}\rangle\,}_{
		\eqdef \,\,\, \widetilde{\mathcal{L}}(\widetilde{\bm x}, \widetilde{\bm z}, \widetilde{\bm\lambda})} ,
	\end{align}
	with a dual-variable substitution  $ \widetilde  {\bm\lambda} =  \bm{\lambda} /\rho$.
	Then, the  dual problem is given  by
	\begin{align}\label{new_dual}
	&\underset{\bm\lambda}{\text{maximize}} \,\,\, \underset{{\bm x}, {\bm z}}{\inf}\,\, \mathcal{L}({\bm x}, {\bm z}, \bm\lambda) , \nonumber\\
	=\,\, & \underset{\widetilde{\bm\lambda}}{\text{maximize}} \,\,\, \underset{\widetilde{\bm x}, \widetilde{\bm z}}{\inf}\,\, \widetilde{\mathcal{L}}(\widetilde{\bm x}, \widetilde{\bm z}, \widetilde{\bm\lambda}), \nonumber\\
	=\,\, &  \underset{\widetilde{\bm\lambda}}{\text{maximize}} \,\, 
	- \bigg( \underset{\widetilde{\bm x}}{\sup}\, \biggl\{\langle -\widetilde{\bm\lambda}, \widetilde{\bm x} \rangle - \widetilde{f}(  \widetilde{\bm x} ) \biggr\}  	
	+\underset{\widetilde{\bm z}}{\sup}\, \biggl\{ \langle \widetilde{\bm\lambda}, \widetilde{\bm z} \rangle - \widetilde{g}   ( \widetilde{\bm z}  )\biggr\}
	+ \langle\widetilde{\bm\lambda},\, \widetilde {\bm c} \rangle 		\bigg) ,\nonumber\\
	=\,\, &   \underset{\widetilde{\bm\lambda}}{\text{maximize}} \,\, 
	- \bigg(  \widetilde{f} ^* (-\widetilde{\bm\lambda})	+   \widetilde{g}^*(\widetilde{\bm\lambda}) 	+ \langle\widetilde{\bm\lambda},\, \widetilde{\bm c} \rangle  \bigg),
	\end{align}
	which can be rewritten into
	\begin{align}
	\underset{\widetilde{\bm\lambda}}{\text{minimize}} \,\,\,\,  \widetilde{f} ^* (-\widetilde{\bm\lambda})	+   \widetilde{g}^*(\widetilde{\bm\lambda}) 	+ \langle\widetilde{\bm\lambda},\, \widetilde{\bm c} \rangle  
	=\underset{\widetilde{\bm\lambda}}{\text{minimize}} \,\,\,\,  \widetilde{f} ^* (-\widetilde{\bm\lambda})	+   \widetilde{g}^*(\widetilde{\bm\lambda}-\widetilde{\bm c}) ,
	\end{align}
	which holds  owing  to Lemma \ref{lem01}. The proof is now concluded.
\end{proof}

\subsubsection{parametrized self-duality}\label{self_d}
Recall  the  domain-parametrized  Moreau identity from  Proposition \ref{pmoreau}, restated here as:
\begin{equation}
\bm{v}  =\,\textbf{Prox}_{f\mathcal{D}^{-1}} (\bm{v}) +  \textbf{Prox}_{f^*{\mathcal{D}^*}} (\bm{v}),
\end{equation}
with $ \bm{v} \in  \mathbb K$,	$ f \in \Gamma_0 (\mathbb H) $,    where $ \mathcal D  \in \mathfrak B (\mathbb H,  \mathbb K) $ is  injective.

Let  $ 	\bar{f}(\bm{y}) \eqdef  f \circ (\rho\mathcal{A})^{-1}(\bm{y})  -  \langle \rho\bm{c}, \bm{y} \rangle, \,\, \forall  \bm{y} \in  \mathbb K$.
Then,
\begin{align}
\qquad&  \quad \big(\mathcal I - \textbf{Prox}_{\bar{f}}\big) \,\,(\bm{v}) \,\,=   \quad\,\textbf{Prox}_{\bar{f}^*} \,\,(\bm{v}),\nonumber\\
\iff
\qquad&  \quad \big(\mathcal I - \textbf{Prox}_{\bar{f}}\big) \,\,(\bm{v}) \,\,= \,      - \,\,\textbf{Prox}_{\bar{f}^*(-\mathcal{I})} \,\,(-\bm{v}),\nonumber\\
\iff
\qquad&\,\,\,   \big(2\textbf{Prox}_{\bar{f}} - \mathcal I\big) \,\,(\bm{v}) \,\,=   \quad\big(2\textbf{Prox}_{\bar{f}^*(-\mathcal{I})} - \mathcal I \big)   \,\,(-\bm{v}).
\end{align}

Let $  \widehat{g}(\bm{y}) \eqdef  g \circ (\rho\mathcal{B})^{-1}(\bm{y}) +   \langle \rho\bm{c}, \bm{y} \rangle ,  \,\, \forall  \bm{y} \in  \mathbb K$.
Then,
\begin{align}\label{fix_pd}
&  \quad\,\,\,\, \big(\mathcal I - \textbf{Prox}_{ \widehat{g}}\big) \,\,\,(\bm{v})\,\,=  \quad\,\textbf{Prox}_{ \widehat{g}^*} \,\,(\bm{v}),\nonumber\\
\iff
&  \quad\,\, \big(\mathcal I - 2\textbf{Prox}_{ \widehat{g}  }\big) \,\,(\bm{v}) \,\,=  \quad \big(\textbf{Prox}_{ \widehat{g}^*}  - \textbf{Prox}_{ \widehat{g}}\big) \,\,(\bm{v}),\nonumber\\
\iff
& \quad\,\,  \big(2\textbf{Prox}_{ \widehat{g}} - \mathcal I\big) \,\,(\bm{v}) \,\,=  -\, \big(2\textbf{Prox}_{ \widehat{g}^*} - \mathcal I \big)   \,\,(\bm{v}).
\end{align}

Combining the above two relations, yields
\begin{align}\label{self_dual_eq}
\underbracket{(2\textbf{Prox}_{\bar{f}} - \mathcal{I})\circ(2\textbf{Prox}_{ \widehat{g}} - \mathcal{I})}_{\mathcal{F}_\text{ADMM}}  
= \underbracket{(2\textbf{Prox}_{\bar{f}^*(-\mathcal{I})} - \mathcal{I})\circ(2\textbf{Prox}_{ \widehat{g}^*} - \mathcal{I})}_{\mathcal{F}^*_\text{ADMM}} . 
\end{align}
In  view  of the primal-dual characterization in Proposition \ref{prop_uni_pd}, we see that the above left-hand-side $ \mathcal{F}_\text{ADMM} $ corresponds to the primal problem and  the right-hand-side $ \mathcal{F}^*_\text{ADMM} $ to the dual. 

Then, the equivalence in \eqref{self_dual_eq} implies the same primal and  dual fixed-point sequences, and  hence the same fixed-point, being
\begin{equation}\label{zeta}
\bm{\zeta}^\star \,\,=\,\,   \rho\mathcal{A}\bm{x}^\star+  \frac{1}{\rho} \bm{\lambda}^\star. 
\end{equation}
Moreover, owing to  the above $ \mathcal{F}_\text{ADMM}  = \mathcal{F}^*_\text{ADMM} $, we have shown \eqref{claim}.  It follows that the above expression \eqref{zeta} is  an unscaled  fixed-point.

\section{Adaptive errors \&  unified  convergence rates}\label{sec_3}
In this section, we consider  the   convergence rate issue. 
Below, we first present a  general error characterization, which then leads to a worst-case convergence rate.
Such a rate is later  optimized  to find the optimal step-size.

We will use the  following     definitions:
\begin{defi}\label{def_Lipschitz}
	Let $ \mathcal{F}: \mathbb H \rightarrow \mathbb H$,  let $  L> 0  $. Then,  $ \mathcal{F} $ is
	
	\vspace{4pt}
	\noindent(i) 
	$1/L$-cocoercive if 
	\begin{equation}\label{coco}
	\frac{1}{L}	\Vert \mathcal{F}\bm{x} - \mathcal{F}\bm{y}  \Vert^2 \leq  \langle \mathcal{F}\bm{x} - \mathcal{F}\bm{y}, \bm{x} - \bm{y} \rangle, \quad\forall \bm{x},\bm{y} \in \mathbb H,
	\end{equation}
	Moreover,  for  $ L =1 $,  $ \mathcal{F} $  is called  being firmly non-expansive.

	\vspace{4pt}
	\noindent(ii) $L$-Lipschitz continuous if
	\begin{equation}\label{L_cont}
	\Vert \mathcal{F}\bm{x} - \mathcal{F}\bm{y}  \Vert \leq L \Vert \bm{x} - \bm{y}  \Vert, \quad \forall \bm{x},\bm{y} \in \mathbb H
	\end{equation}
	Moreover,  for  $ L =1 $, the  above $ \mathcal{F} $  is called  being  non-expansive;  for $ L \in (0,1) $, it  is called a contraction.
\end{defi}

\begin{remark}\label{remark_L}
	By Cauchy-Schwarz  inequality, the above $ (i) $  implies (ii).
\end{remark}

\begin{lem}\label{1/2 averaged operator}
	Let $ \mathcal{F}: \mathbb H \rightarrow \mathbb H$ be non-expansive. Then, the averaging $ \frac{1}{2}\mathcal{I} + \frac{1}{2}\mathcal{F} $ is firmly non-expansive. 
\end{lem}

\begin{lem}\cite[Corollary 23.11]{bauschke2017convex}\label{lem2}
	Let $ \mathcal{F}: \mathbb H \rightarrow \mathbb H$ be firmly non-expansive. Then,  $ \mathcal{I} - \mathcal{F} $ is firmly non-expansive. 
\end{lem}

\subsection{General fixed-point  view}

\subsubsection{Adaptive error characterization}
Here, we  provide an adaptive error characterization. We  obtain  two  insights: (i)  if $ L \in (0,1] $, then the error  is guaranteed strictly decreasing at each iteration, until convergence;  (ii) even $ L  \in (1, +\infty)  $  (i.e., potentially  expansive sequence),  it appears still possible  that the error remains decreasing (a convergence condition).  
\begin{prop}[Adaptive error]\label{prop_error}
	Suppose $ \mathcal{F}: \mathbb H \rightarrow \mathbb H$ is $1/L$-cocoercive with  $ L > 0 $. 
	Let process $ \bm{\zeta}^{k+1} = \mathcal{F} \bm{\zeta}^k$ generate a sequence $ \{\bm{\zeta}^k\} $,  with $ \bm{\zeta}^\star \in \text{Fix}\, \mathcal{F} $ denotes the fixed-point (associated with $ \mathcal{F} $), and   $ k = 0,1,\dots $ being  the iteration number  counter.  
	
	Then,  given an arbitrary initialization $ \bm{\zeta}^0 $, at every $ k $-iteration, the error can be characterized into
	\begin{equation}\label{err}
	\underbrace{\Vert \bm{\zeta}^{k+1} - \bm{\zeta}^\star \Vert^2 - \Vert \bm{\zeta}^{k} - \bm{\zeta}^\star \Vert^2}_{\text{error}} 
	\leq 
	- \underbrace{\bigg( \Vert \bm{\zeta}^{k+1} - \bm{\zeta}^k  \Vert^2 + \bigg(\frac{2}{L}-2\bigg)\Vert \bm{\zeta}^{k+1} - \bm{\zeta}^\star \Vert^2 \bigg)}_{\text{gain}} .
	\end{equation}
	
	Moreover, 
	(i) suppose $ L \in (0,1] $. Then,  sequence $ \{\bm{\zeta}^k\} $ is guaranteed to converge to a  fixed-point (if it exists):
	\begin{equation}\label{re_fix}
	\bm{\zeta}^k \rightarrow  \bm{\zeta}^\star\in \text{Fix}\, \mathcal{F}.
	\end{equation}
	(ii) suppose $ L \in (1, +\infty) $. Then,  sequence $ \{\bm{\zeta}^k\} $ is guaranteed to converge to a  fixed-point (if it exists), if at every iteration  $ k $, the following  holds:
	\begin{equation}\label{err_cond}
	\Vert \bm{\zeta}^{k+1} - \bm{\zeta}^k  \Vert^2  >  \bigg(2-\frac{2}{L}\bigg)\Vert \bm{\zeta}^{k+1} - \bm{\zeta}^\star \Vert^2.
	\end{equation}
\end{prop}
\begin{proof}
	First,  we show  expression \eqref{err}.
	By definition of  $1/L$-cocoercive,  recall \eqref{coco}, one has
	\begin{align*}
	& \Vert \mathcal{F}\bm{\zeta}^{k} - \mathcal{F}\bm{\zeta}^\star  \Vert^2 
	\leq L \langle  \mathcal{F}\bm{\zeta}^{k} - \mathcal{F}\bm{\zeta}^\star , \bm{\zeta}^{k} -\bm{\zeta}^\star\rangle, \nonumber\\
	\iff\,\,& \Vert \bm{\zeta}^{k+1} - \bm{\zeta}^\star \Vert^2 
	\,\,\,\,\leq L\langle  \bm{\zeta}^{k+1} - \bm{\zeta}^\star, \bm{\zeta}^{k} -\bm{\zeta}^\star\rangle ,\nonumber\\
	\iff\,\,& 0
	\leq L\langle  \bm{\zeta}^{k+1} - \bm{\zeta}^\star,  \bm{\zeta}^{k} -\bm{\zeta}^{k+1} \rangle + (L-1) \Vert \bm{\zeta}^{k+1} - \bm{\zeta}^\star \Vert^2 ,\nonumber\\
	\iff\,\,& 0
	\leq L\Vert \bm{\zeta}^{k}  - \bm{\zeta}^\star \Vert^2 -L\Vert \bm{\zeta}^{k+1} - \bm{\zeta}^\star \Vert^2 - L\Vert \bm{\zeta}^{k+1} - \bm{\zeta}^k\Vert^2+ \\
	&\qquad 2(L-1) \Vert \bm{\zeta}^{k+1} - \bm{\zeta}^\star \Vert^2   \label{Pythagoras} ,\nonumber\\		
	\iff\,\,& \Vert \bm{\zeta}^{k+1} - \bm{\zeta}^\star \Vert^2 - \Vert \bm{\zeta}^{k} - \bm{\zeta}^\star \Vert^2 
	\leq -\Vert \bm{\zeta}^{k+1} - \bm{\zeta}^k  \Vert^2    -    \frac{2-2L}{L}\Vert \bm{\zeta}^{k+1} - \bm{\zeta}^\star \Vert^2, \nonumber
	\end{align*}
	which gives  \eqref{err}.
	
	For case   (i),  if $ L \in (0,1] $, then  the error will  be strictly decreasing until  convergence. To see this, first, for $ L=1 $, we have $ {2}/{L}-2 =  0 $  and
	the gain is reduced to a single term $ \Vert \bm{\zeta}^{k+1} - \bm{\zeta}^k  \Vert^2 $. 
	Clearly, as long as $  \bm{\zeta}^{k+1} \neq \bm{\zeta}^k  $, the gain  is non-zero  and therefore the error in  \eqref{err} is strictly decreasing. 
	On the other hand, if $  \bm{\zeta}^{k+1} = \bm{\zeta}^k  $, implying convergence, and $ \bm{\zeta}^{k+1} $ is a fixed-point of $ \mathcal{F} $, which is \eqref{re_fix}.
	Now,  consider $ L \in (0,1) $.  In which case, $ {2}/{L}-2 > 0 $ and  the  error  decreases strictly faster, and therefore also converges.
	
	For case   (ii),  essentially, as  long  as the right-hand-side of \eqref{err}  is strictly  negative for all iterations (until  convergence),  then the error on  its left-hand-side must converge to $ 0 $, which corresponds to our last claim in \eqref{err_cond}.
	
	The proof is now concluded.
\end{proof}

\begin{remark}[simplicity \& abstraction]
	Above, we provide a  simple error characterization.	
	The simplicity  owes to the abstraction,  that we avoid discussing the specific structure of the fixed-point operator. 
	A typical structure  is the sum of two monotone operators,
	see Moursi and Vandenberghe \cite{moursi2019douglas}  for a comprehensive investigation and Giselsson and Boyd \cite{giselsson2016linear} for a different  characterization and  applications.
\end{remark}

\subsubsection{Intrinsic almost-linear rate}
Here, we provide two types of rates. Particularly,   we show that  --- the  sub-linear  rate  $ O(1/(k+1)) $ is intrinsically almost  linear,  in the  sense that it will instantly  become linear  with an arbitrarily small improvement on the constant $ L $.
\begin{prop}[Worst-case rates]\label{prop_convergence}
	For the above defined sequence $ \{\bm{\zeta}^k\}  $: 
	
	(i) Suppose  $ L=1 $, then
	\begin{equation}\label{bd1}
	\Vert \bm{\zeta}^{k+1} - \bm{\zeta}^k  \Vert^2 \leq \frac{1}{k+1} \Vert \bm{\zeta}^\star - \bm{\zeta}^0  \Vert^2.  \tag{sub-linear; intrinsic}
	\end{equation}
	
	(ii) Suppose  $ L \in (0,1) $. Let $ \delta  = \,\, \frac{2}{L}-1 $. Then,
	\begin{align}\label{bd2}
	\Vert \bm{\zeta}^{k+1} - \bm{\zeta}^k  \Vert^2 
	\,\,\leq\,\,& \frac{\delta-1}{\delta^{(k+1)} - 1} \Vert \bm{\zeta}^\star - \bm{\zeta}^0  \Vert^2,  \nonumber\\
	\,\,\leq\,\,&  \delta^{-k} \Vert \bm{\zeta}^\star - \bm{\zeta}^0  \Vert^2 . \tag{linear;  strong  assumption}
	\end{align}
\end{prop} 

\begin{proof}
	The proof for the above two cases are similar.

	$\bullet$ To start, we  prove  \eqref{bd2}, which is slightly more general.
	In view of Proposition \ref{prop_error},
	rearrange \eqref{err} into
	\begin{equation}
	\Vert \bm{\zeta}^{t+1} - \bm{\zeta}^t \Vert^2  \leq \Vert \bm{\zeta}^{t} - \bm{\zeta}^\star \Vert^2 - \delta \Vert \bm{\zeta}^{t+1} - \bm{\zeta}^\star \Vert^2 .	
	\end{equation}
	Sum over all $ t $ iterations  till the $ k $-th one:
	\begin{equation}\label{rate_1}
	\sum_{t=0}^{k} \delta^{t} \Vert \bm{\zeta}^{t+1} - \bm{\zeta}^t  \Vert^2 \leq \Vert \bm{\zeta}^\star - \bm{\zeta}^0  \Vert^2. \tag{result 1}
	\end{equation}

	Since $1/L$-cocoercive is always  $L$-Lipschitz continuous, see Remark \ref{remark_L}. Then,  following  from \eqref{L_cont}, we  arrive at
	\begin{align}
	\Vert \bm{\zeta}^{t+1} - \bm{\zeta}^t \Vert \ \leq \Vert \bm{\zeta}^{t} - \bm{\zeta}^{t-1} \Vert 
	\leq  \Vert \bm{\zeta}^{t-2} - \bm{\zeta}^{t-3} \Vert \leq \cdots 
	\end{align}	
	Moreover, since  $ \delta >  0$,  we  can add it  to both  sides  of  the  inequality.
	That  is, for  any $ t \leq k $,  we have 
	\begin{equation}
	\delta^{t} \Vert \bm{\zeta}^{k+1} - \bm{\zeta}^k  \Vert \leq \delta^{t} \Vert \bm{\zeta}^{t+1} - \bm{\zeta}^t \Vert .
	\end{equation}
	Sum  over $ t= 0,1,\dots,k $,  yielding
	\begin{align}
	&\qquad	\sum_{t=0}^{k} \delta^{t} \Vert \bm{\zeta}^{k+1} - \bm{\zeta}^k  \Vert^2 \leq \sum_{t=0}^{k} \delta^{t} \Vert \bm{\zeta}^{t+1} - \bm{\zeta}^t\Vert^2, \label{eq_L1}\\
	\iff&	\frac{1-\delta^{(k+1)}}{1-\delta} \Vert \bm{\zeta}^{k+1} - \bm{\zeta}^k\Vert^2   \leq \sum_{t=0}^{k} \delta^{t}  \Vert \bm{\zeta}^{t+1} - \bm{\zeta}^t\Vert^2 \label{rate_2}, \tag{result 2}
	\end{align}
	%
	
	Adding (\ref{rate_1}) and (\ref{rate_2}) yields the convergence rate bound
	\begin{equation}
	\Vert \bm{\zeta}^{k+1} - \bm{\zeta}^k\Vert^2   \leq    	\frac{\delta-1}{\delta^{(k+1)} - 1} \Vert \bm{\zeta}^\star - \bm{\zeta}^0  \Vert^2. 
	\end{equation}
	Furthermore,  since $  \delta  = \,\, \frac{2}{L}-1  > 1 $, we have 
	\begin{equation}\label{key}
	\frac{\delta-1}{\delta^{(k+1)} - 1}   \leq  \frac{\delta}{\delta^{(k+1)} } \,=\, \delta^{-k},
	\end{equation}	
	which  gives \eqref{bd2}.

	\vspace{10pt}

	$\bullet$ Now, we  prove \eqref{bd1},  which shares similar arguments  with $ \delta = 1 $.
	To avoid repeating, we  directly start from \eqref{rate_1} and \eqref{eq_L1}, which now become
	\begin{align}
	\sum_{t=0}^{k} \Vert \bm{\zeta}^{t+1} - \bm{\zeta}^t  \Vert^2  &\leq \Vert \bm{\zeta}^\star - \bm{\zeta}^0  \Vert^2 ,        \\
	\sum_{t=0}^{k} \Vert \bm{\zeta}^{k+1} - \bm{\zeta}^k  \Vert^2  &\leq \sum_{t=0}^{k} \Vert \bm{\zeta}^{t+1} - \bm{\zeta}^t\Vert^2.
	\end{align}

	Adding the above two yields
	\begin{equation}
	\Vert \bm{\zeta}^{k+1} - \bm{\zeta}^k  \Vert^2 \leq \frac{1}{k+1}\Vert \bm{\zeta}^\star - \bm{\zeta}^0  \Vert^2,
	\end{equation}
	which concludes the
	proof. 
\end{proof}

\begin{remark}[literature connections \& novelty]\label{reamrk_lit}
	The sub-linear result  itself  as  in  \eqref{bd1}  is well-known  in  the  literature, see typical work from   Ryu and  Yin \cite[Theorem 1, sec. 2.4.2]{ryu2022large} and earlier establishment through variational inequality by Bingsheng and Xiaoming \cite{he2015non}.  For the linear  rate, it can be achieved by various different  strong assumptions, see e.g. \cite{giselsson2016linear,moursi2019douglas}.
	
	Our  novelty  lies on at least the  following two  aspects: (i) a new and simple  proof, which  reveals an  almost-linear nature for the  basic case;
	(ii)  For  the  Lipschitz case, rather  than a linear rate $  \delta^{-k} $, the direct consequence of our proof is a stronger one $ (\delta-1)/(\delta^{(k+1)} - 1)$, and  appears admitting   potential to be  further tightened. 
	Calculating  a  stronger factor could  be  beneficial,  including  important theoretical insights and potentially a new approach to  further  improve the convergence rate. 
\end{remark}

\subsection{Primal-dual  sequences: reciprocal rates}
In  the previous section, we characterize the  general convergence behaviour of a fixed-point sequence. 
Recall that the ADMM fixed-point is a combination  of the primal and dual parts, i.e., $ 	\bm{\zeta}^{k+1} 
=  \rho\mathcal{A}\bm{x}^{k+1} +  \bm{\lambda}^{k} / \rho$, see  \eqref{ADMM_fix}.
Here, we investigate the convergence  behaviour of the primal and dual parts separately.
\begin{prop}\label{itr_speed}
	Consider the  ADMM   \eqref{scaled_from} or \eqref{R-ADMM}. 
	Given arbitrary non-zero step-size $\rho \neq  0$ 
	and arbitrary  initialization $  \bm{\zeta}^0 $. Then,
	\begin{align}\label{itr_p}
	\Vert \mathcal{A}\bm{x}^{k+2} - \mathcal{A}\bm{x}^{k+1} \Vert^2   
	\,\,\leq\,\,    &  \frac{1}{k+1}  \cdot \Vert \bm{\zeta}^\star -  \bm{\zeta}^0  \Vert^2 \cdot  \frac{1}{\rho^2}, \nonumber\\
	\,\,=\,\,    &  \frac{1}{k+1}  \cdot \Vert  \mathcal{A}\bm{x}^\star  + \frac{1}{\rho^2}\bm{\lambda}^\star - \frac{1}{\rho} \bm{\zeta}^0  \Vert^2 ,\tag{primal}
	\end{align}
	and
	\begin{align}\label{itr_d}
	\Vert \bm{\lambda}^{k+1} - \bm{\lambda}^{k} \Vert^2   
	\,\,\leq\,\, & \frac{1}{k+1}   \cdot  \Vert \bm{\zeta}^\star -  \bm{\zeta}^0  \Vert^2  \cdot  \rho^2,  \nonumber \\
	\,\,=\,\, &  \frac{1}{k+1}  \cdot \Vert \rho^2\mathcal{A}\bm{x}^\star  + \bm{\lambda}^\star - \rho\bm{\zeta}^0  \Vert^2    \tag{dual}.
	\end{align}
\end{prop}
\begin{proof}
	The proof follows directly from Proposition \ref{fpi}. We will employ the same scaled variables and functions there, restated here as
	\begin{equation}\label{scaled_var2}
	\widetilde{\bm{x}}^{k} \eqdef \rho\mathcal{A}\bm{x}^{k}, \qquad \widetilde{\bm{z}}^{k} \eqdef \rho\mathcal{B}\bm{z}^{k}, \qquad \widetilde{\bm{\lambda}}^{k} \eqdef \frac{1}{\rho}  \bm{\lambda}^{k}
	\qquad \widetilde{\bm{c}} \eqdef\rho \bm{c},
	\end{equation}
	and 
	\begin{align}
	\bar{f}(\bm{y}) &\eqdef  f \circ (\rho\mathcal{A})^{-1}(\bm{y})  -  \langle \rho\bm{c}, \bm{y} \rangle,  \\  
	\widehat{g}(\bm{y}) &\eqdef g \circ   (\rho\mathcal{B})^{-1}\,(\bm{y}) +   \langle \rho\bm{c}, \bm{y} \rangle,  \quad  \forall\bm{y}  \in  \mathbb{K}.
	\end{align}
	
	First, for the scaled primal iterate, we have
	\begin{align}
	\widetilde{\bm{x}}^{k+1}
	=& \textbf{Prox}_{f  (\rho\mathcal{A})^{-1}} (2\widetilde{\bm{c}} +  2 \widetilde{\bm{z}}^{k}  - \bm{\zeta}^{k}) ,\\
	=& \textbf{Prox}_{ {f}(\rho\mathcal{A})^{-1} }\circ \bigg(\,\widetilde{\bm{c}} + \big(2\textbf{Prox}_{ {g}(\rho\mathcal{B})^{-1}} - \mathcal{I}\big)   (\bm{\zeta}^{k} - \widetilde{\bm{c}} )\bigg) ,\\
	=& \underbracket{\textbf{Prox}_{\bar{f}}\circ(2\textbf{Prox}_{ \widehat{g}} - \mathcal{I})}_{\eqdef\,\, \mathcal{F}_x} \,   \bm{\zeta}^{k},
	\end{align}
	The  last line is a composition of a firmly non-expansive operator and  a non-expansive  operator, which  is at  least non-expansive. 
	Hence,
	\begin{align}
	\Vert \widetilde{\bm{x}}^{k+2} - \widetilde{\bm{x}}^{k+1}\Vert^2 
	=  \Vert \mathcal{F}_x   \bm{\zeta}^{k+1} -  \mathcal{F}_x   \bm{\zeta}^{k}  \Vert^2    
	\leq \Vert    \bm{\zeta}^{k+1} -   \bm{\zeta}^{k}  \Vert^2    
	\leq  \frac{1}{k+1}\Vert \bm{\zeta}^\star - \bm{\zeta}^0  \Vert^2  .
	\end{align}
	Invoke definition \eqref{scaled_var2}, then
	\begin{align}\label{term1}
	\Vert \rho\mathcal{A}\bm{x}^{k+2} - \rho\mathcal{A}\bm{x}^{k+1} \Vert^2   
	\leq     \frac{1}{k+1}     \Vert \bm{\zeta}^\star -  \bm{\zeta}^0  \Vert^2 
	=   \frac{1}{k+1}   \Vert \rho\mathcal{A}\bm{x}^\star  + \bm{\lambda}^\star/\rho - \bm{\zeta}^0  \Vert^2 .
	\end{align}
	Since $\rho\neq 0$, we can  divide both sides above with $ \rho^2 $, which gives \eqref{itr_p}.
	
	
	Similarly, for the  dual iterates, we  have
	\begin{align}
	\widetilde{\bm{\lambda}}^{k}
	\,\,=\,\, \bm{\zeta}^{k}  - \widetilde{\bm{z}}^{k} - \widetilde{\bm{c}} 
	\,\,=\,\,  &   (\bm{\zeta}^{k} - \widetilde{\bm{c}}) - \textbf{Prox}_{ g(\rho\mathcal{B})^{-1}}(\bm{\zeta}^{k} - \widetilde{\bm{c}} ) ,  \nonumber\\
	\,\,=\,\,  &   \underbracket{\bigg( \mathcal{I} - \textbf{Prox}_{ \widehat{g}} \bigg)  }_{\eqdef\,\, \mathcal{F}_\lambda} \, \bm{\zeta}^{k} , 
	\end{align}
	where $ \mathcal{F}_\lambda $ is firmly non-expansive (recall Lemma  \ref{lem2}),  which implies non-expansiveness.
	Following the same arguments  as above, we  have
	\begin{equation}\label{term2}
	\Vert 	\widetilde{\bm{\lambda}}^{k+1} -	\widetilde{\bm{\lambda}}^{k} \Vert^2   \leq  \Vert \mathcal{F}_\lambda  \bm{\zeta}^{k+1} -  \mathcal{F}_\lambda   \bm{\zeta}^{k}  \Vert^2   
	\leq \Vert    \bm{\zeta}^{k+1} -   \bm{\zeta}^{k}  \Vert^2   
	\leq \frac{1}{k+1}    \Vert \bm{\zeta}^\star -  \bm{\zeta}^0  \Vert^2  ,
	\end{equation}	
	where $ \widetilde{\bm{\lambda}}^{k} \eqdef  \bm{\lambda}^{k} /\rho$.
	Since $\rho\neq 0$, we can  scale the above relations with $ \rho^2 $, which gives \eqref{itr_d}.
	The proof is now concluded.
\end{proof}

\begin{remark}[Reciprocal residue  behaviours]\label{rema_recipro}
	In Boyd  et. al. \cite[sec. 3.4.1]{boyd12}, the  authors point  out  that: `\textit{The ADMM update equations suggest that large values of step-sizes place a
		large penalty on violations of primal feasibility and so tend to produce small primal residuals.  Conversely, ... }'
	
	Their interpretation via residues do coincide with our  results. 
	In fact, under our measures, the primal and dual admit  exactly reciprocal behaviours:
	\begin{align}
	\Vert \mathcal{A}\bm{x}^{k+2} - \mathcal{A}\bm{x}^{k+1} \Vert^2   
	\,\,\leq\,\,    &  \frac{1}{k+1}  \cdot  \frac{1}{\rho^2}\cdot\Vert \bm{\zeta}^\star -  \bm{\zeta}^0  \Vert^2 , \\
	\Vert \bm{\lambda}^{k+1} - \bm{\lambda}^{k} \Vert^2   
	\,\,\leq\,\, &  \frac{1}{k+1}  \cdot  \rho^2  \cdot  \Vert \bm{\zeta}^\star -  \bm{\zeta}^0  \Vert^2 .
	\end{align}
	From above we see that a larger choice of step-size $\rho$  tends to speed-up the  primal  iterates  convergence, since the right-hand-side  contains  an additional factor $ {1}/{\rho^2} $.  Similar arguments hold for the dual. 
	
	Meanwhile, we emphasize that this observation is only a general trend, i.e., a larger step-size does not guarantee to improve the primal convergence, since there is another factor	
	$ \bm{\zeta}^\star  =  \rho\mathcal{A}\bm{x}^\star +  \bm{\lambda}^\star / \rho$, which also changes with step-size $\rho$.
\end{remark}

\section{Optimal  step-size selection}\label{sec_4}
In this section, we   determine an optimal step-size choice by minimizing the  worst-case convergence  rates established in the previous section. 

Let us note that the step-size selection scheme should be as simple as  possible, since  otherwise there exists  risk that the overall  efficiency     decreases (calculating  step-size choice adds a cost).
Ideally, we prefer a closed-form  choice, which will be  guaranteed by our proposed  method.

\subsection{Towards a simple general principle}
We consider minimizing the worst-case convergence rates  in Proposition \ref{prop_convergence}, restated here  as
\begin{align}
&\Vert \bm{\zeta}^{k+1} - \bm{\zeta}^k  \Vert^2 \leq \quad \frac{1}{k+1} \quad \Vert \bm{\zeta}^\star - \bm{\zeta}^0  \Vert^2, \qquad   L =1, \label{bound0} \\
&\Vert \bm{\zeta}^{k+1} - \bm{\zeta}^k  \Vert^2 \leq \frac{\delta-1}{\delta^{(k+1)} - 1} \Vert \bm{\zeta}^\star - \bm{\zeta}^0  \Vert^2, \qquad \delta  = \,\, \frac{2}{L}-1 , \quad  L \in (0,1) .
\end{align}

Substituting the ADMM-type fixed-point expression, we obtain the following  two step-size selection schemes: 
\begin{align}
\underset{\rho  \neq  0}{\text{minimize}}\quad  &  \Vert \rho\mathcal{A}\bm{x}^\star +  \bm{\lambda}^\star/{\rho} - \bm{\zeta}^0  \Vert^2, \label{opt_bas}   \tag{basic}  \\
\underset{\rho  \neq  0}{\text{minimize}}\quad  &\frac{\delta-1}{\delta^{(k+1)} - 1} \Vert \rho\mathcal{A}\bm{x}^\star +  \bm{\lambda}^\star/\rho - \bm{\zeta}^0  \Vert^2 \label{opt_lip} \tag{Lipschitz},
\end{align}
where term $ 1/ (k+1) $  is omitted in \eqref{opt_bas} (iteration number counter $ k $ is clearly independent here).
To determine which formulation to  employ, let us note that:

$\bullet$ \eqref{opt_bas} is generally applicable to all ADMM-type convex programs, since $ \mathcal{F}_\text{ADMM}  $ is intrinsically firmly non-expansive  ($ L=1 $), without  need of any additional assumptions.

$\bullet$  \eqref{opt_lip} requires calculating factor $\delta = \frac{2}{L}-1$, i.e., the cocoercive constant $ L < 1$, which is not guaranteed  in general, and therefore this formulation is a tailored method.  
Also, this approach is complicated, due to $ L $ is related to the step-size $\rho$, see e.g. \cite{giselsson2016linear,moursi2019douglas} for the specific relation in the sum of two monotone operators problem setting.
Moreover, if the step-size is changing (needed for practical use, see later Sec. \ref{sec_prac}), we then need adaptively calculating  $ L $, which adds additional  cost compared to \eqref{opt_bas}. 
Therefore, if $ L \approx 1 $,  this formulation is not necessary  better (despite some  structure exploited).
Of course, if   $ L \ll 1$, implying highly structured data, we do expect  efficiency advantages.

Overall, since we aim to establish a  simple and general  selection scheme in this work, we  employ \eqref{opt_bas}.  We leave formulation \eqref{opt_lip} to the future work in the context of some highly-structured applications.

\subsubsection{Conventional failure}\label{sec_fail}
Before  solving \eqref{opt_bas} which employs our unscaled fixed-point, 
we first show that the conventional fixed-point  expressions will fail. 
In the literature, there exist  two fixed-point expressions for ADMM:
\begin{equation}\label{fixs}
\underbracket{\bm{\zeta}_1^\star = \mathcal{A}\bm{x}^\star +\bm{\lambda}^\star/\gamma}_{\text{obtained via primal}}, 
\qquad \,\,
\underbracket{\bm{\zeta}_2^\star  = \gamma \mathcal{A}\bm{x}^\star + \bm{\lambda}^\star}_{\text{obtained via dual}}.
\end{equation}
see $ \bm{\zeta}_1^\star $  from \cite[sec. 3.1]{ryu2022large}, corresponding to the primal problem, and $ \bm{\zeta}_2^\star $  from \cite[sec. 3.2]{ryu2022large} (see also \cite{poon2019trajectory}), corresponding to the  dual problem.

Consider $ \bm{\zeta}_1^\star $ first, for comparison purpose, we denote its step-size as $\gamma_1$.
Substitute it to  bound  \eqref{bound0}: 
\begin{align}\label{p1}
&\underset{\gamma_1 >0}{\text{minimize}}\,\,  	\Vert \mathcal{A}\bm{x}^\star +\bm{\lambda}^\star/\gamma_1  - \bm{\zeta}^0\Vert^2 , \\
= \,\,&\underset{\gamma_1 >0}{\text{minimize}}\,\,  \frac{1}{\gamma_1^2}\Vert\bm{\lambda}^\star \Vert^2	+  \frac{2}{\gamma_1}\langle \mathcal{A}\bm{x}^\star - \bm{\zeta}^0,  \bm{\lambda}^\star\rangle + \Vert \mathcal{A}\bm{x}^\star  - \bm{\zeta}^0\Vert^2. 
\end{align}
with solution
\begin{equation}\label{cas1}
\gamma^\star_1
= \begin{cases}
- \frac{\Vert\bm{\lambda}^\star \Vert^2}{\langle \mathcal{A}\bm{x}^\star- \bm{\zeta}^0 , \bm{\lambda}^\star\rangle} &    \langle \mathcal{A}\bm{x}^\star- \bm{\zeta}^0 , \bm{\lambda}^\star\rangle < 0, \\
\downarrow  0 & 	  \langle \mathcal{A}\bm{x}^\star- \bm{\zeta}^0 , \bm{\lambda}^\star\rangle >  0, \\
+\infty & 	\langle \mathcal{A}\bm{x}^\star- \bm{\zeta}^0 , \bm{\lambda}^\star\rangle =0, \tag{primal case}
\end{cases} 
\end{equation}
where $ \downarrow  0 $  denotes approaching $ 0 $  from the positive orthant.

Consider $ \bm{\zeta}_2^\star $. Similarly, denote its corresponding step-size as $\gamma_2$,  we  obtain
\begin{align}\label{p2}
&\underset{\gamma_2 >0}{\text{minimize}}\,\,  	\Vert \gamma_2\mathcal{A}\bm{x}^\star +\bm{\lambda}^\star - \bm{\zeta}^0 \Vert^2  , \\
= \,\,&\underset{\gamma_2 >0}{\text{minimize}}\,\,    \Vert\bm{\lambda}^\star- \bm{\zeta}^0  \Vert^2	+  2{ \gamma_2}\langle \mathcal{A}\bm{x}^\star,  \bm{\lambda}^\star- \bm{\zeta}^0 \rangle +  \gamma_2^2\Vert \mathcal{A}\bm{x}^\star \Vert^2. 
\end{align}
with solution
\begin{equation}\label{cas2}
\gamma_2^\star
= \begin{cases}
- \frac{\langle \mathcal{A}\bm{x}^\star, \bm{\lambda}^\star- \bm{\zeta}^0 \rangle}{\Vert\mathcal{A}\bm{x}^\star \Vert^2} &    \langle \mathcal{A}\bm{x}^\star, \bm{\lambda}^\star- \bm{\zeta}^0 \rangle < 0, \\
+\infty& 	  \langle \mathcal{A}\bm{x}^\star, \bm{\lambda}^\star- \bm{\zeta}^0 \rangle >  0, \\
\downarrow  0  & 	\langle \mathcal{A}\bm{x}^\star, \bm{\lambda}^\star- \bm{\zeta}^0 \rangle =0. \tag{dual case}
\end{cases} 
\end{equation}

We see that $ \gamma^\star_1 $ and $ \gamma^\star_2 $ suggest different choices. 
However,   in  \eqref{fixs} the step-size is the same, i.e., $ \gamma_1 = \gamma_2 $,  see    derivation details from \cite[sec. 3.1, 3.2]{ryu2022large}.

Additionally, for practical  evidence,  consider   the standard  semidefinite programming.
The solution there is often given as a matrix primal-dual pair $ (\bm{X}^\star,  \bm{\Lambda}^\star) $, satisfying $ 	\langle \bm{X}^\star,  \bm{\Lambda}^\star \rangle =0 $ (where $\mathcal{A}= \mathcal{I}$ vanishes).
Then,  under zero initialization $  \bm{\zeta}^0   =\bm 0 $, the above results  indicate  that the step-size  should always be chosen  $  +\infty  $ for the primal problem and $ \downarrow  0  $ for the dual.
Such extreme choices are clearly against practical experiences.

Overall, the  above \eqref{cas1} and \eqref{cas2} should not  be  the correct answers. Their problematic  behaviours could   be  why the optimal step-size issue is  considered open  in the  literature, as pointed out in  \cite[Sec. 1.2, First-order methods]{stellato2020osqp}, \cite[Sec. 8, Parameter selection]{ryu2022large}.

\subsection{General  optimal step-size  (fixed-point  view)}\label{sec_opt_step}
Here, we solve \eqref{opt_bas}.
For better connection  to the current literature,  we emphasize that the  relation  between  
our domain step-size $\rho\neq 0$  and the 
classical   range step-size  $\gamma > 0$ is, recall  Lemma \ref{lem_rel}:
\begin{equation}
\underset{\text{classical range step-size}}{\gamma}   \quad  = \quad \underset{\text{domain step-size}}{\rho^2}.
\end{equation}

Now, invoke \eqref{opt_bas}, we arrive at:
\begin{equation}\label{opt2}
\underset{\rho \neq 0}{\text{minimize}}\quad  \Vert \rho\mathcal{A}\bm{x}^\star +  \bm{\lambda}^\star/{\rho} - \bm{\zeta}^0  \Vert^2,  \tag{opt.  sel.}
\end{equation}
where  $  \bm{\zeta}^0 = \rho_0\mathcal{A}\bm{x}^0 +\bm{\lambda}^0 /\rho_0$ denotes an arbitrary fixed-point initialization, with  $ \bm{x}^0 $, $ \bm{\lambda}^0 $ and $ \rho_0 \neq 0 $ denote the primal, dual, and step-size initializations, respectively.

Expand all $\rho$-related terms:
\begin{equation}\label{opt4}
\underset{\rho \neq 0}{\text{minimize}}\quad  	\rho^2\Vert \mathcal{A}\bm{x}^\star\Vert^2  +  \frac{1}{\rho^2}\Vert\bm{\lambda}^\star\Vert^2 - 2\rho\langle \mathcal{A}\bm{x}^\star, \bm{\zeta}^0  \rangle  
-\frac{2}{\rho}\langle \bm{\lambda}^\star, \bm{\zeta}^0  \rangle.
\end{equation}
The first-order optimality condition  yields
\begin{equation}\label{opt_cond}
2\rho\Vert \mathcal{A}\bm{x}^\star\Vert^2  -  \frac{2}{\rho^3}\Vert\bm{\lambda}^\star\Vert^2 - 2\langle \mathcal{A}\bm{x}^\star, \bm{\zeta}^0  \rangle +\frac{2}{\rho^2}\langle \bm{\lambda}^\star, \bm{\zeta}^0  \rangle = 0 .
\end{equation}
Since $\rho  \neq 0$, the  above  is  equivalent  to solving the following polynomial:
\begin{equation} \label{quartic}
\rho^4\Vert \mathcal{A}\bm{x}^\star\Vert^2  - \rho^3\langle \mathcal{A}\bm{x}^\star, \bm{\zeta}^0  \rangle + \rho\langle \bm{\lambda}^\star, \bm{\zeta}^0  \rangle -  \Vert\bm{\lambda}^\star\Vert^2 = 0, \quad \rho \neq  0, 
\end{equation}
which  is  at most degree-4, and therefore always admits a closed-form solution, specified below.

\subsubsection{closed-form expressions}\label{closed_forms}
For the sake of convenience, we adopt the following norm measures: 
\begin{align}
\Vert \mathcal{A}\bm{x}^\star\Vert =  0 \, &\iff\, \bm{x}^\star =  \bm 0,   \\
\Vert\bm{\lambda}^\star\Vert  = 0\, &\iff\, \bm{\lambda}^\star =  \bm 0 ,
\end{align}
where $ \mathcal{A} $ is  injective.
Now, we specify the closed-form solution:  \\

$\bullet$   (\textbf{Trivial}.)  Suppose $ \Vert \mathcal{A}\bm{x}^\star\Vert =  0,   \Vert\bm{\lambda}^\star\Vert  = 0 $.
Then,
all feasible choices  are equivalent, i.e.,  $ \rho^\star \in [- \infty , +\infty]/\{0\} $
(since the objective in \eqref{opt2} is a constant). \\

$\bullet$   (\textbf{Partly trivial --  useful primal}.) 
Suppose $ \Vert \mathcal{A}\bm{x}^\star\Vert \neq  0,   \Vert\bm{\lambda}^\star\Vert  = 0 $. 
(i) For zero  initialization $ \bm{\zeta}^0 =  \bm  0 $,    the optimal choice   $ \rho^\star \rightarrow 0 $, i.e., infinitely  close  to zero (absolute-value sense);
(ii) For  non-zero  initialization  $ \bm{\zeta}^0 \neq \bm  0 $,   then $ \rho^\star =  { \langle \mathcal{A}\bm{x}^\star, \bm{\zeta}^0  \rangle }/{ \Vert\mathcal{A}\bm{x}^\star\Vert^2} $. \\

$\bullet$   (\textbf{Partly trivial --  useful dual}.) 
Suppose $ \Vert \mathcal{A}\bm{x}^\star\Vert =  0,   \Vert\bm{\lambda}^\star\Vert  \neq 0 $. 
Similar to  above: (i) Set $ \bm{\zeta}^0 =  \bm  0 $,  then   $ \rho^\star \rightarrow \infty $, i.e., infinitely  close  to infinity (absolute-value sense);
(ii) Set $ \bm{\zeta}^0 \neq \bm  0 $,  then $ \rho^\star =  { \Vert\bm{\lambda}^\star\Vert^2}/{ \langle \bm{\lambda}^\star, \bm{\zeta}^0  \rangle } $. \\

$\bullet$   (\textbf{Non-trivial}.)  Suppose  $\, \Vert \mathcal{A}\bm{x}^\star\Vert \neq  0,  \, \Vert\bm{\lambda}^\star\Vert  \neq 0 $.
(i) Set $ \bm{\zeta}^0 =  \bm  0 $,   then   $ \rho^\star = \pm  \sqrt{\Vert\bm{\lambda}^\star\Vert/\Vert\mathcal{A}\bm{x}^\star\Vert}$;
(ii)  Set $ \bm{\zeta}^0 \neq \bm  0 $,   then  the  solution is specified below.

For light of notations,  we rewrite the  polynomial in \eqref{quartic}  as
\begin{equation} 
a \rho^4 + b\rho^3 +  d\rho + e = 0.
\end{equation}
which, in  our  case, admits the following four roots:
\begin{equation}\label{closed_rho2}
\rho = \begin{cases}
\frac{1}{2} ( -\frac{b}{2a} - {u_4} - \sqrt{u_5 - u_6})  , 		\\
\frac{1}{2} ( -\frac{b}{2a} - {u_4} + \sqrt{u_5 - u_6})  ,		\\
\frac{1}{2} ( -\frac{b}{2a} + {u_4} - \sqrt{u_5 + u_6}) , 	\\
\frac{1}{2} ( -\frac{b}{2a} + {u_4} + \sqrt{u_5 + u_6})  ,
\end{cases} 
\end{equation}
where
\begin{equation*}
u_4 = \sqrt{\frac{b^2}{4a^2}+u_3},\quad
u_5 = \frac{b^2}{2a^2} - u_3,\quad
u_6 = -\frac{\frac{b^3}{a^3} + \frac{8d}{a} }{4u_4},
\end{equation*}
and where
\begin{equation*}
u_1 = \frac{\sqrt{27}}{2}(ad^2 + b^2  e),
u_2 = u_1+ \sqrt{(bd - 4ae)^3 + u_1^2},
u_3 = \frac{1}{\sqrt{3}a} (\sqrt[3]{u_2} - \frac{bd - 4ae}{ \sqrt[3]{u_2}} ).	
\end{equation*}
At this  stage, we  do not realize a direct  way to distinguish  the  above 4  roots.  
We would  suggest: first remove complex roots  (since we consider a real Hilbert space; typically,  two real roots  and two complex roots), then find the best one via 
\begin{equation}
\rho^\star = \underset{\rho \in \eqref{closed_rho2}}{\text{argmin}} \,\, \Vert \rho\mathcal{A}\bm{x}^\star +  \bm{\lambda}^\star/{\rho} - \bm{\zeta}^0  \Vert^2,
\end{equation}
which adds  negligible cost  owing to the  closed-form expressions. 

\begin{remark}[complex variables]
	For  some applications, the primal-dual solutions may be  complex. In which  case, one may simply solve:
	\begin{equation} 
	\rho^4\Vert \mathcal{A}\bm{x}^\star\Vert^2  - \rho^3\langle \mathcal{A}\bm{x}^\star, \bm{\zeta}^0 \rangle_\mathbb{R} + \rho\langle \bm{\lambda}^\star, \bm{\zeta}^0 \rangle_\mathbb{R} -  \Vert\bm{\lambda}^\star\Vert^2 = 0,
	\end{equation}
	where $ \langle\cdot,\cdot\rangle_\mathbb{R} $ denotes the real part of an inner product.
\end{remark}

\subsection{Optimal rates}

\subsubsection{general (zero initialization): balanced  rates}\label{sec_itr_conv1}
Below we will see that,
under  zero  initialization, employ the optimal step-size. Then,  the worst-case convergence rates of  the  normalized   primal, dual and  fixed-point sequences are the same.
\begin{prop}\label{pro_non}
	Suppose  $ \Vert \mathcal{A}\bm{x}^\star\Vert \neq  0,  \, \Vert\bm{\lambda}^\star\Vert  \neq 0 $.  Set $ \bm{\zeta}^0  =  \bm{0} $ and  choose  step-size $   \rho^\star = \pm\sqrt{{\Vert\bm{\lambda}^\star\Vert}/{\Vert\mathcal{A}\bm{x}^\star\Vert}}$.
	Let  
	\begin{equation}\label{in_angle}
	\theta  =  \arccos  \,  \frac{ \langle \mathcal{A}\bm{x}^\star, \bm{\lambda}^\star\rangle }{\Vert\mathcal{A}\bm{x}^\star \Vert \Vert \bm{\lambda}^\star \Vert}.  \tag{intrinsic angle}
	\end{equation}
	Then, 
	\begin{align}
	\frac{\Vert \mathcal{A}\bm{x}^{k+2} - \mathcal{A}\bm{x}^{k+1} \Vert^2   }{\Vert\mathcal{A}\bm{x}^\star \Vert^2 }
	\,\,\leq\,\,  &  \,\,\frac{2}{k+1}   \,   (1  + \cos \theta ), \label{cc1}   \tag{primal}\\
	\frac{\Vert \bm{\lambda}^{k+1} - \bm{\lambda}^{k} \Vert^2   }{\Vert\bm{\lambda}^\star \Vert^2 }
	\quad\,\, \leq\,\, &  \,\,\frac{2}{k+1}   \,  (1  + \cos \theta ),
	\label{cc2}    \tag{dual}\\
	\frac{\Vert \bm{\zeta}^{k+1} - \bm{\zeta}^{k} \Vert^2   }{\Vert\mathcal{A}\bm{x}^\star \Vert\Vert\bm{\lambda}^\star \Vert }
	\quad\,\, \leq\,\, &  \,\,\frac{2}{k+1}   \,  (1  + \cos \theta ). \tag{fixed-point} 
	\end{align}
\end{prop}
\begin{proof}
	Following from   Proposition \ref{itr_speed}, given   $ \bm{\zeta}^0   =  \bm{0} $ and  the   optimal step-size $\rho^\star$, 
	we  arrive  at
	\begin{align}
	\Vert \mathcal{A}\bm{x}^{k+2} - \mathcal{A}\bm{x}^{k+1} \Vert^2   
	\,\,\leq\,\,   &   \frac{1}{k+1}   \,  \Vert  \mathcal{A}\bm{x}^\star  + \frac{1}{(\rho^\star)^2}\bm{\lambda}^\star - \bm{0} \Vert^2 . \\
	\Vert \bm{\lambda}^{k+1} - \bm{\lambda}^{k} \Vert^2   
	\,\,\leq\,\, &  \frac{1}{k+1}  \,  \Vert (\rho^\star)^2\mathcal{A}\bm{x}^\star \, + \bm{\lambda}^\star - \bm{0}   \Vert^2 ,\\
	\Vert \bm{\zeta}^{k+1} - \bm{\zeta}^k  \Vert^2 
	\,\,\leq\,\,   &  \frac{1}{k+1} \Vert \rho^\star\mathcal{A}\bm{x}^\star \, + \bm{\lambda}^\star/\rho^\star - \bm{0} \Vert^2  .
	\end{align}
	Invoke  $   \rho^\star = \pm\sqrt{{\Vert\bm{\lambda}^\star\Vert}/{\Vert\mathcal{A}\bm{x}^\star\Vert}}$ and $ \theta $ defined in \eqref{in_angle},
	and rearrange the terms   concludes  the  proof.
\end{proof}

\begin{remark}[intrinsic   efficiency]
	Given  a certain convex program, let  us note that its  optimal  solution  is   fixed, hence $ \theta $  is   fixed.	
	Then, from above, we see that   the efficiency  of  different convex  programs (under zero initialization) can be  characterized  by $ \cos  \theta \in [-1,1]$. 
\end{remark}

\begin{remark}[Successor  to  the literature adaptive principle?]
	As pointed  out in Boyd et. al. \cite[sec. 3.4.1]{boyd12},   a simple adaptive step-size selection scheme that often works well is (e.g. \cite{he2000alternating,wang2001decomposition}):
	\begin{equation}
	\gamma^{k+1}  = \begin{cases}
	\tau^{\text{incr}}\,\gamma^{k} &   \text{if} \,\, \Vert \bm{r}^k \Vert  > \mu \Vert \bm{s}^k \Vert, \\
	\gamma^{k}/\tau^{\text{decr}} & 	 \text{if}  \,\,  \Vert \bm{r}^k \Vert < \mu \Vert \bm{s}^k \Vert, \\
	\gamma^{k} & 	\text{otherwise}, 
	\end{cases}  
	\end{equation}
	where  $ \bm{r}^k $ and $ \bm{s}^k $  are primal and  dual residues,  respectively.  Typical choices
	might be $ \mu = 10 $, $ \tau^{\text{incr}}  =  \tau^{\text{decr}} = 2 $.
	`The idea behind this penalty
	parameter update is to try to keep the primal and dual residual norms
	within a factor of $ \mu $ of one another as they both converge to zero.' \cite{boyd12}

	Such  an idea appears consistent here --- the  optimal step-size is  the one that  makes the (normalized) primal and dual iterates converge in the same speed. Compared to the above heuristic adaptive scheme, our    step-size  choice controls the primal and  dual  sequence to be the same speed, rather than  theirs within factor $ \mu $ (typically $ \mu = 10 $). In this sense, ours seems a better choice.
\end{remark}

\begin{remark}[single stopping threshold]
	Owing  to the balanced rates,  we  may  simplify  the stopping criteria for ADMM into a single error threshold (under zero initialization). This is a simplification to  the multiple thresholds in the current literature, see  \cite{boyd12}. 
\end{remark}

\subsubsection{partly trivial (zero  initialization):  fixed non-trivial rate} \label{sec_special}
Below, we show  that if either the primal or dual solution is a zero vector (trivial), then the convergence rate of the
non-trivial iterate is independent of the step-size selection.
\begin{prop}\label{prop_zero}
	(i)  Suppose  $ \Vert \mathcal{A}\bm{x}^\star\Vert \neq  0,  \, \Vert\bm{\lambda}^\star\Vert  = 0 $. Set $ \bm{\zeta}^0  =  \bm{0} $.
	Then,
	\begin{align}
	\qquad	{\Vert \mathcal{A}\bm{x}^{k+2} - \mathcal{A}\bm{x}^{k+1} \Vert^2  }
	\,\,\leq\,\,\,\,    & \frac{1}{k+1} \cdot \Vert  \mathcal{A}\bm{x}^\star \Vert^2, \tag{primal} \label{01} \\
	\qquad	{\Vert \bm{\lambda}^{k+1} - \bm{\lambda}^{k} \Vert^2   }
	\,\, \leq  \,\,\,\,  &  \frac{1}{k+1} \cdot  \Vert\mathcal{A}\bm{x}^\star\Vert^2   \cdot  \rho^4     ,  \tag{trivial dual} \label{02}  \\
	\qquad		\Vert \bm{\zeta}^{k+1} - \bm{\zeta}^k  \Vert^2 
	\,\, \leq  \,\,\,\,   &  \frac{1}{k+1} \cdot  \Vert\mathcal{A}\bm{x}^\star\Vert^2  \cdot  \rho^2   .  \tag{fixed-point} 
	\end{align}

	(ii)  Suppose  $ \Vert \mathcal{A}\bm{x}^\star\Vert =  0,  \, \Vert\bm{\lambda}^\star\Vert  \neq 0$. Set $ \bm{\zeta}^0  =  \bm{0} $.
	Then,
	\begin{align}
	{\Vert \bm{\lambda}^{k+1} - \bm{\lambda}^{k} \Vert^2 } 
	\,\, \leq  \,\,\,\,  & \frac{1}{k+1}\cdot \Vert \bm{\lambda}^{\star}  \Vert^2, \tag{dual} \label{d01} \\
	\Vert \mathcal{A}\bm{x}^{k+2} - \mathcal{A}\bm{x}^{k+1} \Vert^2 
	\,\, \leq  \,\,\,\,  & \frac{1}{k+1} \cdot  \Vert\bm{\lambda}^\star\Vert^2 \cdot \frac{ 1}{\rho^4}   ,  \tag{trivial primal}\\
	\Vert \bm{\zeta}^{k+1} - \bm{\zeta}^k  \Vert^2 
	\,\, \leq  \,\,\,\,  & \frac{1}{k+1} \cdot  \Vert\bm{\lambda}^\star\Vert^2 \cdot \frac{ 1}{\rho^2}    .  \tag{fixed-point} 
	\end{align}
\end{prop}
\begin{proof}
	The proof follows  straightforwardly from Proposition \ref{itr_speed}.
\end{proof}

\begin{remark}[independence \& optimality interpretation]
	Interestingly, under zero initialization, the above result implies that the   step-size selection is independent of the  non-trivial sequence worst-case rate.
	This will be verified numerically  in Sec. \ref{feasi}, which shows a complete independence in  practice.
	
	Moreover,  
	we see that the previously derived   step-size choice  is optimal, in the sense  that  the trivial sequence convergence rate is most accelerated.
	Specifically,  for the above case (i), recall the optimal step-size  being $\rho^\star \rightarrow 0$.
	Substituting  it to   the above,  only the dual and fixed-point sequence  one is affected, with  their  upper  bounds  being:
	\begin{equation}
	\underset{\rho^\star \rightarrow 0}{\lim} \,  (\rho^\star)^4\cdot \frac{  \Vert\mathcal{A}\bm{x}^\star\Vert^2 }{k+1} =  0, 
	\quad\,\,	
	\underset{\rho^\star \rightarrow 0}{\lim} \,  (\rho^\star)^2\cdot \frac{  \Vert\mathcal{A}\bm{x}^\star\Vert^2 }{k+1}  =  0.
	\end{equation}
	Similar arguments hold for case (ii). 
\end{remark}

\subsubsection{partly trivial (non-zero  initialization): initialization angle}\label{sec_tri}
Here,  we consider an arbitrary non-zero initialization, where  the convergence rate of the non-trivial sequence  is now related to the step-size choices and can be improved compared to the previous fixed  one in Proposition \ref{prop_zero}. 
The maximum  improvement depends on the  initialization angle only. 
\begin{prop}[balanced/fixed-point view]\label{pro0}
	(i) Suppose  $ \Vert \mathcal{A}\bm{x}^\star\Vert \neq  0,  \, \Vert\bm{\lambda}^\star\Vert  = 0$.  Set  any  initialization $ \bm{\zeta}^0  \neq \bm{0} $.  
	Let 
	\begin{equation}
	\omega_1  =  \arccos  \, \frac{\langle \mathcal{A}\bm{x}^\star, \bm{\zeta}^0  \rangle }{\Vert\mathcal{A}\bm{x}^\star\Vert\Vert \bm{\zeta}^0\Vert}. \tag{initial primal angle}
	\end{equation}
	Choose  step-size  to be $ \rho^\star=  \frac{\Vert\bm{\zeta}^0 \Vert}{\Vert\mathcal{A}\bm{x}^\star\Vert} \cos \omega_1 $.
	Then,
	\begin{align}
	{\Vert \mathcal{A}\bm{x}^{k+2} - \mathcal{A}\bm{x}^{k+1} \Vert^2  }
	\,\, \leq  \,\,\,\,  &  \frac{1}{k+1} \cdot \Vert  \mathcal{A}\bm{x}^\star \Vert^2 \cdot \tan^2  \omega_1,    \tag{non-trivial}     \\
	{\Vert \bm{\lambda}^{k+1} - \bm{\lambda}^{k} \Vert^2  }
	\,\, \leq  \,\,\,\,  &  \frac{1}{k+1} \cdot  \frac{ \Vert  \bm{\zeta}^0\Vert^4}{ 4\Vert \mathcal{A}\bm{x}^\star \Vert^2}\cdot \sin ^2  2\omega_1  .  \tag{trivial dual} \\
	{\Vert \bm{\zeta}^{k+1} - \bm{\zeta}^k  \Vert^2}
	\,\, \leq  \,\,\,\,  & \frac{1}{k+1} \cdot   {\Vert\bm{\zeta}^0\Vert^2}\cdot  \sin^2  \omega_1 .  \tag{fixed-point} 
	\end{align}

	(ii) Suppose  $ \Vert \mathcal{A}\bm{x}^\star\Vert =  0,  \, \Vert\bm{\lambda}^\star\Vert  \neq 0$.
	Set  any  initialization $ \bm{\zeta}^0  \neq \bm{0} $. 
	Let 
	\begin{equation}
	\omega_2  =  \arccos  \, \frac{ \langle \bm{\lambda}^\star, \bm{\zeta}^0  \rangle }{\Vert\bm{\lambda}^\star\Vert\Vert \bm{\zeta}^0\Vert}. \tag{initial dual  angle}
	\end{equation}
	Choose  step-size  to be $ \rho^\star= 1/ \big( \frac{\Vert\bm{\zeta}^0 \Vert}{\Vert\bm{\lambda}^\star\Vert} \cos \omega_2 \big) $.
	Then,
	\begin{align}
	{\Vert \bm{\lambda}^{k+1} - \bm{\lambda}^{k} \Vert^2 }
	\,\, \leq  \,\,\,\,  &     \frac{1}{k+1}  \cdot   \Vert \bm{\lambda}^{\star}  \Vert^2 \cdot\tan^2  \omega_2			\tag{non-trivial} ,\\
	\Vert \mathcal{A}\bm{x}^{k+2} - \mathcal{A}\bm{x}^{k+1} \Vert^2 
	\,\, \leq  \,\,\,\,  &    \frac{1}{k+1} \cdot   \frac{ \Vert  \bm{\zeta}^0\Vert^4}{ 4 \Vert \bm{\lambda}^\star \Vert^2}  \cdot\sin ^2  2\omega_2 .  \tag{trivial primal}\\
	\Vert \bm{\zeta}^{k+1} - \bm{\zeta}^k  \Vert^2 
	\,\, \leq  \,\,\,\,  & \frac{1}{k+1} \cdot \Vert  \bm{\zeta}^0\Vert^2 \cdot   \sin^2  \omega_2 .  \tag{fixed-point} 
	\end{align}
\end{prop}
\begin{proof}
	For case  (i),  Proposition \ref{itr_speed}  reduces to
	\begin{align}
	\Vert \mathcal{A}\bm{x}^{k+2} - \mathcal{A}\bm{x}^{k+1} \Vert^2   
	\,\,\leq\,\,    &  \frac{1}{k+1}   \,  \Vert  \mathcal{A}\bm{x}^\star  + \frac{1}{\rho^2}\bm  0 - \frac{1}{\rho}\, \bm{\zeta}^0  \Vert^2 .
	\end{align}
	Invoke  the step-size choice $ \rho^\star $, which yields
	\begin{align}
	\Vert  \mathcal{A}\bm{x}^\star  \Vert^2   -  \frac{2}{\rho^\star}   \langle \mathcal{A}\bm{x}^\star, \bm{\zeta}^0  \rangle + \frac{1}{(\rho^\star)^2} \Vert  \bm{\zeta}^0\Vert^2
	=\, &  -\Vert  \mathcal{A}\bm{x}^\star  \Vert^2   +    \frac{\Vert  \mathcal{A}\bm{x}^\star  \Vert^2 }{\cos^2  \omega} ,\\ \nonumber
	=\, & \Vert  \mathcal{A}\bm{x}^\star  \Vert^2   \tan^2  \omega_1 .
	\end{align}
	For the dual iterates,   note that  the primal and dual sequence rates are related via a factor  of $ \rho^4 $, recall Proposition \ref{itr_speed}. Hence,
	\begin{align}
	\Vert \bm{\lambda}^{k+1} - \bm{\lambda}^{k} \Vert^2   
	\,\,\leq\,\, &  \frac{1}{k+1}  \,\cdot (\rho^\star)^4  \Vert  \mathcal{A}\bm{x}^\star  \Vert^2  \tan ^2  \omega ,\\
	\,\,=\,\, &  \frac{1}{k+1}  \,\cdot  \frac{\Vert\bm{\zeta}^0 \Vert^4}{\Vert\mathcal{A}\bm{x}^\star\Vert^4} \cos^4 \omega_1 \cdot \Vert  \mathcal{A}\bm{x}^\star  \Vert^2 \cdot\frac{\sin ^2  \omega_1 }{\cos ^2  \omega_1 }, \\
	\,\,=\,\, &  \frac{1}{k+1}  \,\cdot  \frac{\Vert\bm{\zeta}^0 \Vert^4}{\Vert\mathcal{A}\bm{x}^\star\Vert^2} \cdot \sin ^2  \omega_1 \cos ^2  \omega_1 .
	\end{align}
	Invoke relation $ \sin  2\omega_1  = 2 \sin\omega_1\cos\omega_1$ yields the result.
	Similarly,  the fixed-point admits  a  factor  of $ \rho^2 $ relation  to  the primal, which gives
	\begin{align}
	\Vert \bm{\zeta}^{k+1} - \bm{\zeta}^k  \Vert^2 
	\,\,\leq\,\, &  \frac{1}{k+1}  \,\cdot (\rho^\star)^2  \Vert  \mathcal{A}\bm{x}^\star  \Vert^2  \tan ^2  \omega ,\\
	\,\,=\,\, &  \frac{1}{k+1}  \,\cdot  \frac{\Vert\bm{\zeta}^0 \Vert^2}{\Vert\mathcal{A}\bm{x}^\star\Vert^2} \cos^2 \omega_1 \cdot \Vert  \mathcal{A}\bm{x}^\star  \Vert^2 \cdot\frac{\sin ^2  \omega_1 }{\cos ^2  \omega_1 } ,\\
	\,\,=\,\, &  \frac{1}{k+1}  \,\cdot  \Vert\bm{\zeta}^0 \Vert^2 \cdot\sin ^2  \omega_1.
	\end{align}
	
	At last, case  (ii)  is a  symmetric  situation  and  the  same  arguments  apply, which  are  omitted  to avoid repeating.
	The  proof  is   now concluded.
\end{proof}

\begin{remark}[good  angle range: monotonic improvements]
	Recall from the previous section that the non-trivial sequence admits a fixed $ 1/(k+1) $  rate. Here,  we see that  such a  rate  is  improved  if
	\begin{equation}
	\tan^2  \omega< 1,
	\end{equation}
	which  corresponds  to $ \omega \in  (-\pi/4, \pi/4) \, \bigcup\, (3\pi/4, 5\pi/4)  $. 
	Meanwhile, suppose $ \omega $ is moving from (i) $ -\pi/4 $ towards $ 0 $, or (ii) from $ \pi/4  $  towards $ 0$, or (iii) from $ 3\pi/4 $ towards $ \pi $, or (iv) from $ 5\pi/4 $   towards $ \pi $. 
	
	Then,  we  have  $ \tan^2  \omega $, $ \sin^2  2\omega  $, and $ \sin^2  \omega  $ monotonically decreasing, implying the rate bounds in Proposition \ref{pro0} are monotonically improving.
\end{remark}

The above result is optimal for the fixed-point sequence, which  is a balance of  primal and dual rates. 
However, for the partly trivial settings here, it  appears  that we do not  need  a balance. 
Instead, it seems  natural  to promote the non-trivial one as much as possible, meanwhile completely ignore the  rest. This leads to the result below.
\begin{prop}[non-trivial sequence only]\label{pro_part}
	(i) Suppose  $ \Vert \mathcal{A}\bm{x}^\star\Vert \neq  0$,  $\Vert\bm{\lambda}^\star\Vert  = 0$.  Set  any  initialization $ \bm{\zeta}^0  \neq \bm{0} $.  
	Let 
	\begin{equation}
	\omega_1  =  \arccos  \, \frac{\langle \mathcal{A}\bm{x}^\star, \bm{\zeta}^0  \rangle }{\Vert\mathcal{A}\bm{x}^\star\Vert\Vert \bm{\zeta}^0\Vert}. \tag{initial  angle}
	\end{equation}
	Choose  step-size  to be $ \rho^\star_\text{pri} =  \frac{\Vert\bm{\zeta}^0 \Vert}{\Vert\mathcal{A}\bm{x}^\star\Vert} \frac{1}{\cos \omega_1} $.
	Then,
	\begin{align}\label{re01}
	{\Vert \mathcal{A}\bm{x}^{k+2} - \mathcal{A}\bm{x}^{k+1} \Vert^2  }
	\,\, \leq  \,\,\,\,  &     \frac{1}{k+1} \cdot \Vert  \mathcal{A}\bm{x}^\star \Vert^2   \cdot \sin^2 \omega_1 ,  \tag{non-trivial}   \\
	{\Vert \bm{\lambda}^{k+1} - \bm{\lambda}^{k} \Vert^2   }
	\,\, \leq  \,\,\,\,  &  \frac{1}{k+1} \cdot  \frac{ \Vert  \bm{\zeta}^0\Vert^4}{ \Vert \mathcal{A}\bm{x}^\star \Vert^2}  \cdot \frac{ \tan^2  \omega_1 }{\cos^2 \omega_1},  \tag{trivial dual} \\
	\Vert \bm{\zeta}^{k+1} - \bm{\zeta}^k  \Vert^2 
	\,\, \leq  \,\,\,\,  & \frac{1}{k+1} \cdot \Vert  \bm{\zeta}^0\Vert^2 \cdot { \tan^2  \omega_1 }.  \tag{fixed-point} 
	\end{align}
	
	(ii) Suppose  $ \Vert \mathcal{A}\bm{x}^\star\Vert =  0,  \, \Vert\bm{\lambda}^\star\Vert  \neq 0$.
	Set  any  initialization $ \bm{\zeta}^0  \neq \bm{0} $. 
	Let 
	\begin{equation}
	\omega_2  =  \arccos  \, \frac{ \langle \bm{\lambda}^\star, \bm{\zeta}^0  \rangle }{\Vert\bm{\lambda}^\star\Vert\Vert \bm{\zeta}^0\Vert}. \tag{initial  angle}
	\end{equation}
	Choose  step-size  to be $ \rho^\star_\text{dual} = 1/ \big( \frac{\Vert\bm{\zeta}^0 \Vert}{\Vert\bm{\lambda}^\star\Vert} \frac{1}{\cos \omega_2} \big) $.
	Then,
	\begin{align}\label{re02}
	{\Vert \bm{\lambda}^{k+1} - \bm{\lambda}^{k} \Vert^2 }
	\,\, \leq  \,\,\,\,  &    \frac{1}{k+1} \cdot {\Vert \bm{\lambda}^{\star}  \Vert^2 } \cdot \sin^2 \omega_2 	\tag{non-trivial} ,\\
	\Vert \mathcal{A}\bm{x}^{k+2} - \mathcal{A}\bm{x}^{k+1} \Vert^2 
	\,\, \leq  \,\,\,\,  &    \frac{1}{k+1} \cdot   \frac{ \Vert  \bm{\zeta}^0\Vert^4}{  \Vert \bm{\lambda}^\star \Vert^2}\cdot \frac{ \tan^2  \omega_2 }{\cos^2 \omega_2},  \tag{trivial primal}\\
	\Vert \bm{\zeta}^{k+1} - \bm{\zeta}^k  \Vert^2 
	\,\, \leq  \,\,\,\,  & \frac{1}{k+1} \cdot \Vert  \bm{\zeta}^0\Vert^2\cdot  \tan^2  \omega_2 .  \tag{fixed-point} 
	\end{align}
\end{prop}
\begin{proof}
	For case  (i),  Proposition \ref{itr_speed}  reduces to
	\begin{align}
	\Vert \mathcal{A}\bm{x}^{k+2} - \mathcal{A}\bm{x}^{k+1} \Vert^2   
	\,\,\leq\,\,    &  \frac{1}{k+1}   \,  \Vert  \mathcal{A}\bm{x}^\star  + \frac{1}{\rho^2}\bm  0 - \frac{1}{\rho}\, \bm{\zeta}^0  \Vert^2 .
	\end{align}
	Hence, the  step-size that optimize the  above bound   is  
	\begin{align}
	\rho^\star_\text{pri}  = \underset{\rho \neq 0}{\text{argmin}}  \,\, \Vert \mathcal{A}\bm{x}^\star  + \bm 0  - \frac{1}{\rho}\, \bm{\zeta}^0 \Vert^2 
	=  { \Vert \bm{\zeta}^0 \Vert^2}/{ \langle \mathcal{A}\bm{x}^\star, \bm{\zeta}^0  \rangle }
	=   \Vert \bm{\zeta}^0 \Vert/ (\Vert \mathcal{A}\bm{x}^\star \Vert \cos \omega_1). 
	\end{align}
	The rest  arguments are   identical to Proposition \ref{pro0}, and are omitted to avoid repeating.   The proof is now concluded.
\end{proof}

\begin{remark}[comparison: improvement]
	Compare the above  to the previous fixed-point view as  in  Proposition \ref{pro0}. We see that the non-trivial sequence convergence  rate is guaranteed improved:
	\begin{equation}
	\tan^2  \omega  \in [0, +\infty]  \,\,\geq\,\,     \sin^2 \omega \in [0, 1].
	\end{equation}
	Moreover, the  upper  rate bound is now guaranteed no worse than $ 1/(k+1) $, i.e., no worse than the previous  zero-initialization case.  The  worst case  happens if and only if  either $ \omega = \pi/2 $, or $ \omega = 3\pi/2 $, i.e., when the initialization  angle happens to be orthogonal to the ground-truth.	
\end{remark}

\begin{remark}[cost]
	The above improvement comes with a cost --- an additional factor $ 1/\cos^6 \omega $ for the trivial iterate sequence, and factor $ 1/\cos^2 \omega $ for the fixed-point sequence. 
\end{remark}

\begin{remark}[optimal initialization preview]
	For all rate results in this section, the upper bounds take value $ 0 $  when $ \omega = 0 $ or $ \omega =  \pi  $. This is the optimal initialization, and  will be  discussed below in  a more general context.
\end{remark}

\subsection{Warm start: quality-based initial step-sizes}\label{sec_opt_ini}
Across all  previous sections, we  see that the  initialization makes  a significant contribution to the convergence rate, meanwhile can be arbitrarily chosen.
Then, a natural  question is  --- what  is the optimal choice?

This question  needs to be answered separately: (i) for  the primal  and dual iterates, there is no extra freedom, and the optimal initialization is the optimal  solution itself, i.e., a unique choice $ \bm{x}^0  =  \bm{x}^\star, \, \bm{\lambda}^0 =\bm{\lambda}^\star  $;  (ii)  for the fixed-point, there exists  a perturbation freedom, i.e., any  element from  the  following  set is optimal:
\begin{equation}\label{fix_01}
\bm{\zeta}^0_{\text{opt}} \in  \big\{\rho_{0}\mathcal{A}\bm{x}^\star + \bm{\lambda}^\star/ \rho_{0}\,|\,  \rho_{0} \in  [-\infty, +\infty] /\{0\}   \big\},
\end{equation} 
where $ \rho_{0} $  is the initial step-size.

Despite all  $ \rho_{0}  \in  [-\infty, +\infty] /\{0\}$ are equivalent in the above optimal case, they do make a difference  in practice regarding the warm-start (not exactly optimal).
Below,  we  provide a rough  estimation to guide the initial step-size selection. In case that an estimation is not possible, we set $ \rho_{0} = 1$.
\begin{prop}[warm-start  initial step-size]\label{prop_warm}
	Given a warm start as  
	\begin{equation}
	\bm{x}^0  = \bm{x}^\star + \epsilon_{\text{err}1}, \quad\bm{\lambda}^0  = \bm{\lambda}^\star + \epsilon_{\text{err}2}.
	\end{equation}
	Let 
	\begin{equation}
	\eta  =  \arccos  \, \frac{\langle \mathcal{A}\epsilon_{\text{err}1}, \epsilon_{\text{err}2} \rangle}{\Vert\mathcal{A}\epsilon_{\text{err}1}\Vert\Vert\epsilon_{\text{err}2}\Vert}.
	\end{equation} 
	Fix the step-size to be:
	\begin{equation}\label{re_ini}
	\qquad      \rho_k = \rho_0  = \pm \sqrt{ \frac{\Vert \epsilon_{\text{err}2} \Vert}{ \Vert\mathcal{A}\epsilon_{\text{err}1}\Vert } }, \qquad  \forall k =1, 2 \dots
	\end{equation}
	Then, the following hold:
	\begin{align}
	\frac{\Vert \mathcal{A}\bm{x}^{k+2} - \mathcal{A}\bm{x}^{k+1} \Vert^2  }{\Vert\mathcal{A}\epsilon_{\text{err}1}\Vert^2 }  
	\,\leq\,\,    &  \frac{2}{k+1}   \,   (1  + \cos \eta ), \label{pri_s1}\\
	\frac{\Vert \bm{\lambda}^{k+1} - \bm{\lambda}^{k} \Vert^2 }{\Vert\epsilon_{\text{err}2}\Vert^2 } 
	\quad\leq\,\, & \frac{2}{k+1}   \,   (1  + \cos \eta ),\label{dual_s1} \\
	\frac{\Vert \bm{\zeta}^{k+1} - \bm{\zeta}^{k} \Vert^2  }{\Vert\mathcal{A}\epsilon_{\text{err}1}\Vert \Vert\epsilon_{\text{err}2}\Vert }  
	\quad\leq\,\,   &  \frac{2}{k+1}   \,   (1  + \cos \eta ). \label{fix_s1} 
	\end{align}
\end{prop}
\begin{proof}
	The proof follows straightforwardly from  Proposition \ref{itr_speed}. 
\end{proof}
\begin{remark}[quality-based initial choice]
	In view  of \eqref{re_ini}, it is a norm  ratio, and  hence possible to  be roughly estimated.
	This result can be interpreted as --- one  may choose the initial step-size based on the expected quality of the primal and dual warm-starts.  The default choice can be $ \rho_0 = 1 $, implying no  a priori knowledge at all. 
\end{remark}
\begin{remark}[warm start  specification]
	Above,  we did  not specify how to achieve a warm start, due  to it being  a separate issue. As aforementioned, this issue has been  intensively studied  in the  literature, see e.g. \cite{sambharya2023end,sambharya2023learning}.
\end{remark}

\subsection{Practical use}\label{sec_prac}
In the  previous sections, we have established the  theoretical optimal  step-size selection principle. However,  it is built on the optimal solutions, which are not known   a priori.
For practical use,  we may replace the optimal solutions with the  $ (k+1)$-th  iterate  information and update the step-size adaptively,  i.e., solve the  following polynomial instead of \eqref{quartic}:
\begin{equation} 
\rho^4\Vert \mathcal{A}\bm{x}^{k+1}\Vert^2  - \rho^3\langle \mathcal{A}\bm{x}^{k+1}, \bm{\zeta}^0  \rangle + \rho\langle \bm{\lambda}^{k+1}, \bm{\zeta}^0  \rangle -  \Vert\bm{\lambda}^{k+1}\Vert^2 = 0, \quad \rho \neq  0. 
\end{equation}

Let $ \bm{\zeta}^0 = \bm{0}$, we obtain a  simple version:
\begin{algorithm}[H]
	\caption{ ADMM  with  adaptive step-sizes (zero  initialization)}
	\label{Algo_zero}
	\begin{algorithmic}[1]
		\REQUIRE  Set $\,\bm{z}^0= \bm{0}, \,\bm{\lambda}^0 = \bm{0},\,   \rho_{0}  = 1 $.
		\WHILE{iteration $ k  =  0, 1, 2, \dots   $} \STATE
		\begin{align}
		\bm{x}^{k+1} =\,\,&  \, \underset{\bm{x}}{\text{argmin}} \,\, f(\bm{x})+ \frac{1}{2}\Vert  \rho_k (\mathcal{A}\bm{x} -   \mathcal{B}\bm{z}^k  - \bm{c})  + \bm{\lambda}^k/\rho_k  \Vert^2, \\				
		\bm{z}^{k+1} =\,\, &\, \underset{\bm{z}}{\text{argmin}} \,\, g(\bm{z})+ \frac{1}{2}\Vert  \rho_k(\mathcal{A}\bm{x}^{k+1} -   \mathcal{B}\bm{z}  - \bm{c})  + \bm{\lambda}^k/\rho_k  \Vert^2, \\
		\bm{\lambda}^{k+1} =\,\,  &\,\bm{\lambda}^{k} + (\rho_k)^2 \big( \mathcal{A}\bm{x}^{k+1} - \mathcal{B}\bm{z}^{k+1} - \bm{c} \big),\\
		\rho_{k+1}  	=\,\,  &\,    \pm \sqrt{\Vert\bm{\lambda}^{k+1}\Vert/\Vert\bm{x}^{k+1}\Vert} \quad  //\, \text{two choices are equivalent}.
		\end{align}	
		\ENDWHILE
		\ENSURE primal  and dual solutions $ \bm x^\star $  and  $ \bm \lambda^\star $, respectively.
	\end{algorithmic}
\end{algorithm} 
If one can  access  a good $ \bm{\zeta}^0 $, then the following version may be employed:
\begin{algorithm}[H]
	\caption{ ADMM  with adaptive step-sizes (warm start)}
	\label{Algo_1}
	\begin{algorithmic}[1]
		\REQUIRE Iterates warm-start $\bm{\lambda}^0\in  \mathbb K, \,\bm{z}^0 \in  \mathbb P$. Initial step-size $ \rho_{0}  $  by  estimating \eqref{re_ini}  (if not feasible,  set $ \rho_{0}  =1 $).
		Let $ \mathcal{A}\bm{x}^0 = \mathcal{B}\bm{z}^0 + \bm{c} $ and let
		$ \bm{\zeta}^0  = \rho_{0}\mathcal{A}\bm{x}^0 +   \bm{\lambda}^0/\rho_{0} $.
		\WHILE{iteration $ k  =  0, 1, 2, \dots   $} \STATE 
		\begin{align}
		\bm{x}^{k+1} =\,\,&  \, \underset{\bm{x}}{\text{argmin}} \,\, f(\bm{x})+ \frac{1}{2}\Vert  \rho_k (\mathcal{A}\bm{x} -   \mathcal{B}\bm{z}^k  - \bm{c})  + \bm{\lambda}^k/\rho_k  \Vert^2, \nonumber\\				
		\bm{z}^{k+1} =\,\, &\, \underset{\bm{z}}{\text{argmin}} \,\, g(\bm{z})+ \frac{1}{2}\Vert  \rho_k(\mathcal{A}\bm{x}^{k+1} -   \mathcal{B}\bm{z}  - \bm{c})  + \bm{\lambda}^k/\rho_k  \Vert^2, \nonumber\\
		\bm{\lambda}^{k+1} =\,\,  &\,\bm{\lambda}^{k} + (\rho_k)^2 \big( \mathcal{A}\bm{x}^{k+1} - \mathcal{B}\bm{z}^{k+1} - \bm{c} \big). 	
		\end{align}		
		\STATE Compute the $ (k+1)$-th  step-size choice $  \rho_{k+1} \neq  0 $ as a root of
		\begin{equation}\label{root0}
		\rho^4\Vert \mathcal{A}\bm{x}^{k+1}\Vert^2  - \rho^3\langle \mathcal{A}\bm{x}^{k+1}, \bm{\zeta}^0  \rangle + \rho\langle \bm{\lambda}^{k+1}, \bm{\zeta}^0  \rangle -  \Vert\bm{\lambda}^{k+1}\Vert^2 = 0,  
		\end{equation}
		with  a closed-from expression specified in  Sec. \ref{closed_forms}. 
		\ENDWHILE
		\ENSURE primal  and dual solutions $ \bm x^\star $  and  $ \bm \lambda^\star $, respectively.
	\end{algorithmic}
\end{algorithm}

\begin{remark}[adaptive convergence]
	Numerically, we find  that  the above two algorithms always converge in all our experiments. 	By a closer look, we find that the  convergence  owes to that the adaptive step-size will soon converge  to the  fixed optimal step-size.
	
	Theoretically, we may  establish  convergence by adding a  criterion to force the  algorithm stop  updating the step-size after some iterations. Then,  the ADMM convergence result for a fixed step-size can apply, see Proposition \ref{fpi}.   
	(Meanwhile, such an early stop can save some  runtime.)
\end{remark}

\begin{remark}[literature connection]
	The  above adaptive estimation of the optimal point  information has been successfully applied in the importance sampling from the machine learning field,  see  \cite[Sec. IV. C]{9904868}, \cite{7738878}.
\end{remark}

\begin{remark}[non-trivial only view]
	The above is a balanced/fixed-point view. Similar adaptive strategy also applies to the partly  trivial  case in Proposition \ref{pro_part}. For example, consider its non-trivial primal case (i).
	Let 
	$ 	\omega_1^k  =  \arccos  \, {\langle \mathcal{A}\bm{x}^k, \bm{\zeta}^0  \rangle }/({\Vert\mathcal{A}\bm{x}^k\Vert\Vert \bm{\zeta}^0\Vert}) $. Then,
	adaptively set $ \rho_k =  {\Vert\bm{\zeta}^0 \Vert}/  ({\Vert\mathcal{A}\bm{x}^k\Vert}{\cos \omega_1^k}) $.
	Same arguments apply  to dual case  (ii). See numerical details  in Sec. \ref{sec_5_2_2}.
\end{remark}

\section{Numerical examples}\label{sec_5}
Below,  we test our results via  two   numerical  examples. For the general non-trivial case,  we  consider the popular  Lasso application; for the  partly trivial  case, we study the feasibility problem.
\subsection{Lasso}\label{num_lasso}
Consider the following problem:
\begin{equation}
\text{minimize}    \,\, 		\frac{1}{2}\Vert\bm{A}\bm{x} - \bm{b}\Vert^2 +  \alpha\Vert\bm{z}\Vert_1,  \quad
\text{subject to}  \,\, \,		\bm{x} = \bm{z}.
\end{equation}
with variable $ \bm{x}\in \mathbb{R}^{1000}, \bm{b} \in \mathbb{R}^{300}, \bm{A} \in \mathbb{R}^{ 300 \times 1000} $.   	
We   generate elements of $  \bm{A}, \bm{x} $ from a normal distribution $ \mathcal{N}(0,1) $.  
We  set $ \bm{b} = \bm{A}\bm{x}_0 +  \epsilon_{\text{noise}} $, where  $ \bm{x}_0  $ is  a random sparse vector, with density factor $ 1/3 $,  and  where $ \epsilon_{\text{noise}}   $ is some  random  noise  from  $ \mathcal{N}(0,0.1) $. At last, we set  regularization parameter $ \alpha = 0.1  \Vert \bm{A}^T \bm b \Vert_{\infty}$.

Unless specified,  the figures   are  plotted  against  the  following normalized error measures:
\begin{align}
\frac{\Vert \bm{\zeta}^{k+1} - \bm{\zeta}^{k} \Vert^2   }{\Vert\bm{x}^\star \Vert\Vert\bm{\lambda}^\star \Vert },
\quad
\frac{\Vert\bm{x}^{k+2} - \bm{x}^{k+1} \Vert^2   }{\Vert\bm{x}^\star \Vert^2 },
\quad
\frac{\Vert \bm{\lambda}^{k+1} - \bm{\lambda}^{k} \Vert^2   }{\Vert\bm{\lambda}^\star \Vert^2 }.  \tag{normalized errors}
\end{align}

\begin{remark}[data conditioning]
	Above, we did not apply  column-normalization  for the sensing matrix $ \bm{A} $.
	We find  that  if  normalized, the optimal  choice $\rho^\star$  is roughly  at $ \pm 0.9 $,  i.e., close  to the unit  length, 
	implying the  data is originally well-conditioned, where the step-size selection  becomes less important.
\end{remark}

\subsubsection{zero initialization}
\begin{figure}[H]
	\centering
	\begin{subfigure}[t]{0.42\textwidth}
		\includegraphics[scale=0.42]{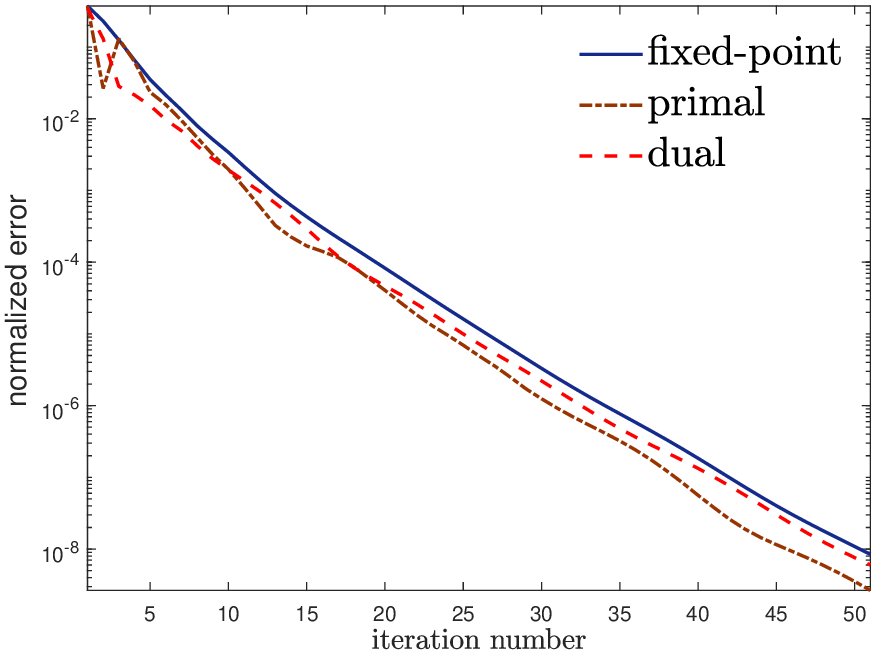}\caption{Step-size: worst-case optimal.}\label{figthe}
	\end{subfigure}
	\begin{subfigure}[t]{0.42\textwidth}
		\includegraphics[scale=0.42]{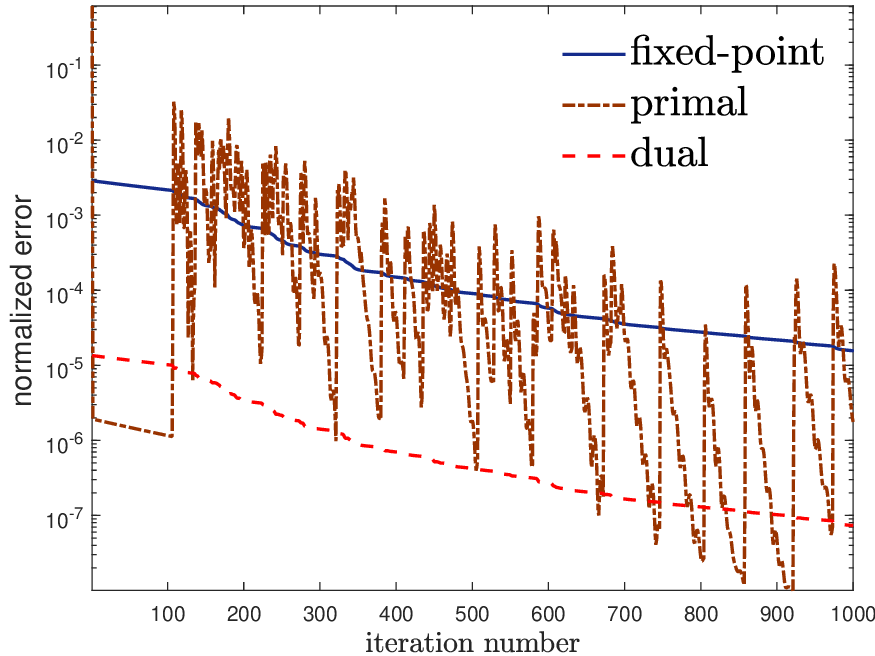}\caption{Step-size: unit length.}\label{fi01}
	\end{subfigure}
	\caption{(Fixed step-size) Theoretical results, see Proposition  \ref{pro_non}.}
	\begin{subfigure}[t]{0.42\textwidth}
		\includegraphics[scale=0.42]{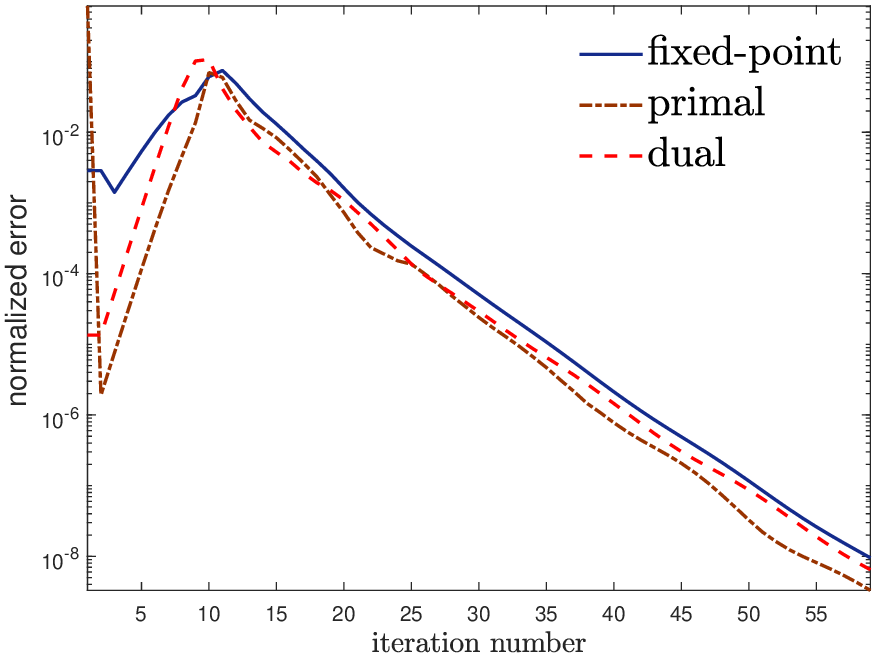}\caption{Convergence: 3 sequences.}\label{est_a}
	\end{subfigure}
	\begin{subfigure}[t]{0.42\textwidth}
		\includegraphics[scale=0.42]{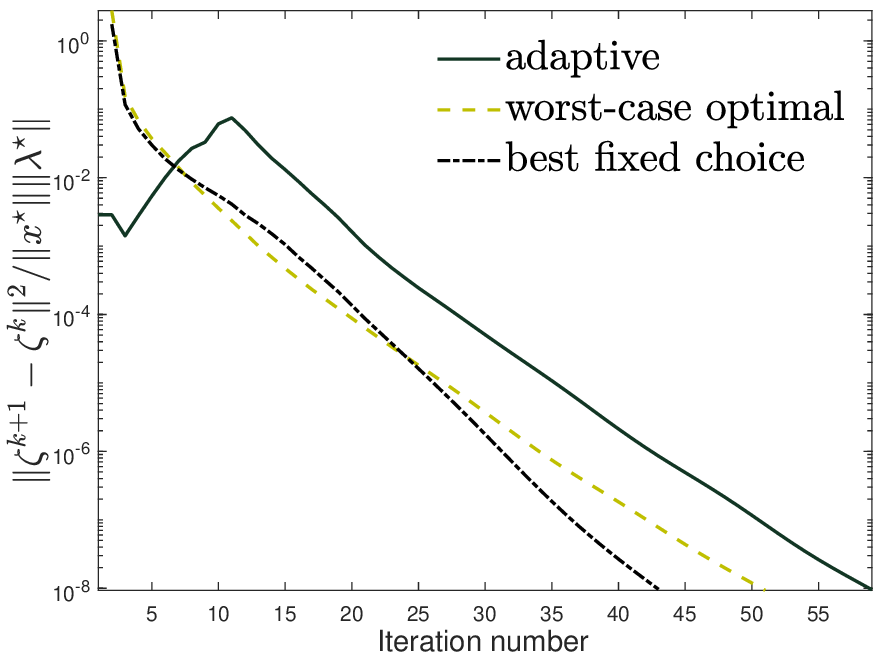}\caption{Comparison: fixed-point sequence.}\label{est_b}
	\end{subfigure}
	\caption{(Adaptive step-size)  Practical  use via Algorithm \ref{Algo_zero}.}
\end{figure}

\begin{remark}[balanced rates]
	The above numerical results  do align with our theoretical conclusions  --- the  optimal  step-size would set all sequences (fixed-point, primal and dual) to the  same speed.  
\end{remark}
\begin{remark}[adaptive]
	Compare Fig. \ref{est_a} and Fig. \ref{figthe}, we see that the adaptive strategy  performs quite differently at the beginning stage, but soon (roughly 13 iterations) becomes stable  and similar to the theoretical one.
\end{remark}
\begin{remark}[near-limit performance]
	In  Fig.  \ref{est_b}, we investigate how much potential is left for further improvement, since ours is only optimal in the  worst-case sense  (no  tailored structure exploited).
	To this end, we fixed the random number generator with seed $ 0 $, and 
	found the underlying best fixed step-size choice  via an exhaustive grid search.	
	Quite remarkably, we see  that our  choice is  similar to  the underlying best fixed one.  
\end{remark}

\subsubsection{warm start}\label{sec_warm}
Here, we test the warm start  case corresponding  to Sec. \ref{sec_opt_ini} with     settings:  
(i)   Full warm-start: $ \bm{x}^0 = \bm{x}^\star + \epsilon_{\text{err}1},  \,\, \bm{\lambda}^0  =  \bm{\lambda}^\star + \epsilon_{\text{err}2}  $;
(ii) Partial warm-start:  $ \bm{x}^0 =  \bm{x}^\star + \epsilon_{\text{err}1},  \,\,\bm{\lambda}^0  = \bm{0} $;
(iii) Random initialization: $ \bm{x}^0  \sim  \mathcal{N}(0, 1) ,  \,\, \bm{\lambda}^0 \sim  \mathcal{N}(0, 1) $.

\begin{figure}[H]
	\centering
	\begin{subfigure}[t]{0.42\textwidth}
		\includegraphics[scale=0.42]{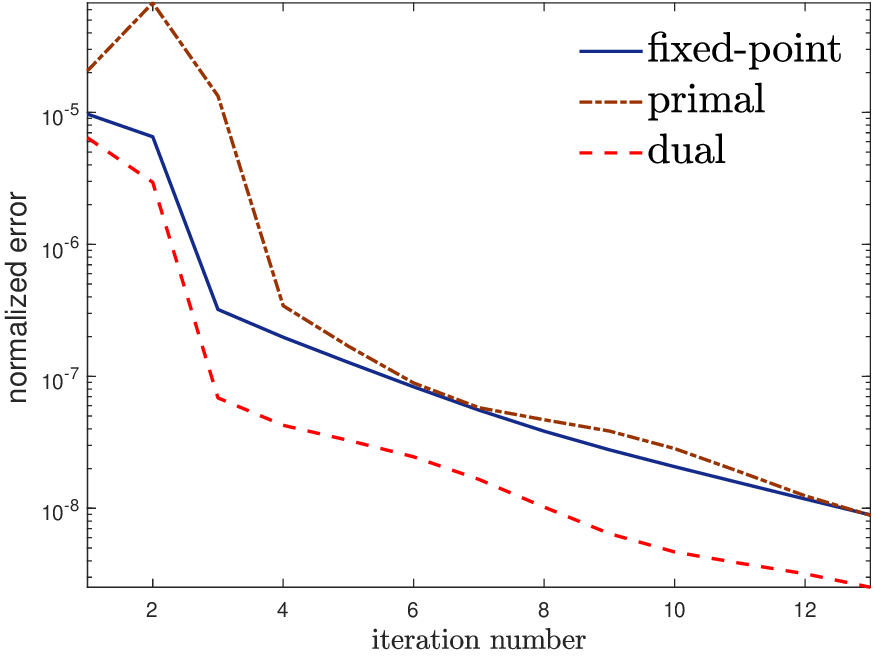}\caption{Algorithm \ref{Algo_1} with $\rho_{0} = 10$ (estimated).}
	\end{subfigure}
	\begin{subfigure}[t]{0.42\textwidth}
		\includegraphics[scale=0.42]{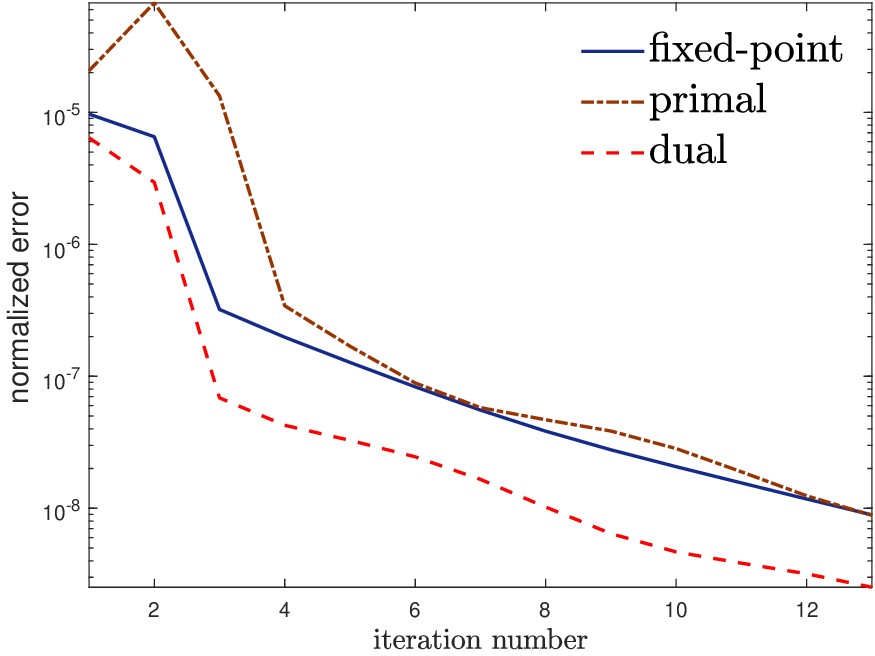}\caption{Fixed $\rho_k= \rho_0  = 10, \forall k$ (estimated).}
	\end{subfigure}
	\begin{subfigure}[t]{0.42\textwidth}
		\includegraphics[scale=0.42]{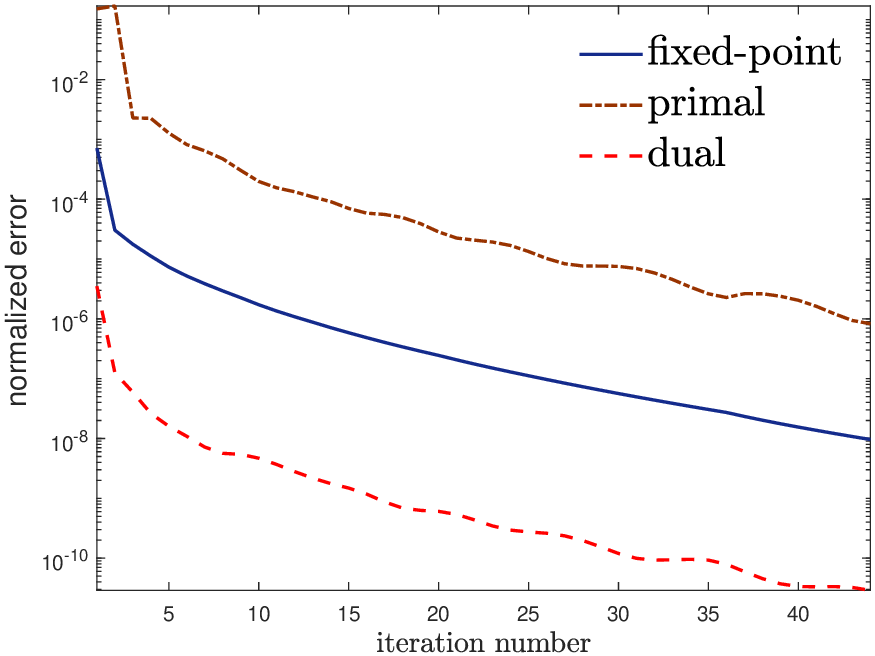}\caption{Algorithm \ref{Algo_1} with $\rho_{0} = 1$.\label{fig_3.1}}
	\end{subfigure}
	\begin{subfigure}[t]{0.42\textwidth}
		\includegraphics[scale=0.42]{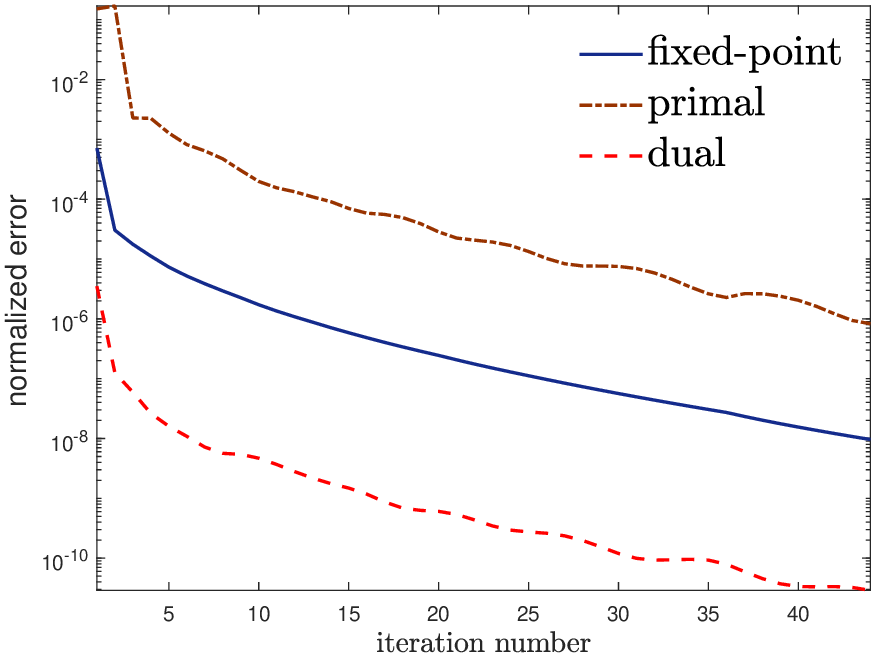}\caption{Fixed $\rho_k= \rho_0  = 1, \forall k$.\label{fig_3.2}}
	\end{subfigure}
	\caption{Full warm-start with $ \epsilon_{\text{err}1}  \sim \mathcal{N}(0, 10^{-3})  $, $ \epsilon_{\text{err}2} \sim  \mathcal{N}(0, 10^{-1})  $.\label{fig_c1}}
\end{figure}

\begin{figure}[H]
	\centering
	\begin{subfigure}[t]{0.42\textwidth}
	\includegraphics[scale=0.42]{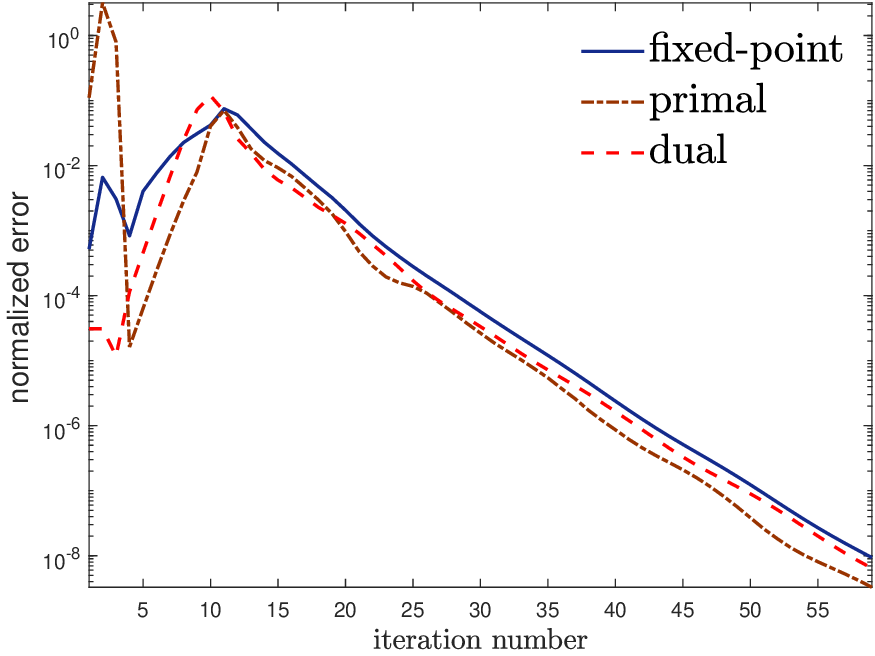}\caption{Algorithm \ref{Algo_1} with $ \rho_0 = 1 $.}
\end{subfigure}
\begin{subfigure}[t]{0.42\textwidth}
	\includegraphics[scale=0.42]{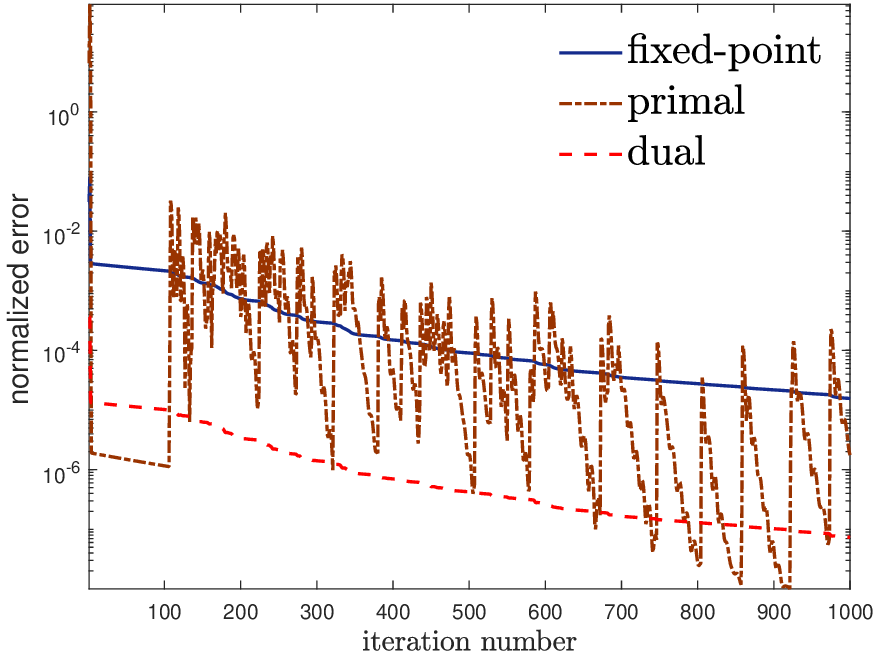}\caption{Fixed $\rho_k= \rho_0  = 1, \forall k$.}
\end{subfigure}
\caption{Partial warm-start with $ \epsilon_{\text{err}1}  \sim \mathcal{N}(0, 1)  $. \label{fig_partial}}
\begin{subfigure}[t]{0.42\textwidth}
	\includegraphics[scale=0.42]{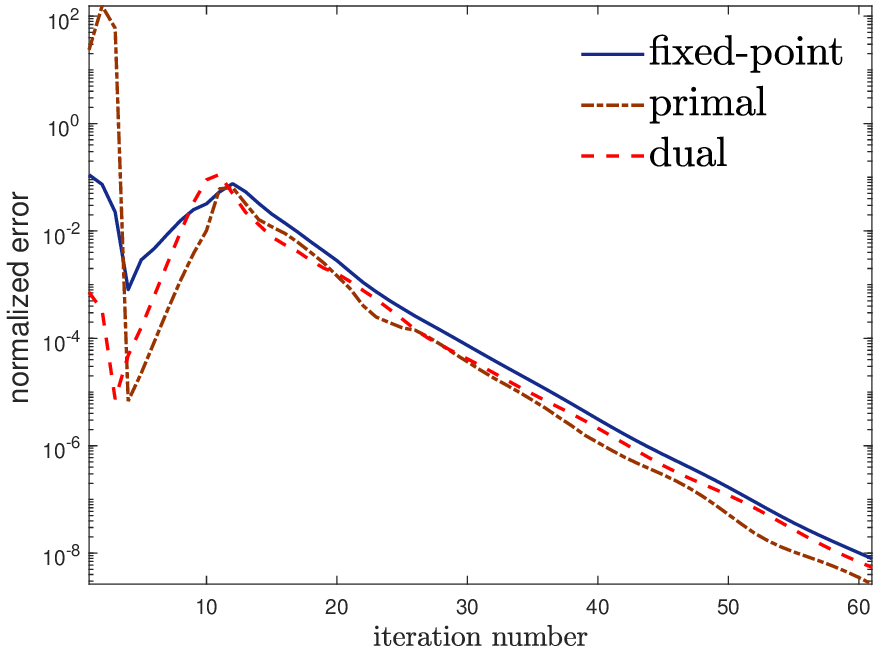}\caption{Algorithm \ref{Algo_1} with $ \rho_0 = 1 $. }
\end{subfigure}
\begin{subfigure}[t]{0.42\textwidth}
	\includegraphics[scale=0.42]{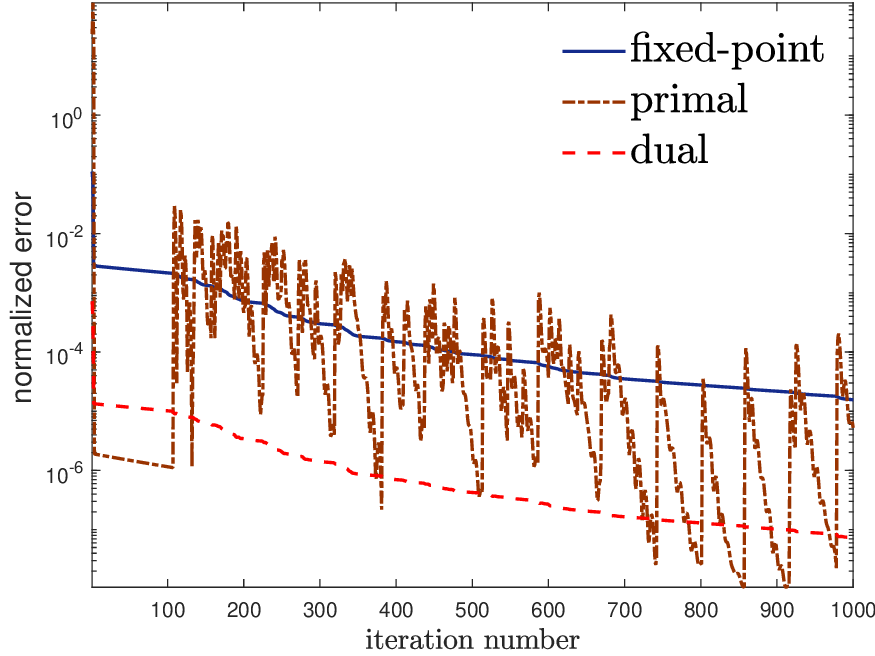}\caption{Fixed $\rho_k= \rho_0  = 1, \forall k$.}
\end{subfigure}
\caption{Random initialization: $ \bm{x}^0  \sim  \mathcal{N}(0, 1) ,  \quad \bm{\lambda}^0 \sim  \mathcal{N}(0, 1) $. \label{fig_rand}}
\end{figure}

\begin{remark}[full warm start]
	In  view of Fig. \ref{fig_c1}, our Algorithm \ref{Algo_1} is roughly identical to the fixed  step-size strategy. 
	By a closer look,  we find that  our adaptive step-size update  in \eqref{root0} outputs  almost-identical value as $\rho_0$.
	In fact,  this  is not too surprising by viewing \eqref{fix_01}. 
	That is, given a full warm-start setting, the best fixed-point is likely  the one corresponds to  $\rho_0$. 
	
	Additionally, due to the warm-start is not exactly optimal. We can change $\rho_0$ to accelerate certain  sequences. Here,  by  $\rho_0  =  10$, we are promoting the fixed-point sequence rate, which  is roughly $ 4 \times $  faster than setting  $\rho_0  =  1$.
\end{remark}

\begin{remark}[other initializations]
	We see that our  algorithm exhibits great advantages in  Fig \ref{fig_partial} and Fig \ref{fig_rand}, regarding  partial warm start and random initialization, respectively.
	
	Additionally, we omit the  case $  \bm{x}^0  = \bm{0},  \bm{\lambda}^0  = \bm{\lambda}^\star + \epsilon_{\text{err}2}  $, since  for Lasso,  the optimal  primal $ \bm{x}^\star $ is highly  sparse,   similar  to a zero  vector. 
	We observe that such a setting is similar to  the full warm start case.  
\end{remark}

\subsection{Partly  trivial: feasibility problem}\label{feasi}
Here, we consider a special application where the dual solution  is trivial,  i.e.,  a  zero vector $ \bm{\lambda}^\star = \bm{0} $.
Consider the following  feasibility problem:
\begin{align}
\text{Find}    \quad			&\bm{x} \in \mathbb{H} ,\nonumber \\
\text{subject to}  \quad 		&\bm{x} \in C,   \nonumber\\
&\bm{x} \in D,
\end{align}
where subsets $ C , D   $ are  specified as 
\begin{align}\label{sets}
C =  \{  \bm{y}\in  \mathbb{R}^{200} \, |\, \bm{A}_1 \bm{y} =  \bm{b}_1 \} , \quad 
D =  \{  \bm{z}\in  \mathbb{R}^{200}\, |\, \bm{A}_2 \bm{z} =  \bm{b}_2 \} ,
\end{align}
with $ \bm{A}_1\in \mathbb{R}^{80\times200} $, $ \bm{b}_1 \in  \mathbb{R}^{80}$, and $ \bm{A}_2\in \mathbb{R}^{100\times200} $, $ \bm{b}_2 \in  \mathbb{R}^{100}$.

It can be rewritten into an ADMM form:
\begin{equation}
\text{minimize}    \,\, 		\bm{1}_C(\bm{x}) + \bm{1}_D(\bm{z}) \quad 
\text{subject to}  \,\, \,		\bm{x} = \bm{z},
\end{equation}
with $ \bm{1}_C $ and $ \bm{1}_D $ the  indicator functions corresponding to \eqref{sets}.

The figures  in this section  are  plotted  against  the  following normalized error measures, corresponding to   Proposition  \ref{pro_non}:
\begin{align}
{\Vert \bm{\zeta}^{k+1} - \bm{\zeta}^{k} \Vert^2   },
\quad
{\Vert\bm{x}^{k+2} - \bm{x}^{k+1} \Vert^2   },
\quad
{\Vert \bm{\lambda}^{k+1} - \bm{\lambda}^{k} \Vert^2   }.  \tag{error measures}
\end{align}

\subsubsection{zero initialization: fixed non-trivial  rate}
For the above problem,   set $ \bm{\zeta}^0  =  \bm{0} $. Then,  Proposition \ref{prop_zero}  reduces to: 
\begin{align}
\qquad	{\Vert\bm{x}^{k+2} - \bm{x}^{k+1} \Vert^2   }
\,\,\leq\,\,\,\,    & \frac{1}{k+1} {  \Vert\bm{x}^\star\Vert^2 }    , \tag{primal}\label{zer1} \\
\qquad{\Vert \bm{\lambda}^{k+1} - \bm{\lambda}^{k} \Vert^2   }
\,\, \leq  \,\,\,\,  &  \frac{1}{k+1} {  \Vert\bm{x}^\star\Vert^2 }\cdot  \rho^4      ,  \tag{trivial dual}   \\
\qquad	{\Vert \bm{\zeta}^{k+1} - \bm{\zeta}^{k} \Vert^2   }
\,\, \leq  \,\,\,\,   &  \frac{1}{k+1} {  \Vert\bm{x}^\star\Vert^2 }\cdot  \rho^2    .  \tag{fixed-point} \label{zer2}
\end{align}
Here, we numerically verify this result. For clarity, below we only plot the primal and fixed-point sequences.
\begin{figure}[H]
	\centering
	\includegraphics[scale=0.5]{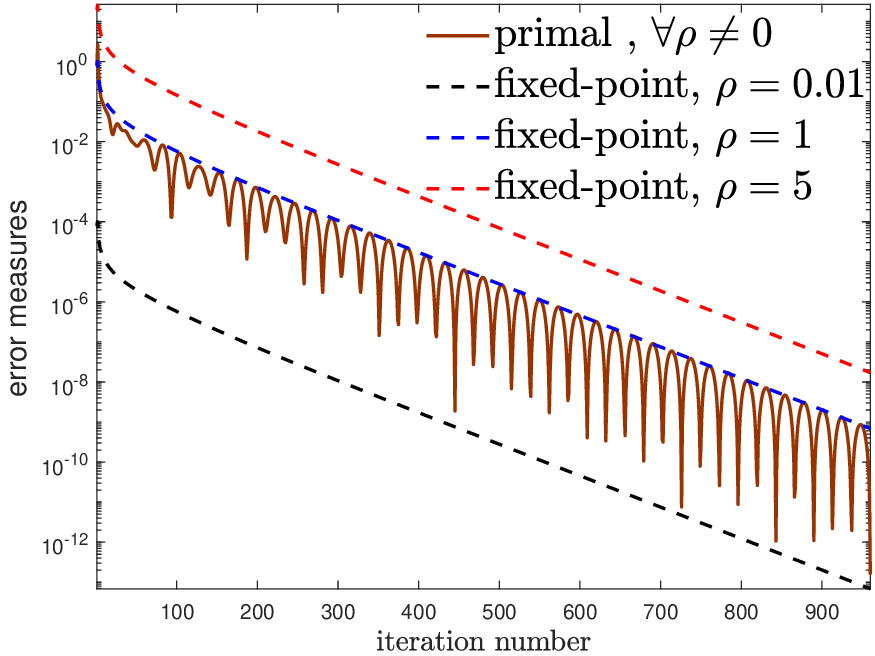}\caption{Fixed rate for primal iterate, independent of step-size  choices.}
\end{figure}

\begin{remark}
	From above,  we see that primal  sequence (non-trivial) has exactly the  same practical convergence rate,  regardless of the step-size choices. Meanwhile, the fixed-point rate does  improve as the step-size $ \rho $  approaching $ 0 $, which  is aligned  with the above theory \eqref{zer2}.
\end{remark}

\subsubsection{non-zero initialization}\label{sec_5_2_2}
For the above problem,   set   any $ \bm{\zeta}^0  \neq  \bm{0} $.
Let $ \omega  =  \arccos  \, {\langle \bm{x}^\star, \bm{\zeta}^0  \rangle }/{\Vert\bm{x}^\star\Vert\Vert \bm{\zeta}^0\Vert}  $. 

$\bullet$ (i) (balanced view) 
Set $ \rho^\star=  \frac{\Vert\bm{\zeta}^0 \Vert}{\Vert\bm{x}^\star\Vert} \cos \omega $.
Then,    Proposition \ref{pro0} reduces to
\begin{align}
{\Vert \bm{x}^{k+2} - \bm{x}^{k+1} \Vert^2  }
\,\, \leq  \,\,      \frac{1}{k+1} \Vert  \bm{x}^\star \Vert^2 \cdot \tan^2 \omega , \label{r_comp1}  
\end{align}

$\bullet$  (ii) (non-trivial only)
Set $ \rho^\star_\text{pri} =  \frac{\Vert\bm{\zeta}^0 \Vert}{\Vert\bm{x}^\star\Vert} \frac{1}{\cos \omega} $. Then,   Proposition \ref{pro_part} gives
\begin{align}
{\Vert \bm{x}^{k+2} - \bm{x}^{k+1} \Vert^2  }
\,\, \leq  \,\,      \frac{1}{k+1} \Vert  \bm{x}^\star \Vert^2 \cdot \sin^2 \omega . \label{r_comp2} 
\end{align}

For practical use,  let
\begin{equation}
\omega_{k+1}  =  \arccos  \, \frac{\langle \bm{x}^{k+1}, \bm{\zeta}^0  \rangle }{\Vert\bm{x}^{k+1}\Vert\Vert \bm{\zeta}^0\Vert}, \qquad k = 0,1,\dots,
\end{equation}
and adaptively update the step-size via
\begin{align}
\qquad	\rho_{k+1}  	=\,\,  &\,    \frac{\Vert\bm{\zeta}^0 \Vert}{\Vert\bm{x}^{k+1}\Vert} \cdot   \cos \omega_{k+1}   ,  \tag{balanced view} \label{adapt_step1}\\
\qquad	\rho_{k+1}  	=\,\,  &\,   \frac{\Vert\bm{\zeta}^0 \Vert}{\Vert\bm{x}^{k+1}\Vert} \cdot   \frac{1}{\cos \omega_{k+1} } \tag{primal iterate view} \label{adapt_step2}.
\end{align}

\begin{figure}[H]
	\centering
	\begin{subfigure}[t]{0.45\textwidth}
		\includegraphics[scale=0.42]{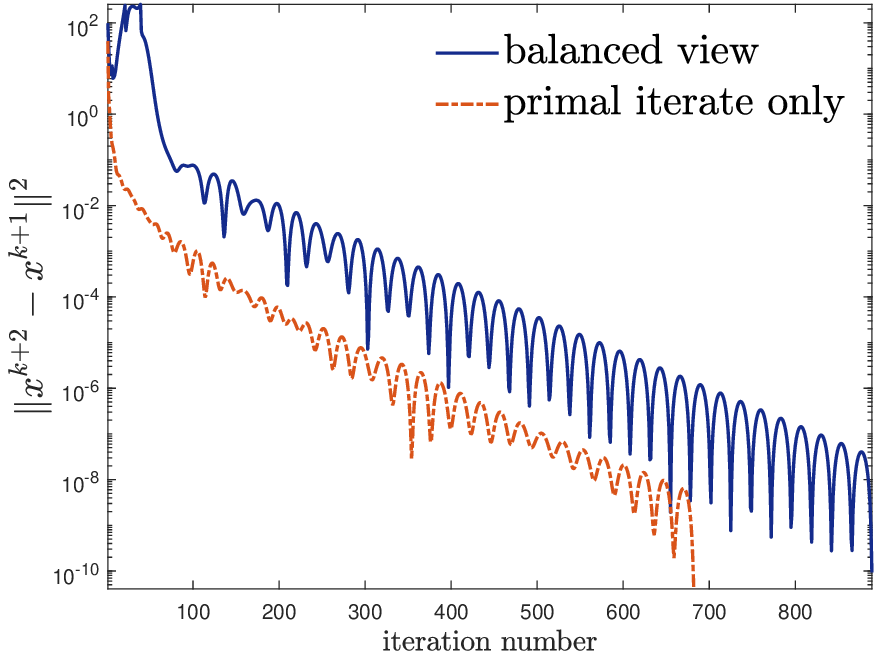}\caption{ $ \epsilon_{\text{err}1}  \sim \mathcal{N}(0, 1) $.}\label{ffig1}
	\end{subfigure}
	\begin{subfigure}[t]{0.45\textwidth}
		\includegraphics[scale=0.42]{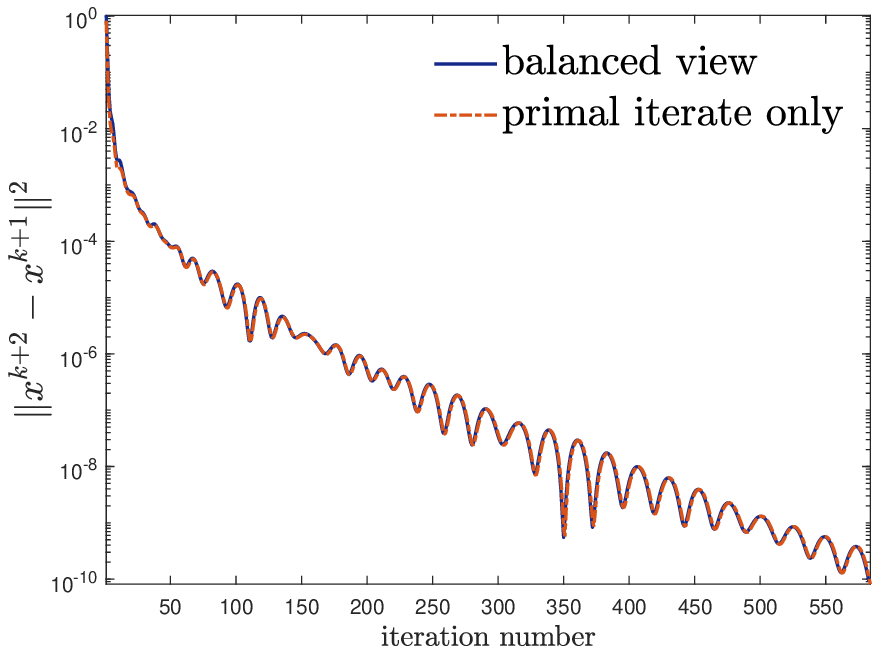}\caption{$ \epsilon_{\text{err}1}  \sim \mathcal{N}(0, 0.1) $.}\label{ffig2}
	\end{subfigure}
	\caption{(Adaptive) Primal warm-start: $ \bm{x}^0 =  \bm{x}^\star + \epsilon_{\text{err}1},  \, \bm{\lambda}^0  = \bm{0} $.}
\end{figure}

\begin{remark}
	In  view of \eqref{r_comp1} and \eqref{r_comp2}, we see  that the two worst-case rates admit a factor difference of $ \cos^2 \omega $.
	(i) For  Fig. \ref{ffig1},  we find   $ \cos^2  \omega \approx 0.023  $.  
	Theoretically,  we expect different  convergence  rates, which is indeed the case numerically.
	(ii)  For  Fig. \ref{ffig2},  we find  $ \cos^2  \omega \approx 0.78 $, which is close to $ 1 $, we expect them to be  similar,  and numerically they admit the same  performance.
\end{remark}

\subsubsection{Comparison}
Here, we make   comparison  between  zero and non-zero initialization. 
We are   motivated by the following theoretical results:
\begin{align}
\qquad	{\Vert\bm{x}^{k+2} - \bm{x}^{k+1} \Vert^2   }                                                                   
\,\,\leq\,\,\,\,    & \frac{1}{k+1} {  \Vert\bm{x}^\star\Vert^2 }   , \tag{zero ini., see \eqref{zer1}}\\
\qquad	{\Vert\bm{x}^{k+2} - \bm{x}^{k+1} \Vert^2   }                                                                   
\,\,\leq\,\,\,\,    & \frac{1}{k+1} {  \Vert\bm{x}^\star\Vert^2 }   \cdot \sin^2 \omega.\tag{non-zero ini., see  \eqref{r_comp2}}
\end{align}
\begin{figure}[H]
	\centering
	\includegraphics[scale=0.5]{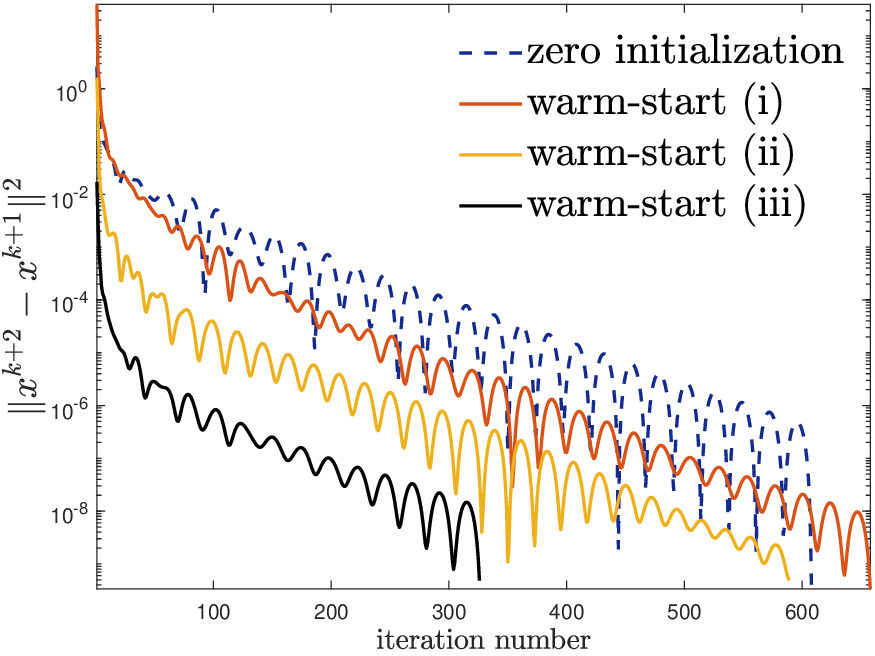}\caption{(Adaptive) Comparison: warm-start  $ \bm{x}^0 =  \bm{x}^\star + \epsilon_{\text{err}1},  \, \bm{\lambda}^0  = \bm{0} $, with
		(i) $ \epsilon_{\text{err}1} \sim \mathcal{N}(0,1) $; 
		(ii) $ \epsilon_{\text{err}1} \sim \mathcal{N}(0,0.1) $; 
		(iii) $ \epsilon_{\text{err}1} \sim \mathcal{N}(0,0.01) $.}
\end{figure}
\begin{remark}
	Numerically, the  3 warm-start settings  above admit angle  values: (i) $ \sin^2 \omega \approx 0.98 $; (ii) $ \sin^2 \omega \approx 0.22 $; (iii) $ \sin^2 \omega \approx 0. 003 $.
	For a smaller $ \sin^2 \omega $, we  do observe an improved convergence rate, which aligns with the theoretical result.
\end{remark}

\section{Conclusion}
In this work, we present a general, simple, worst-case optimal step-size selection principle for ADMM-type algorithms, which solves a long-standing open issue. There are two main types of choices: (i) step-size that accelerates the fixed-point sequence, which yields some balanced rates; (ii) accelerating only the primal or dual iterate, at the cost of slowing down all other sequences.

The key to our success owes to a novel domain-type parametrization, symmetric to the range-type in  the current literature, which avoids an implicit step-size related scaling/shifting. 
Interestingly, we reveal that  the algorithm's  efficiency (optimal convergence rate) is intrinsically related to the primal-dual solution and initialization  angles.  This discovery appears useful and  may  be worth further investigation in the future.

\bibliographystyle{unsrt}
\bibliography{Reference/ref1,Reference/Ref_S,Reference/ML_application,Reference/sr_applications}

\appendix
\section{Additional  experiments}
This material aims to  support the generality claim of our step-size selection scheme. 
All experiments  below use the standard data settings  from the review paper by  Stephen Boyd et al.  \cite{boyd12}, with open-sourced code from their website (relatively  well-conditioned data).

All figures  in this section  are  plotted  against  the  following measures (unless specified),  which  we referred to as the normalized errors:
\begin{align}
\frac{\Vert \bm{\zeta}^{k+1} - \bm{\zeta}^{k} \Vert^2   }{\Vert\mathcal{A}\bm{x}^\star \Vert\Vert\bm{\lambda}^\star \Vert },
\quad
\frac{\Vert \mathcal{A}\bm{x}^{k+2} - \mathcal{A}\bm{x}^{k+1} \Vert^2   }{\Vert\mathcal{A}\bm{x}^\star \Vert^2 },
\quad
\frac{\Vert \bm{\lambda}^{k+1} - \bm{\lambda}^{k} \Vert^2   }{\Vert\bm{\lambda}^\star \Vert^2 }.  
\end{align}

For non-zero initialization, we test the following settings:\\
(i)  $ \,\, $ Full warm-start:  $ \bm{x}^0 = \bm{x}^\star + \epsilon_{\text{err}1},  \quad \bm{\lambda}^0  =  \bm{\lambda}^\star + \epsilon_{\text{err}2}  $;\\
(ii) $ \, $ Partial primal warm-start: $ \bm{x}^0 =  \bm{x}^\star + \epsilon_{\text{err}1},  \quad \bm{\lambda}^0  = \bm{0} $;\\
(iii) Partial dual warm-start: $ \bm{x}^0  = \bm{0},  \quad \bm{\lambda}^0  = \bm{\lambda}^\star + \epsilon_{\text{err}2} $.\\
(iv) Random initialization: $ \bm{x}^0  \sim  \mathcal{N}(0, 1) ,  \quad \bm{\lambda}^0 \sim  \mathcal{N}(0, 1) $.

\subsection{Linear  programming}
Consider the following linear  program:
\begin{equation}
\text{minimize}     \,\,  			\bm{c}^T\bm{x} \quad
\text{subject to}  \,\, 		\bm{Ax=b}, \,\,  \bm{x}\geq 0,
\end{equation}
with variable $ \bm{x},\bm{c} \in \mathbb{R}^{500},  \bm{b}\in \mathbb{R}^{400}, \bm{A} \in \mathbb{R}^{400\times 500} $. Vector $ \bm{c}  $ is generated from a uniform distribution on interval $ [0.5, 1.5] $. The elements of $ \bm{x}, \bm{A} $ are generated from a normal distribution $ \mathcal{N}(0,1) $ and then set to be positive by taking the absolute value. $ \bm{b} $ is set  equal to $ \bm{A}\bm{x} $.

\subsubsection{zero  initialization}
\begin{figure}[H]
	\centering
	\begin{subfigure}[t]{0.42\textwidth}
		\includegraphics[scale=0.42]{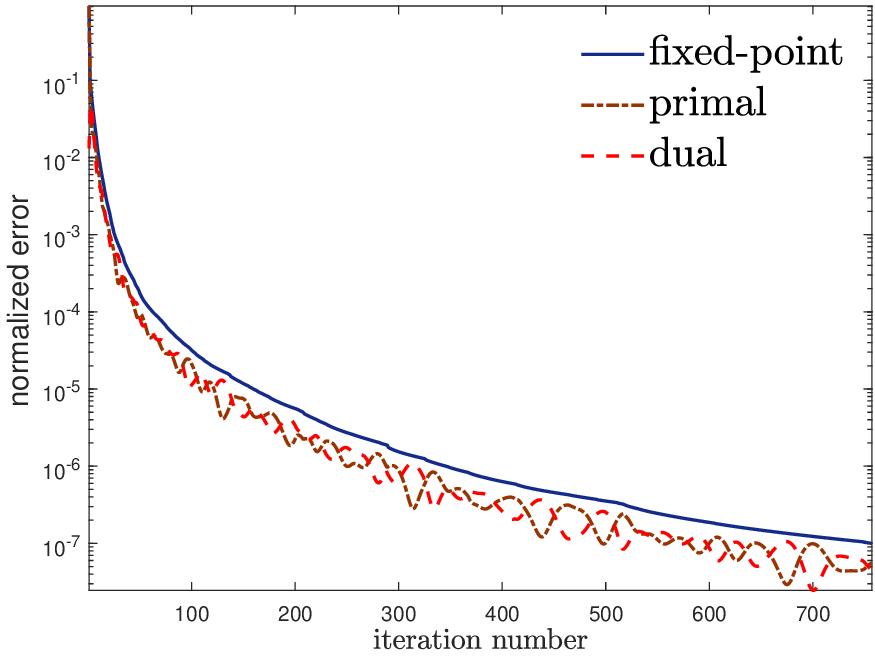}\caption{Step-size: optimal  choice $\rho^\star \approx \pm 0.6$.}\label{f01}
	\end{subfigure}
	\begin{subfigure}[t]{0.42\textwidth}
		\includegraphics[scale=0.42]{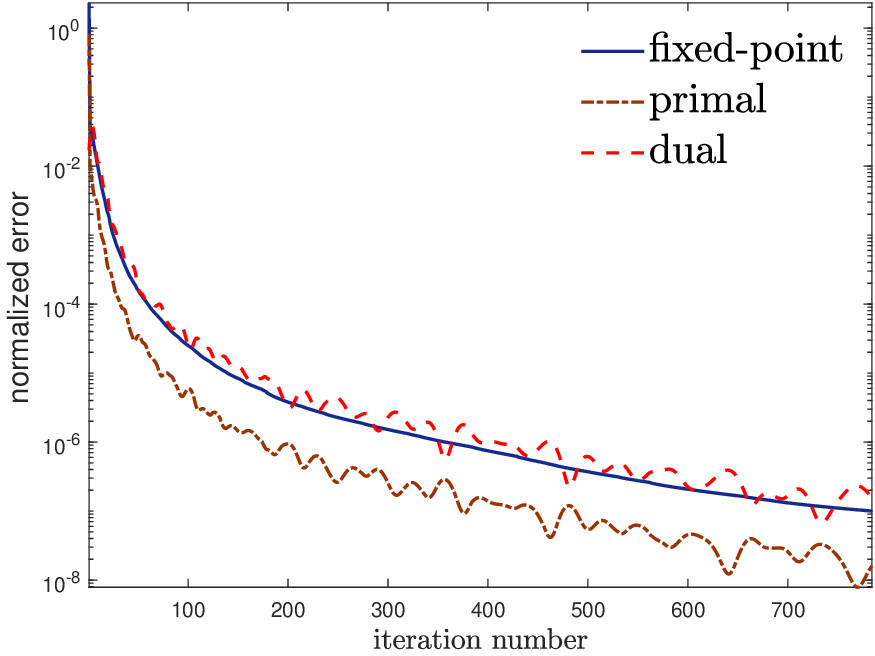}\caption{Step-size: unit length.}
	\end{subfigure}
	\caption{Theoretical results test (zero initialization).}
	\begin{subfigure}[t]{0.42\textwidth}
		\includegraphics[scale=0.42]{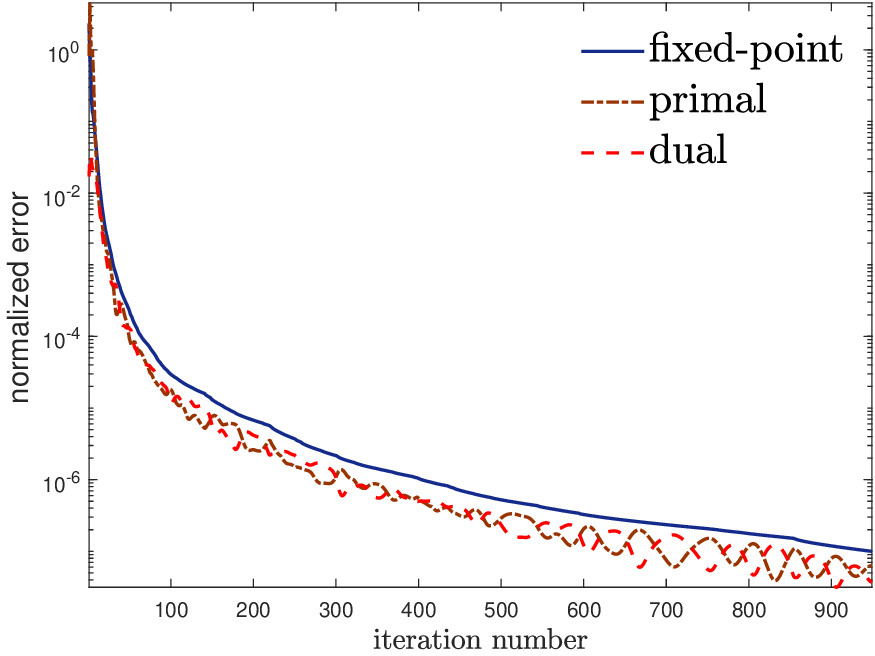}\caption{Estimated step-size choice (adaptive).}\label{fig_lp2}
	\end{subfigure}
	\begin{subfigure}[t]{0.42\textwidth}
		\includegraphics[scale=0.42]{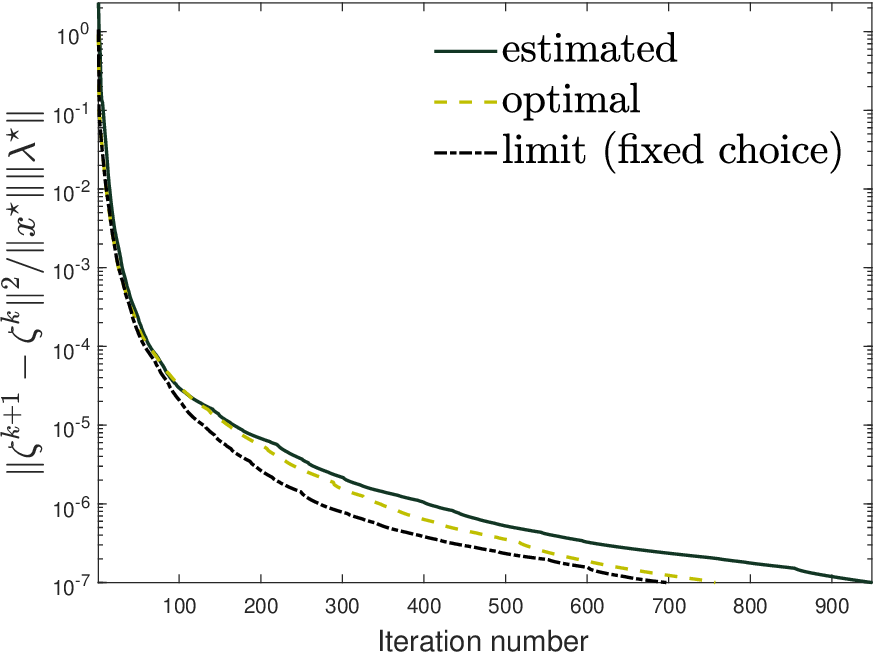}\caption{Fixed-point convergence comparison.}
	\end{subfigure}
	\caption{Practical  use  test (adaptive; zero initialization).}
\end{figure}

\begin{remark}
	The above  two step-size choices exhibit similar performances. With  a closer look, we find that the optimal choice in Fig. \ref{f01}  does yield a marginally faster fixed-point sequence convergence (requires 754 iterations, compared  to 778). Meanwhile, its primal iterates convergence is slightly slower than the unit-length step-size.
\end{remark}

\begin{remark}
	We see that the adaptive step-size does well-approximate  the previous  theoretical  result, and  exhibits near-limit (for fixed step-size) performance.
\end{remark}

\subsubsection{non-zero  initialization}
Set non-zero $\bm{\zeta}^0 \neq  \bm{0}$ according to the 4 strategies at the beginning of the appendix.
\begin{figure}[H]
	\centering
	\begin{subfigure}[t]{0.42\textwidth}
		\includegraphics[scale=0.42]{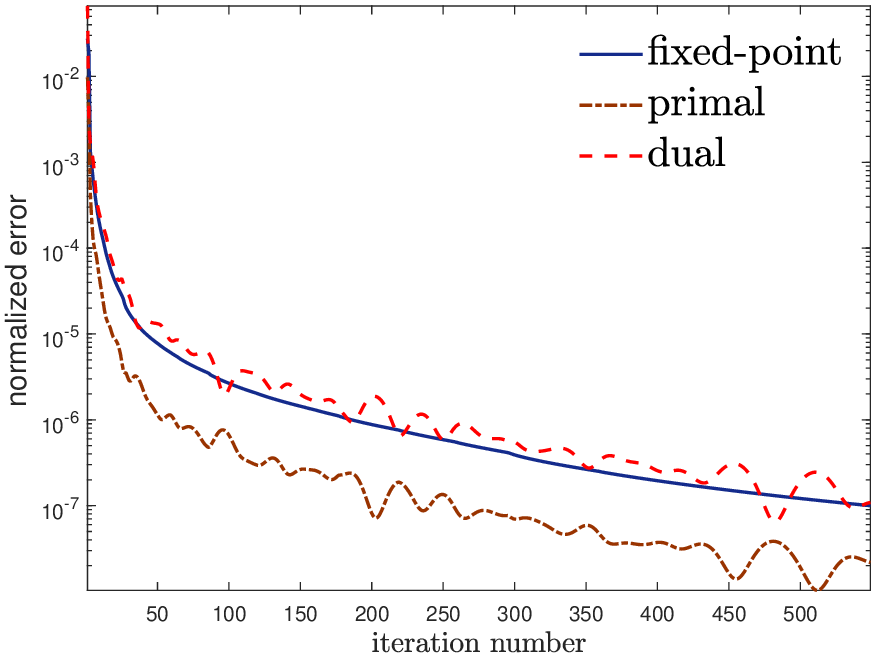}\caption{warm-start: $  \epsilon_{\text{err}1},\epsilon_{\text{err}2} \sim  \mathcal{N}(0, 0.1) $.}
	\end{subfigure}
	\begin{subfigure}[t]{0.42\textwidth}
		\includegraphics[scale=0.42]{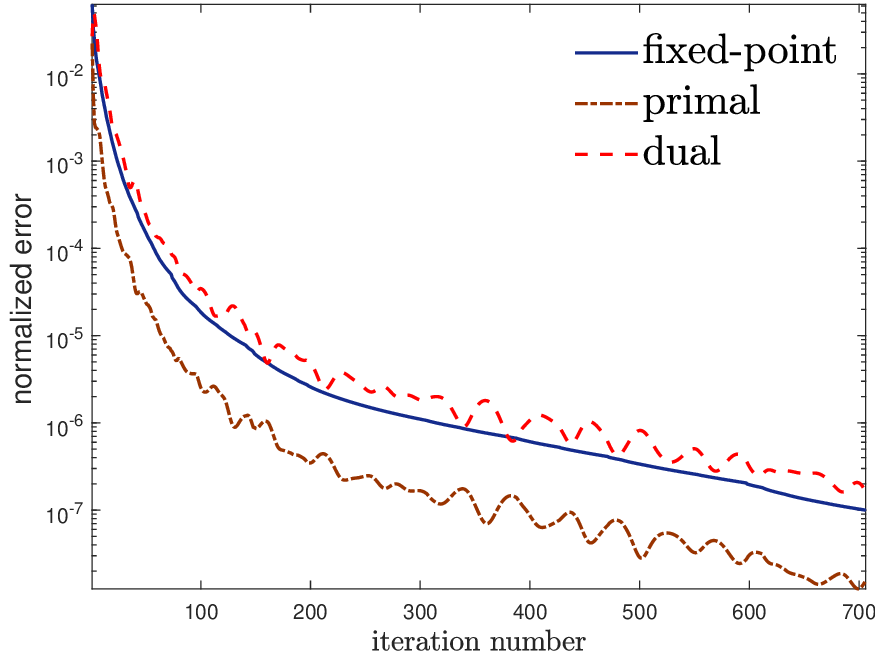}\caption{primal warm-start:  $  \epsilon_{\text{err}1} \sim  \mathcal{N}(0, 0.1) $.}
	\end{subfigure}	
	\begin{subfigure}[t]{0.42\textwidth}
		\includegraphics[scale=0.42]{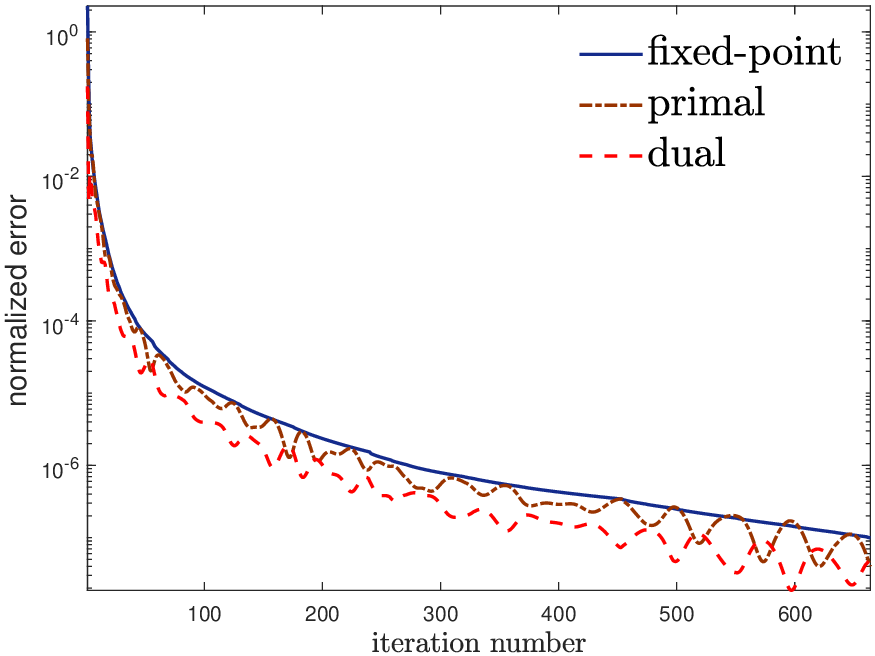}\caption{dual warm-start: $  \epsilon_{\text{err}2} \sim  \mathcal{N}(0, 0.1) $.}
	\end{subfigure}
	\begin{subfigure}[t]{0.42\textwidth}
		\includegraphics[scale=0.42]{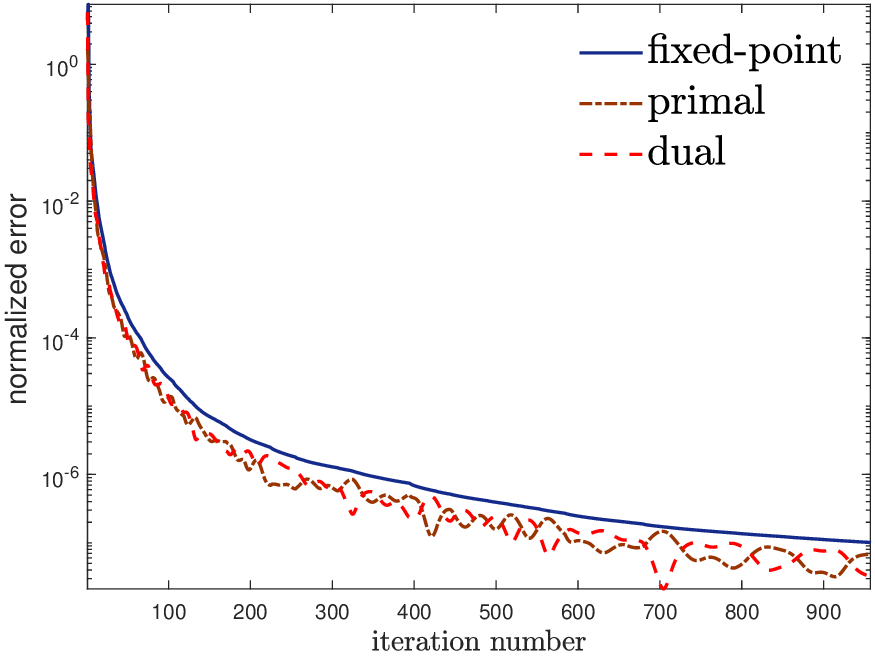}\caption{random initialization $ \mathcal{N}(0, 1) $.}
	\end{subfigure}
	
	\caption{Practical  use  test (adaptive,  $\rho_{0} =  1$).}
\end{figure}

\begin{remark}
	Compared the above to the previous zero initialization case,  we see  that a  warm start does improve the performance. Interestingly, random initialization exhibits a similar performance compared to Fig. \ref{fig_lp2}.
\end{remark}

\subsection{Quadratic  programming}
Consider the following quadratic  program:
\begin{align}
\text{minimize}     \,\,  		\frac{1}{2}\bm{x}^T\bm{P}\bm{x} + \bm{q}^T\bm{x} + r \quad
\text{subject to}  \,\, 		\bm{a} \leq \bm{x} \leq \bm{b},  
\end{align}
with variable $r \in \mathbb{R}, \bm{x}, \bm{q}\in \mathbb{R}^{100}, \bm{P} \in \mathbb{S}^{100}_{++} $  (positive  definite   matrix). 
The elements of $ \bm{P}  $ are first generated from a uniform distribution on interval $ [0, 1] $. Then, its eigenvalues are increased by a uniformly distributed random number from $ [0, 1] $. 
All other data is generated from $ \mathcal{N}(0,1) $.
\begin{figure}[H]
	\centering
	\begin{subfigure}[t]{0.42\textwidth}
		\includegraphics[scale=0.42]{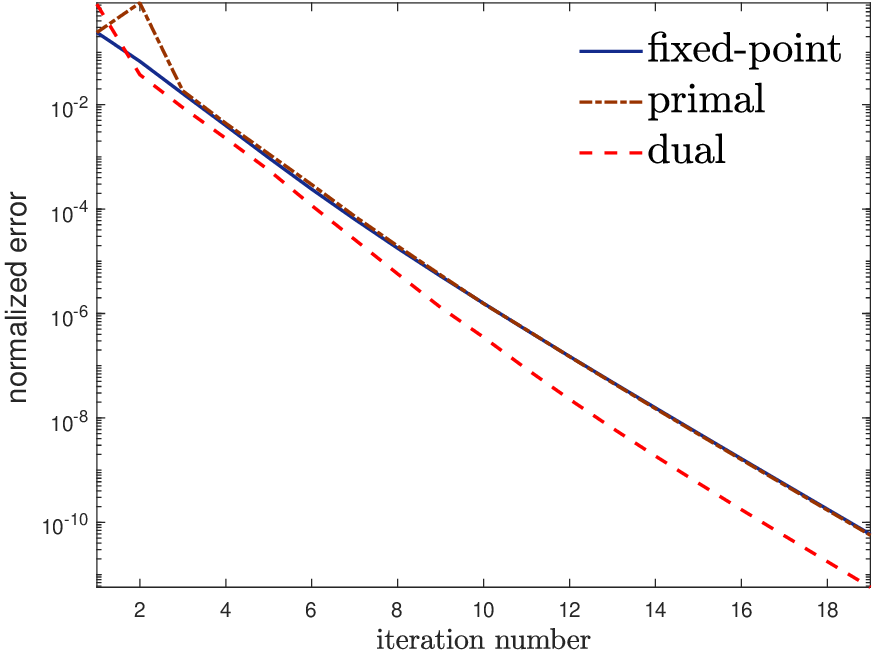}\caption{Step-size: optimal  choice $\rho^\star \approx \pm 1.34$.}
	\end{subfigure}
	\begin{subfigure}[t]{0.42\textwidth}
		\includegraphics[scale=0.42]{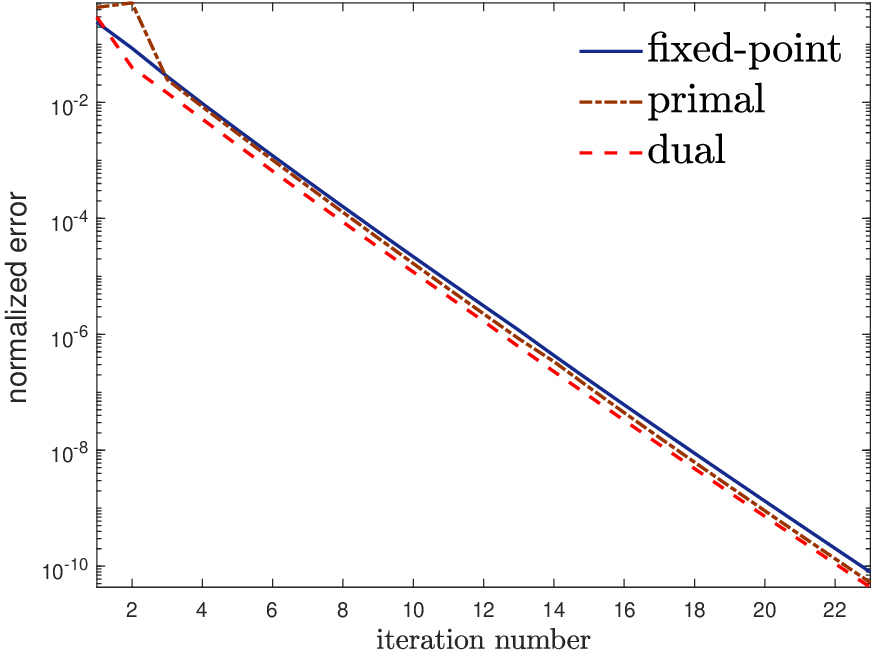}\caption{Step-size: unit length.}
	\end{subfigure}
	\caption{Theoretical results test (zero initialization).}
	\begin{subfigure}[t]{0.42\textwidth}
		\includegraphics[scale=0.42]{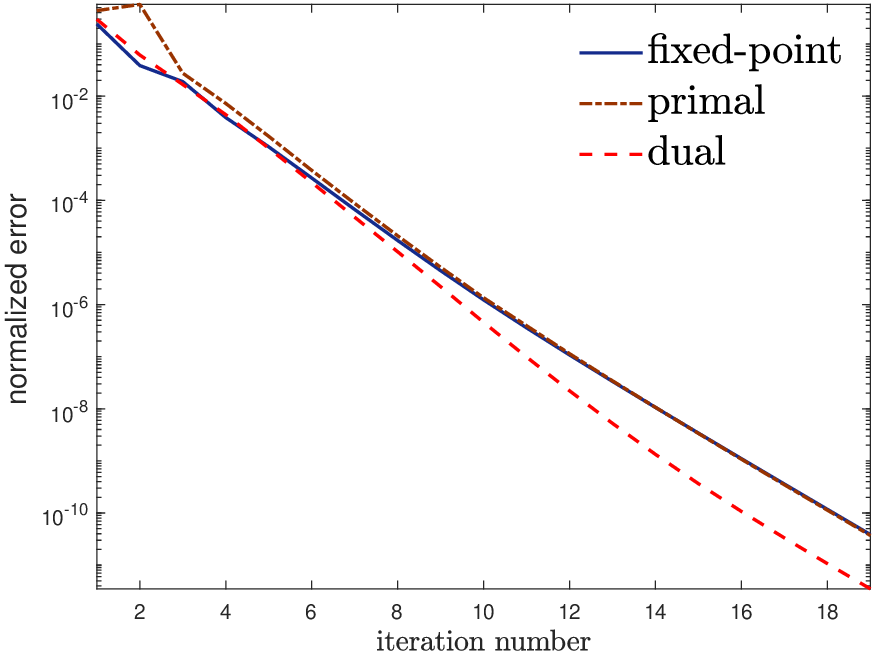}\caption{Estimated step-size choice (adaptive).}
	\end{subfigure}
	\begin{subfigure}[t]{0.42\textwidth}
		\includegraphics[scale=0.42]{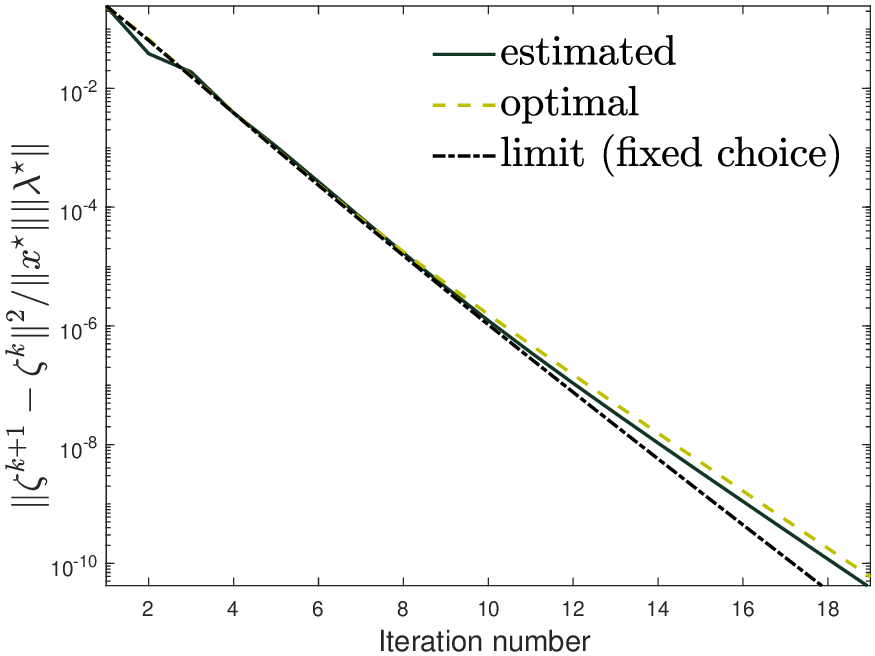}\caption{Fixed-point convergence comparison.}
	\end{subfigure}
	\caption{Practical  use  test (zero initialization).}
	\begin{subfigure}[t]{0.42\textwidth}
		\includegraphics[scale=0.42]{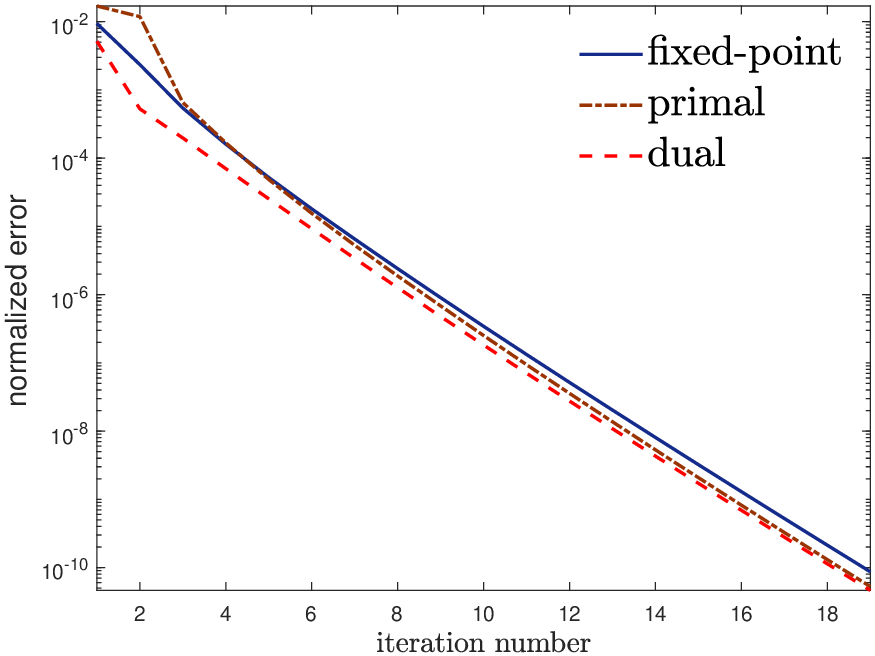}\caption{warm-start: $  \epsilon_{\text{err}1},\epsilon_{\text{err}2} \sim  \mathcal{N}(0, 0.1) $.}
	\end{subfigure}
	\begin{subfigure}[t]{0.42\textwidth}
		\includegraphics[scale=0.42]{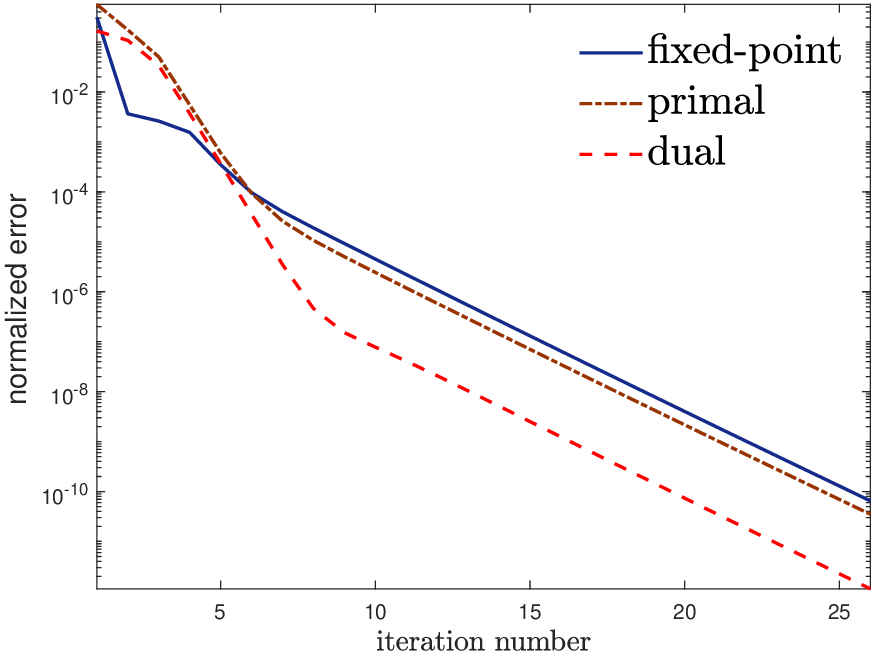}\caption{primal warm-start:  $  \epsilon_{\text{err}1} \sim  \mathcal{N}(0, 0.1) $.}
	\end{subfigure}	
	\begin{subfigure}[t]{0.42\textwidth}
		\includegraphics[scale=0.42]{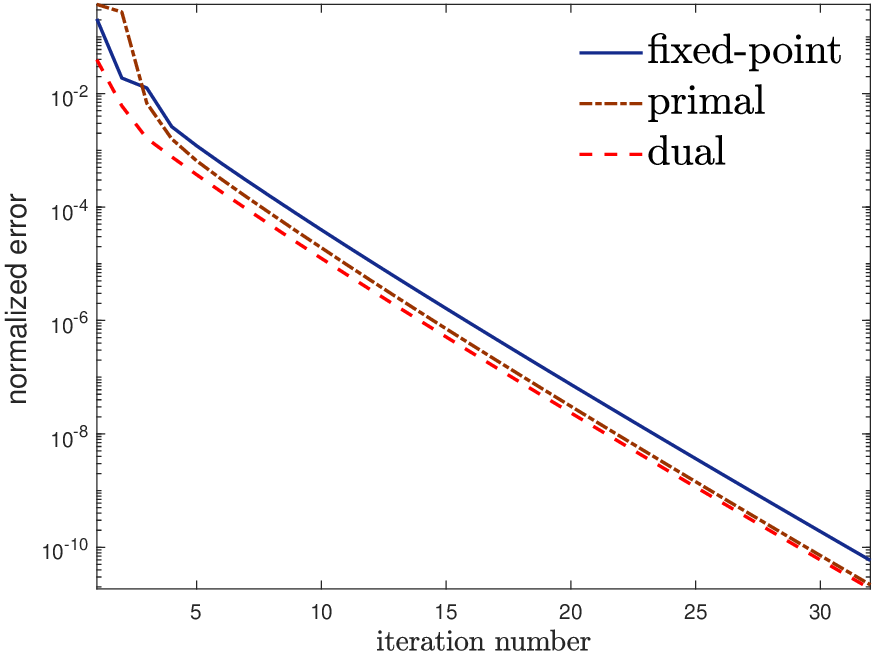}\caption{dual warm-start: $  \epsilon_{\text{err}2} \sim  \mathcal{N}(0, 0.1) $.}
	\end{subfigure}
	\begin{subfigure}[t]{0.42\textwidth}
		\includegraphics[scale=0.42]{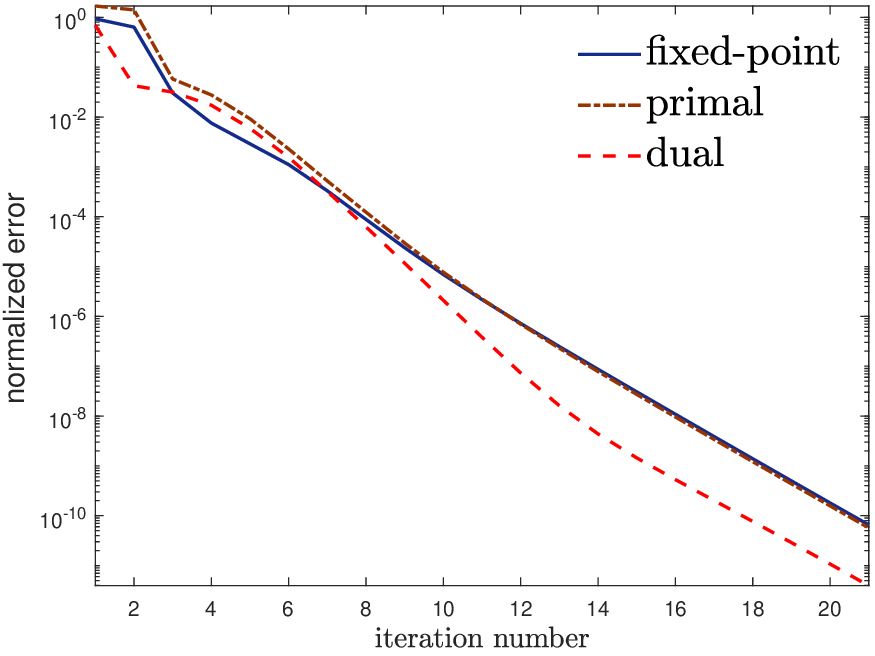}\caption{random initialization $ \mathcal{N}(0, 1) $.}
	\end{subfigure}
	\caption{Practical  use  test (non-zero  initialization,  $\rho_{0} =  1$).}
\end{figure}

\begin{remark}
	Under this specific setting  for QP, we note that different settings perform similarly, implying a well-conditioned nature.
\end{remark}

\subsection{Total  variation}
Consider the following total  variation  problem:
\begin{equation}
\text{minimize}    \,\, 		 \frac{1}{2}\Vert\bm{x-b}\Vert^2 +  \alpha \sum_i |x_{i+1}-x_i|
\end{equation}
with variable $ \bm{x},  \bm{b} \in \mathbb{R}^{100}$.  $  \bm{x} $ is first set to be a ones vector and then randomly scaled on randomly selected entries. $ \bm{b} $ is set to equal to $ \bm{x} + \bm{\epsilon}$, where $ \bm{\epsilon}  $ is an additive Gaussian noise from  $ \mathcal{N}(0,1) $.
Regularization parameter  is set $\alpha =  5$. Moreover,  term $ \sum_i |x_{i+1}-x_i| $  can  be compactly  written into $ \bm{Dx} $, where  we use the popular notation $ \bm{D} $  to denote a difference  matrix, corresponding to  our operator $\mathcal{A}$.

\subsubsection{zero  initialization}
Set initialization $\bm{x}^0 =  \bm{0}$, $\bm{\lambda}^0 =  \bm{0}$.
\begin{figure}[H]
	\centering
	\begin{subfigure}[t]{0.42\textwidth}
		\includegraphics[scale=0.42]{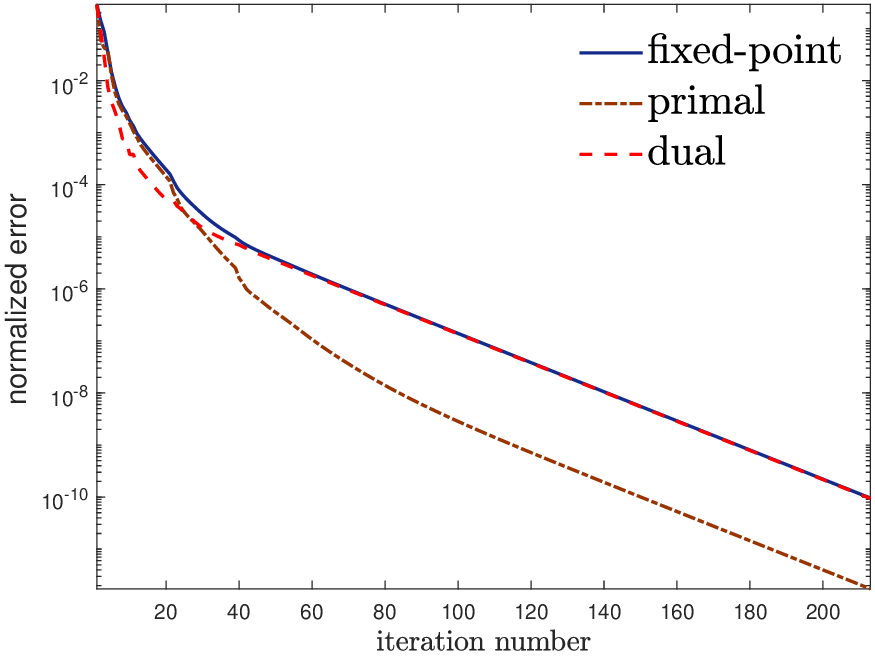}\caption{Step-size: optimal  choice $\rho^\star \approx \pm 2.1$.}
	\end{subfigure}
	\begin{subfigure}[t]{0.42\textwidth}
		\includegraphics[scale=0.42]{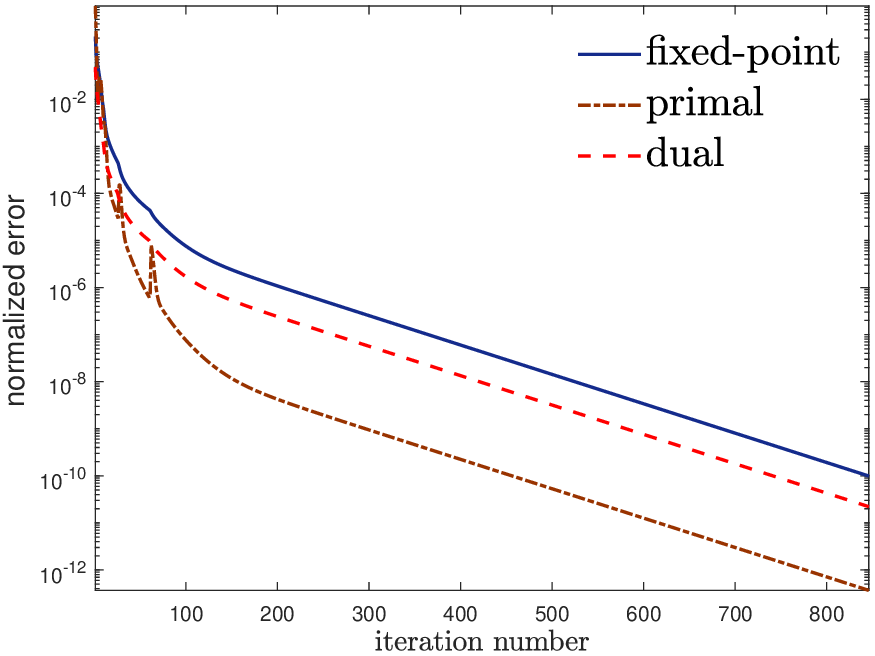}\caption{Step-size: unit length.}
	\end{subfigure}
	\caption{Theoretical results test (zero initialization).}
	\begin{subfigure}[t]{0.42\textwidth}
		\includegraphics[scale=0.42]{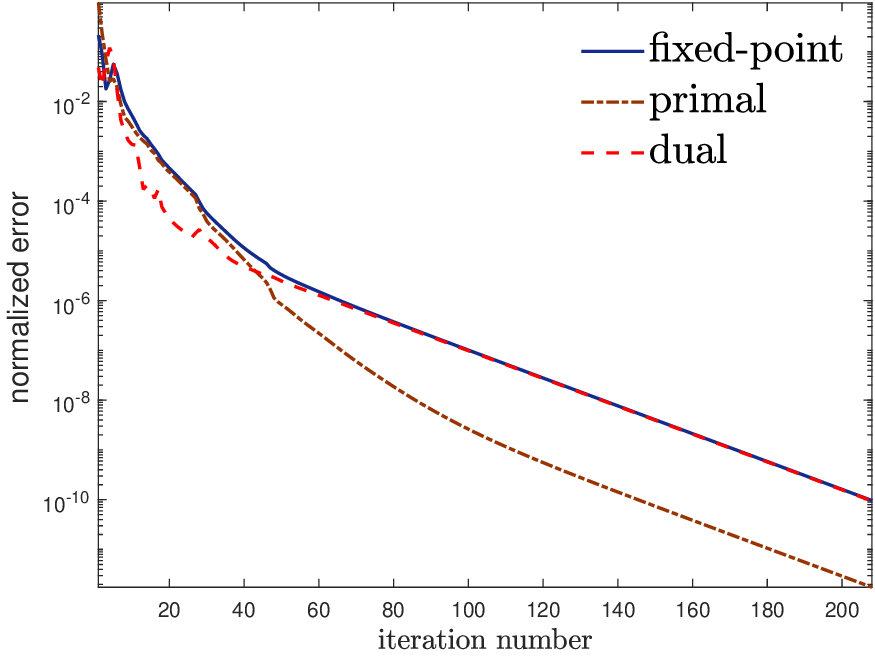}\caption{Estimated step-size choice (adaptive).}\label{fig_tv}
	\end{subfigure}
	\begin{subfigure}[t]{0.42\textwidth}
		\includegraphics[scale=0.42]{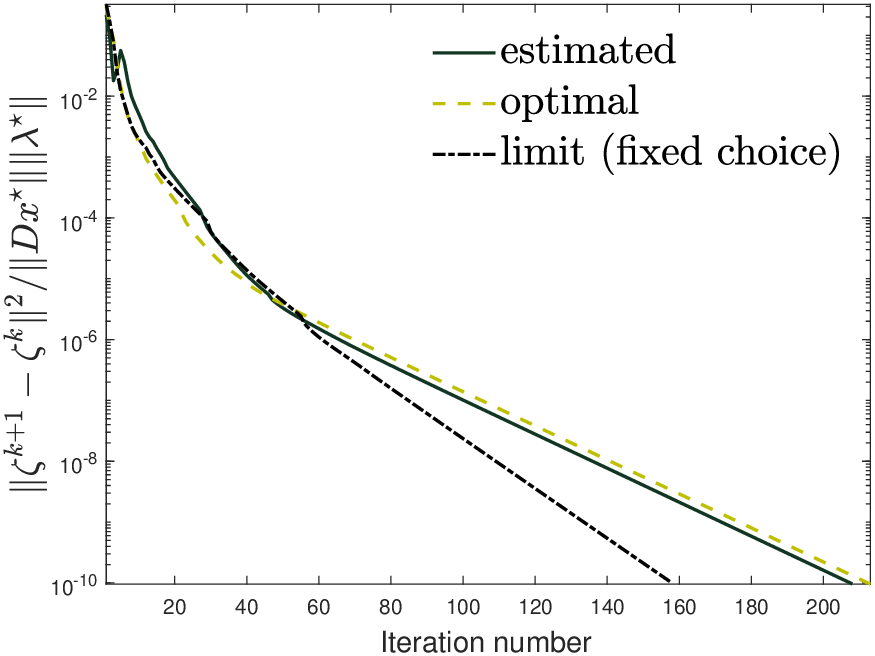}\caption{Fixed-point convergence comparison.}\label{fig_tv2}
	\end{subfigure}
	\caption{Practical  use  test (zero initialization).} 
\end{figure}

\begin{remark}
	While the optimal choice is quite close  to the unit length, it  makes a  big difference for this  application,  yielding a $ 4  \times$  acceleration.
	This is a big improvement due  to the data is relatively  well-conditioned. 
	Indeed,  from Fig. \ref{fig_tv2},  we see that the limit (for fixed step-size) is with iteration 158, and we  achieved  213 without  any  tailored structure exploited.  
\end{remark}

\subsubsection{non-zero  initialization}
Set non-zero $\bm{\zeta}^0 \neq  \bm{0}$ according to the 4 strategies at the beginning of the appendix.
\begin{figure}[H]
	\centering
	\begin{subfigure}[t]{0.42\textwidth}
		\includegraphics[scale=0.42]{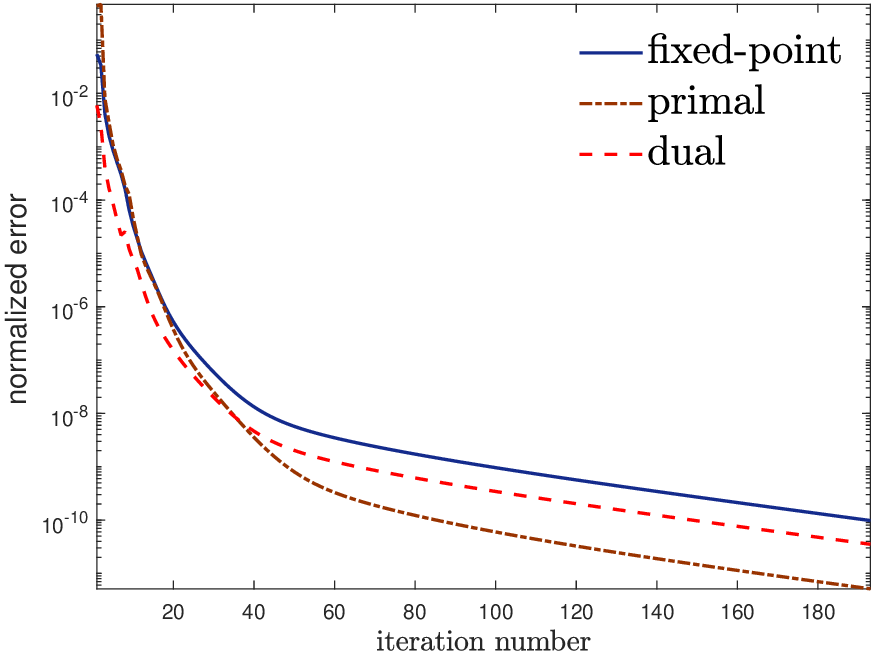}\caption{warm-start: $  \epsilon_{\text{err}1},\epsilon_{\text{err}2} \sim  \mathcal{N}(0, 0.1) $.}
	\end{subfigure}
	\begin{subfigure}[t]{0.42\textwidth}
		\includegraphics[scale=0.42]{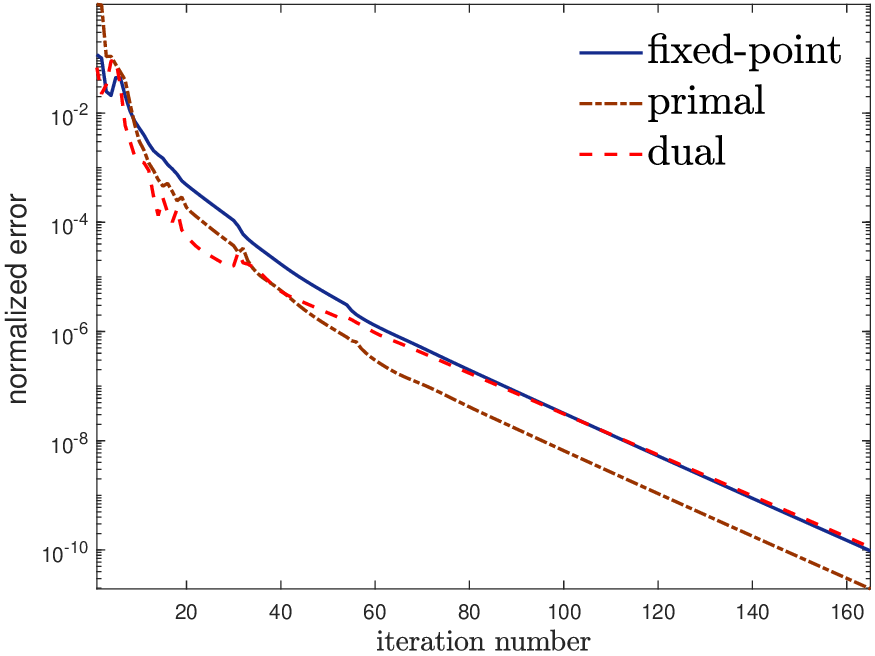}\caption{primal warm-start:  $  \epsilon_{\text{err}1} \sim  \mathcal{N}(0, 0.1) $.}
	\end{subfigure}	
	\begin{subfigure}[t]{0.42\textwidth}
		\includegraphics[scale=0.42]{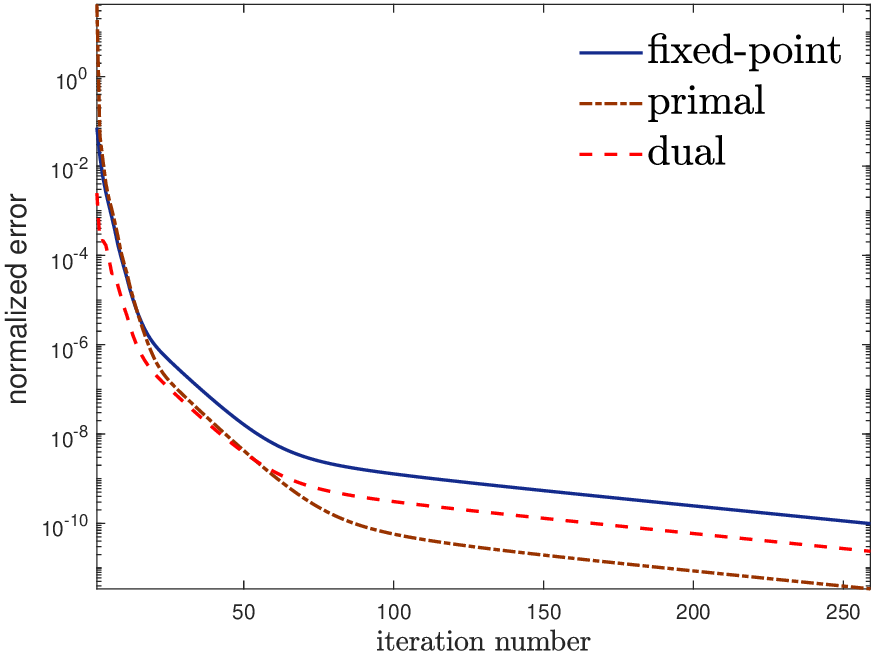}\caption{dual warm-start: $  \epsilon_{\text{err}2} \sim  \mathcal{N}(0, 0.1) $.}
	\end{subfigure}
	\begin{subfigure}[t]{0.42\textwidth}
		\includegraphics[scale=0.42]{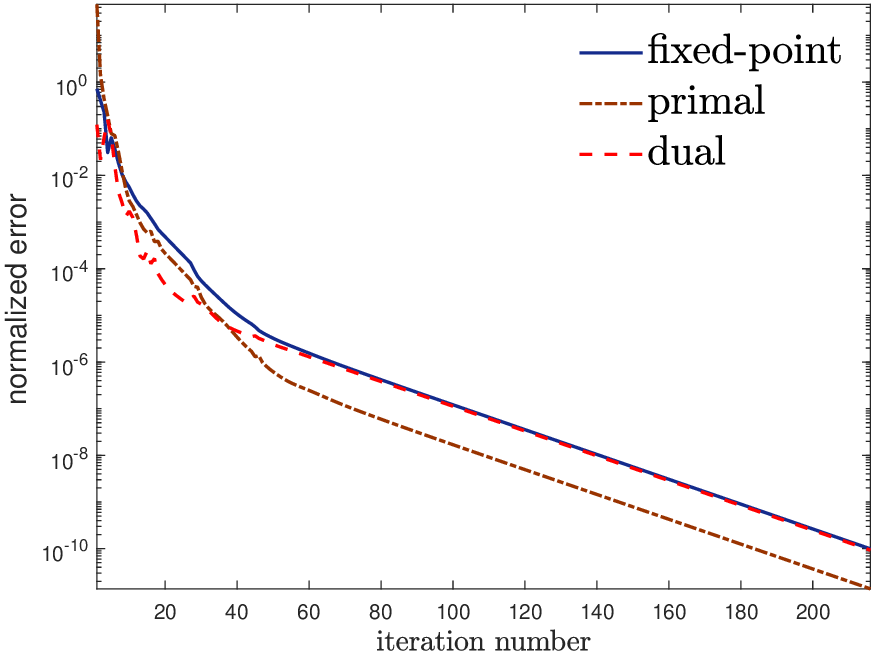}\caption{random initialization $ \mathcal{N}(0, 1) $.}
	\end{subfigure}
	
	\caption{Practical  use  test (non-zero  initialization,  $\rho_{0} =  1$).}
\end{figure}

\begin{remark}
	Compared the above to the previous zero initialization case in Fig. \ref{fig_tv},  we see  that the $ \mathcal{N}(0, 0.1) $-type warm start does slightly improve the performance. Moreover, under the current setting,  we observe that  the primal warm-start strategy performs the best.
\end{remark}

\end{document}